\documentclass[11pt]{article}

\usepackage[flushleft]{threeparttable}
\usepackage{multirow}
\usepackage{shortcuts}
\usepackage{stmaryrd}
\usepackage{lipsum}% Table of contents
\usepackage{graphicx}
\usepackage{subfig}
\usepackage{setspace}
\usepackage{mathtools}
\usepackage{pifont}
\usepackage{indentfirst}
\newcommand{\cmark}{\ding{51}}%
\newcommand{\xmark}{\ding{55}}%
\usepackage{tikz}

\graphicspath{{./Figures/}} 
\usepackage{hyperref}
\usepackage{titlefoot}
\usepackage{tikz}
\usetikzlibrary{patterns}
\hypersetup{
    colorlinks,
    citecolor=black,
    filecolor=black,
    linkcolor=black,
    urlcolor=black
}

\newcommand{\nux}{{\nu^x}}

\newcommand{\cc}{c}
\newcommand{\quand}{\quad \text{and} \quad}

\newcommand{\wh}[1]{{\widehat{#1}}}

\newcommand{\dx}{{\Delta x}}
\newcommand{\dt}{{\Delta t}}

\def\!{\mskip-\thinmuskip}
\newcommand{\hb}[1]{{\{\!\!\{#1\}\!\!\}_\gamma}}
\newcommand{\hbn}[1]{{\left\{\!\!\left\{#1\right\}\!\!\right\}}}
\newcommand{\jp}[1]{{[#1]}}

\newcommand{\sgnb}{{\mathrm{sgn}(\jp{b})}}

\newcommand{\tU}{\widetilde{U}}

\newcommand{\tOmega}{\widetilde{\Omega}}
\newcommand{\mm}{m}
\newcommand{\dd}{\mathrm{d}}
\usepackage{amsthm}

\theoremstyle{definition}
\newtheorem{ex}{Example}[section]
\numberwithin{figure}{section}

\title{On a numerical artifact of solving shallow water equations with a discontinuous bottom: Analysis and a nontransonic fix}
\author{
	Zheng Sun\footnote{Department of Mathematics, The University of Alabama,
		Tuscaloosa, AL 35487, USA. E-mail: zsun30@ua.edu. The work of this author is partially supported by the NSF grant DMS-2208391. }\and Yulong Xing\footnote{Department of Mathematics, The Ohio State University,
		Columbus, OH 43210, USA. E-mail: xing.205@osu.edu. The work of this author is partially supported by the NSF grants DMS-1753581 and DMS-2309590.}} 
\date{}
\newcommand{\nm}[1]{\| #1 \|}

\newcommand{\halfcmark}{\cmark\kern-1.3ex\raisebox{.3ex}{\xmark}}

\allowdisplaybreaks 

\begin{document}

\unmarkedfntext{2010 \emph{Mathematics Subject Classification.} Primary 65M06, 65M12, 65M22, 35L65, 35L67.}
\unmarkedfntext{\emph{Key words and phrases.} Shallow water equations. Nonconservative hyperbolic systems. First-order numerical schemes. Lax--Wendroff theorem. Discontinuous bottom topography. Riemann problems.} 
\maketitle 

\textbf{Abstract.} In this paper, we study a numerical artifact of solving the nonlinear shallow water equations with a discontinuous bottom topography. For various first-order schemes, the numerical solution of the momentum will form a spurious spike at the discontinuous points of the bottom, which should not exist in the exact solution. The height of the spike cannot be reduced even after the mesh is refined. For subsonic problems, this numerical artifact may cause the wrong convergence to a function far away from the exact solution. To explain the formation of the spurious spike, we perform a convergence analysis by proving a Lax--Wendroff type theorem. It is shown that the spurious spike is caused by the numerical viscosity in the computation of the water height at the discontinuous bottom. The height of the spike is proportional to the magnitude of the viscosity constant in the Lax--Friedrichs flux. Motivated by this conclusion, we propose a modified scheme by adopting the central flux at the bottom discontinuity in the equation of mass conservation, and show that this numerical artifact can be removed in many cases. For various numerical tests with nontransonic Riemann solutions, we observe that the modified scheme is able to retrieve the correct convergence.  

\section{Introduction}\label{sec:introduction}
\setcounter{equation}{0}
\setcounter{figure}{0}
\setcounter{table}{0} 

The nonlinear shallow water equations (SWEs) are used to model the free surface flow in rivers and coastal areas for which the horizontal length scale is much greater than the vertical length scale. They have wide applications in atmospheric and oceanic sciences and hydraulic engineering, such as  prediction of tsunami and storm surges, simulation of dam break and flooding, etc. In this paper, we study the numerical solutions of the one-dimensional SWEs with a discontinuous bottom topography
\begin{equation}\label{eq-swe}
{U}_t + {F}({U})_x = {S}({U})b_x, \quad U(x,0) = U_0,
\end{equation}
where
\begin{equation}\label{eq-swe-def}
{U} = \left(
\begin{array}{c}
h\\
\mm
\end{array}\right), \quad {F}({U}) = \left(
\begin{array}{c}
\mm\\
\frac{\mm^2}{h} + \hf gh^2 
\end{array}\right),\quand {S}({U}) = \left(
\begin{array}{c}
0\\
-gh
\end{array}\right).
\end{equation}
Here $h$ is the water height, $m = hu$ is the momentum, $u$ is the velocity, $g$ is the (known) gravity acceleration constant, and $b$ is the (known) bottom topography function. To address the main issue, we avoid the discussion of the dry bed problems by assuming $h\geq \cc_0>0$ to be uniformly positive. Furthermore, we assume $b$ is smooth except for a single jump discontinuity at $x = 0$.  

In the case that the bottom topography is flat, the source term vanishes, and the SWEs \eqref{eq-swe} become a homogeneous system of hyperbolic conservation laws. This system and  associated numerical methods have been well studied in the literature \cite{serre1999systems,leveque1992numerical,hesthaven2017numerical}. The complication arises as the bottom becomes nonflat, especially when $b$ is discontinuous. In this case, $b$ can be treated as an additional unknown, and \eqref{eq-swe} is then augmented as a nonconservative hyperbolic system %\cite{lefloch1989shock}
\begin{equation}\label{eq-ncsv}
	\tU_t + A(\tU) \tU_x = 0,
\quad \text{ with }\quad 
	\tU = \left(
	\begin{array}{c}
		U\\
		b
	\end{array}\right)
	%= \left(
	%\begin{array}{c}
	%h\\
	%q\\
	%b
	%\end{array}\right), 
	\quand A(\tU) = 
	\left(
	\begin{array}{cc}
		\frac{\partial F}{\partial U}&-S\\
		0&0
	\end{array}\right).
	%=
	%\left(
	%\begin{array}{ccc}
	%0&1&0\\
	%-u^2+gh&2u&gh\\
	%0&0&0
	%\end{array}\right)
\end{equation}
Its weak solution can be defined through the theory developed by Dal Maso, LeFloch and Murat \cite{maso1995definition} for the nonconservative product, which requires \eqref{eq-ncsv} to be satisfied in the sense of Borel measures. This definition relies on a prescribed path connecting the two states in the phase plane and different choices of the path may lead to different weak solutions. For the SWEs, the choice of the path is related to the definition of the hydrostatic pressure at the bottom discontinuity \cite{cozzolino2011numerical}. Besides the weak solution, people have also studied the exact solutions to the Riemann problems for the SWEs. Among these works, there have been some controversial discussions on how to characterize the relations connecting the flow variables across the bottom discontinuity. Two approaches have been pursued in the literature, and they may lead to different exact solutions for the same Riemann problem. One approach is based on the mass and energy conservation, and the derived Riemann solution preserves the Riemann invariants \cite{alcrudo2001exact, lefloch2007riemann,lefloch2011godunov,han2014exact,aleksyuk2019uniqueness}. The other approach is based on the mass and momentum conservation, and the derived solution satisfies the generalized Rankine--Hugoniot condition \cite{bernetti2008exact,rosatti2010riemann}. We briefly discuss the difference between these two approaches in Section \ref{sec:riemann} and refer to \cite{rosatti2010riemann} for further details. In this paper, the exact reference solution is generated through the second approach by the exact Riemann solvers in \cite{bernetti2008exact,rosatti2010riemann} as will be explained in Section \ref{sec:riemann}. 
%Specifically for the SWEs, there have also been various studies of the exact solutions to Riemann problems by considering the intersections of Hugoniot loci. See \cite{alcrudo2001exact, lefloch2007riemann,bernetti2008exact,lefloch2011godunov,han2014exact,aleksyuk2019uniqueness} and references therein. One difficulty is that the system \eqref{eq-ncsv} may be nonstrictly hyperbolic at sonic states, which leads to resonant waves and nonunique solutions \cite{han2014exact}. Appropriate conditions may be needed to single out the physically relevant solution \cite{aleksyuk2019uniqueness}.

As for the numerical discretization of the SWEs \eqref{eq-swe} with a nonflat bottom, different computational methods have been developed in the past decades. See, for example, \cite{bermudez1994upwind,kurganov2002central,audusse2004fast,XS2005,RB2009,XZS2010,fjordholm2011well,berthon2016fully,HXX2022}. We refer to \cite{xing2017numerical} and references therein for further related works. However, the convergence of the numerical methods to the exact solution is usually not guaranteed in general settings. To address the dependence of the weak solution on the prescribed path, Par\'es introduced the so-called path-conservative schemes in \cite{pares2006numerical} to preserve the formal consistency. See also \cite{castro2017well} and references therein. However, a detailed numerical investigation in \cite{abgrall2010comment} shows that the path-conservative schemes may not guarantee the correct convergence. The issue can be further explained by the work of \cite{hou1994nonconservative} and \cite{castro2008many}. Especially in \cite{castro2008many}, a convergence analysis shows that after mesh refinement the limit of the numerical solution does not satisfy \eqref{eq-ncsv}, but admits an inhomogeneous hyperbolic system containing a Borel measure source term. This source term vanishes only when certain strong convergence assumptions can be made, which unfortunately may not hold in general. Despite these results, a convergence theorem was proved by Mu\~{n}oz-Ruiz and Par\'es  in \cite{munoz2011convergence}. In that paper, they specifically analyzed a balance law in the form of \eqref{eq-swe} and showed that when $b\in W^{1,1}$ the path-conservative scheme converges to a weak solution of \eqref{eq-swe} under classical assumptions. However, by the embedding theorem, $b \in W^{1,1}$ implies that $b$ is absolutely continuous. Although the correct convergence can be assured with smooth bottom topography, the SWEs with an abrupt riverbed (when $b$ contains a jump discontinuity) remain uncovered. 

In this paper, we analyze numerical solutions to the SWEs with a discontinuous bottom. We consider several first-order schemes, which include the non-well-balanced Lax--Friedrichs scheme (LxF scheme), the well-balanced LxF scheme (wbLxF scheme), the well-balanced scheme with hydrostatic reconstruction by Audusse et al. \cite{audusse2004fast} (HR scheme), and the first-order version of the well-balanced scheme with flux and source modification by Xing and Shu \cite{xing2006high} (XS scheme). We note that all these numerical schemes suffer a similar numerical artifact:
\begin{itemize}
	\item \textit{The exact momentum of the SWEs should be continuous across the bottom discontinuity \cite[the first equation in (17)]{bernetti2008exact}. But the numerical solutions of all the above-mentioned first-order schemes will form a spurious spike at the bottom discontinuity. Furthermore, the height of the spike does not decrease as one refines the mesh. }
\end{itemize}
Accompanying this numerical artifact, the numerical solutions of both $h$ and $m$ near the discontinuity could be far from the exact solutions for many tested subsonic problems, and they may converge to a wrong solution after the mesh refinement. We refer to Sections \ref{sec:artifact} and \ref{sec:num} for detailed numerical tests. 

To explain the cause of the numerical artifact and the wrong convergence, we rewrite these different numerical schemes into a unified form \eqref{eq:pcs} and establish a Lax--Wendroff type theorem on two sides of the riverbed separated by the bottom discontinuity. With this Lax--Wendroff type theorem, we can further deduce the following result on the limit of the numerical solution: 
\begin{itemize}
	\item \textit{The height of the spurious spike is proportional to the numerical viscosity in the equation of mass conservation at the bottom discontinuity, and the jump of $h$ for the non-well-balanced scheme or the jump of $h+b$ for the well-balanced schemes. }
\end{itemize}  
This observation motivates us to look into a central-Lax--Friedrichs scheme (cLxF), which applies the central flux in the equation of mass conservation at the bottom discontinuity (with zero numerical dissipation) and the LxF flux everywhere else. A similar modification can also be applied to well-balanced methods. Although we cannot prove that the scheme will converge to a correct weak solution, the following partial results can be proven:
\begin{enumerate}
	\item The cLxF scheme will not form a one-sided spurious spike in the numerical momentum at the bottom discontinuity.
	\item Under certain assumptions, the numerical water height of the cLxF scheme does not have transition points at the bottom discontinuity. 
\end{enumerate}
Numerically, we observe the cLxF scheme converges to the correct solution for all the nontransonic Riemann problems we tested. It is especially notable that at the same time the classical LxF scheme will converge to a totally wrong solution for subsonic tests. But for transonic (or resonant) problems, the cLxF scheme may converge to a wrong solution possibly with different wave patterns. This issue may be related to the incapability of capturing the entropy solution, and will be left for future investigation. 

We remark that our analysis is closely related to the work \cite{castro2008many} on the convergence analysis of the path-conservative schemes. In \cite{castro2008many}, the authors also attribute the wrong convergence to the numerical viscosity. They considered a generic nonconservative hyperbolic problem and derived the equivalent equation satisfied by the LxF scheme. It is shown that the vanishing viscosity in the equivalent equation may differ from  that associated with the prescribed path. In our work, we give a further characterization of this issue for \eqref{eq-swe} from a different perspective and show that the wrong convergence will accompany the formation of the spurious spike, whose height is directly determined by the strength of the numerical viscosity. %Furthermore, due to the close connection of this work and \cite{castro2008many}, we remark that further improvement of the cLxF scheme in the current framework can be difficult. As is pointed out in \cite{castro2017well}, ``even the Godunov method based on the right weak solutions of the Riemann problems fails, in general, to converge to the right solution: this is due to the numerical viscosity introduced in the averaging step." One may have to seek other numerical strategies such as front tracking schemes, random choice schemes, and controlled dissipation schemes. See \cite{lefloch2014numerical} and references therein. 

In addition to the aforementioned numerical artifact, we would like to point out another factor that will effect the solution limit of a numerical scheme. As has been mentioned and will be further explained in Section \ref{sec:ws}, the choice of the path in the definition of the weak solution is closely related to the definition of the hydrostatic pressure in the SWEs model at the bottom discontinuity. Typically, the pressure is assumed to be proportional to the water depth. However, due to the presence of the bottom step, the corresponding water depth is double-valued, and one needs to introduce a parameter $\gamma$ (see \eqref{eq-hb} for details) to indicate how the water depth is computed when evaluating the pressure along the step. We also refer to \cite{cozzolino2011numerical} for further details. Different existing numerical methods may correspond to different hydrostatic pressure terms. For example, the hydrostatic reconstruction scheme \cite{audusse2004fast} comes with a hydrostatic pressure with the parameter $\gamma = \sgnb$, corresponding to the pressure on the lower side of bottom step. But the XS scheme in \cite{xing2006high} comes with a hydrostatic pressure with the parameter $\gamma = 0$, corresponding to the arithmetic average of the pressure on both sides of the bottom step. As different values of $\gamma$ correspond to different weak solutions, one can expect that the two schemes will also converge to different solutions. This fact indicates that some numerical schemes, although not originally designed as a path-conservative scheme, are indeed inherently associated with certain choices of the path. Due to the aforementioned numerical artifact, preserving the path along does not guarantee correct convergence. But for schemes corresponding to different paths, even if they converge as expected, their limit solutions would be different.

The rest of the paper is organized as follows. In Section \ref{sec:artifact}, we use a simple example to illustrate the numerical artifact and the wrong convergence of the first-order schemes when solving the SWEs with a discontinuous bottom. In Section \ref{sec:ws}, we revisit the definition of the weak solution to \eqref{eq-swe} and explain the connection between the path and the hydrostatic pressure. Two different definitions of the Riemann solutions will also be briefly explained. In Section \ref{sec:schemes}, we rewrite several first-order numerical schemes into a unified form. The main theoretical results of the paper are given in Section \ref{sec:ca}. We prove a Lax--Wendroff type theorem on the convergence of the solution limit (Subsection \ref{sec:glw}), explain the formation of the spurious spike (Subsection \ref{sec:pwlimit}), and present the cLxF scheme along with its theoretical properties (Subsection \ref{sec:cLxF}). Detailed numerical tests with the cLxF versus LxF schemes are given in Section \ref{sec:num}. Conclusions are given in Section \ref{sec:conclusion}. 

\section{Numerical artifact and wrong convergence}\label{sec:artifact}
\setcounter{equation}{0}

We consider a dam-break problem \eqref{eq-swe} with the following settings: 
\begin{equation}\label{eq:db}
U^0 = \begin{pmatrix}
h^0\\ m^0
\end{pmatrix}, \quad 
h^0 = \left\{
\begin{array}{cc}
1& x<0\\
0.1& x\geq 0
\end{array}\right., \quad  
m^0 = 0, \quad b = \left\{\begin{array}{cc}
0& x<0\\
0.7& x\geq 0
\end{array}\right., \quand g = 9.81.
\end{equation}
Its exact solution obtained through the procedure in \cite[Section 3.6]{bernetti2008exact} admits the states
\begin{equation} \label{sec2:exact}
U_L = \begin{pmatrix}
1\\
0
\end{pmatrix}\xrightarrow{\text{1-rarefaction}}\begin{pmatrix}
0.9458\\
0.1629
\end{pmatrix}\xrightarrow{\text{0-wave}}
\begin{pmatrix}
0.1964\\
0.1629
\end{pmatrix}\xrightarrow{\text{2-shock}}
\begin{pmatrix}
0.1\\0
\end{pmatrix} =U_R.
\end{equation}
Here we address that we follow the derivation in  \cite{bernetti2008exact,rosatti2010riemann} for the exact Riemann solution, which is different from the solution in \cite{alcrudo2001exact, lefloch2007riemann,lefloch2011godunov,han2014exact,aleksyuk2019uniqueness}. The exact solution \eqref{sec2:exact} admits the generalized Rankine--Hugoniot condition, and its Riemann invariants are not constant across the bottom discontinuity. We refer to Section \ref{sec:riemann} for further details and explanations. 

In our numerical simulation, we set the computational domain as $[-5,5]$ and compute to $T = 1$. The CFL number in this test is taken as $0.5$. We apply $N = 100$, $200$, and $25600$ grid points in the spatial discretization. We consider $N = 25600$ to be an extremely refined mesh such that the corresponding numerical solution can be viewed as the solution limit of the schemes. 

\begin{figure}[!h]
	\centering
	\subfloat[LxF, $h$.]{%
		\includegraphics[width=0.25\textwidth,trim={2cm 0 1cm 0}]{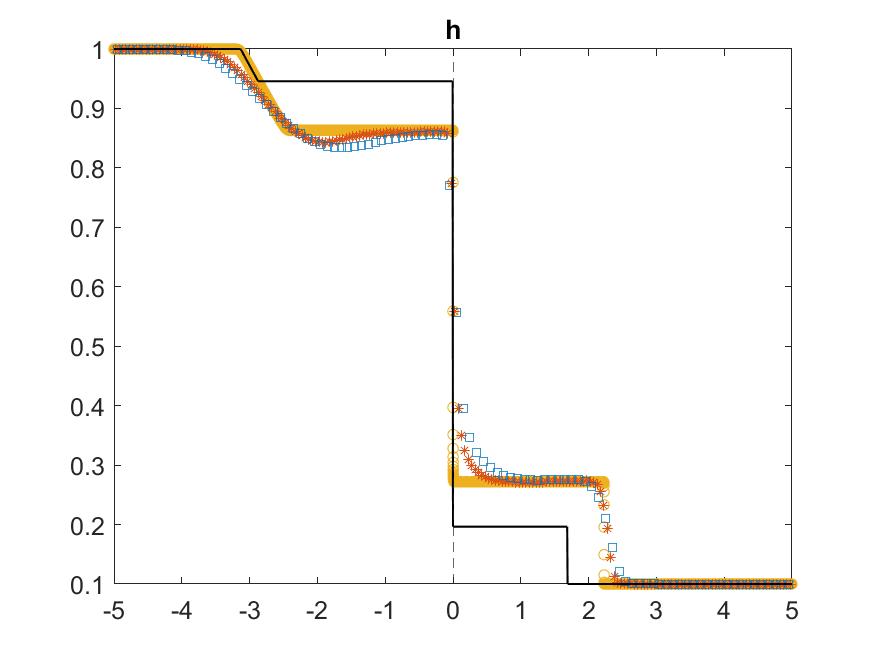}
		\label{fig:db-h-LxF}}
	\subfloat[LxF, $m= hu$.]{%
		\includegraphics[width=0.25\textwidth,trim={2cm 0 1cm 0}]{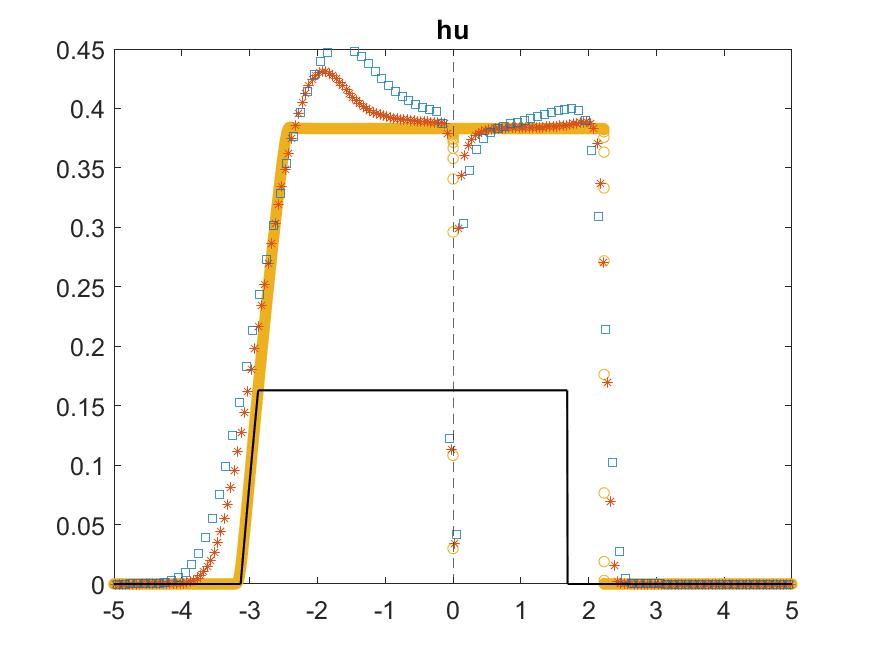}
		\label{fig:db-m-LxF}}
	\subfloat[HLLC, $h$.]{%
		\includegraphics[width=0.25\textwidth,trim={2cm 0 1cm 0}]{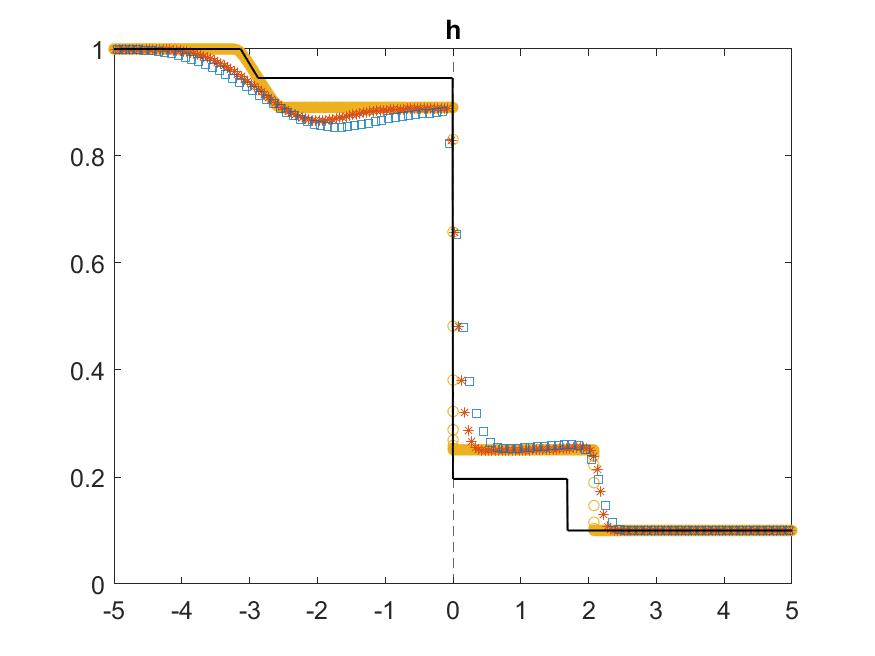}
		\label{fig:db-h-HLLC}}
	\subfloat[HLLC, $m= hu$.]{%
		\includegraphics[width=0.25\textwidth,trim={2cm 0 1cm 0}]{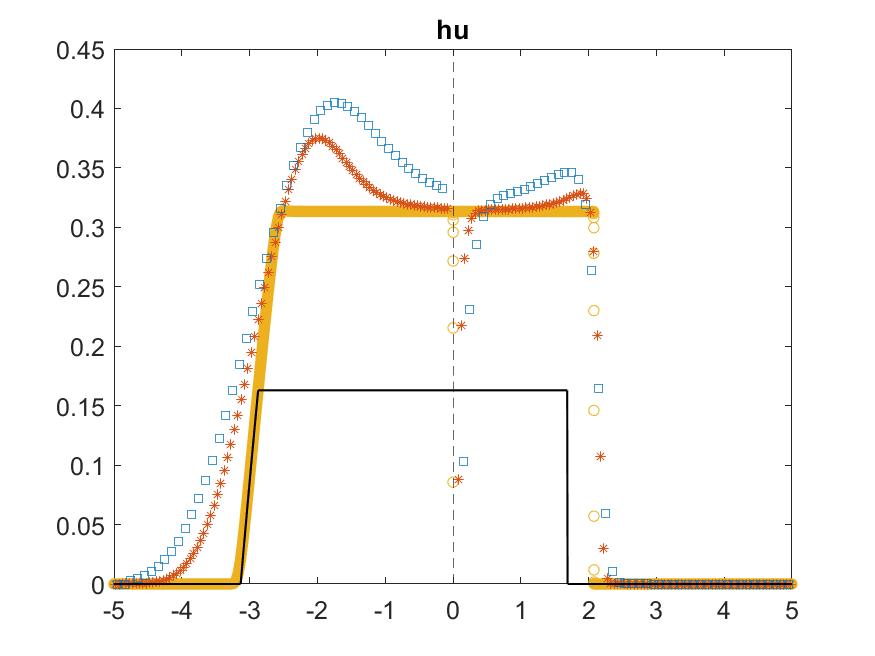}
		\label{fig:db-m-HLLC}}
	\caption{Solutions of the LxF and HLLC schemes to the dam-break problem \eqref{eq:db}. Black solid line: exact solution; blue squares: $N = 100$; red stars: $N = 200$, yellow circles: $N = 25600$.}\label{fig:db-nwb}
\end{figure}

We first solve the problem with the (local) LxF scheme (defined in \eqref{eq:LxF}), which can also be viewed as a slight variation of the path-conservative scheme in \cite{castro2008many}. Numerical solutions are given in Figure \ref{fig:db-nwb}. The exact solution \eqref{sec2:exact} is also included for comparison. We note that the LxF scheme converges to a wrong solution on very fine meshes, both in $h$ and $m$. Furthermore, an obvious numerical artifact is observed: the exact $m$ should be continuous across the bottom discontinuity due to the conservation of mass \cite[the first equation in (17)]{bernetti2008exact}. However, the numerical momentum forms a spurious spike at the bottom discontinuity. The height of the spike is not diminished as we refine the mesh. In addition, we apply the HLLC scheme to solve the same problem. In this numerical scheme, we use the same source discretization as the LxF scheme, and apply the HLLC numerical flux (\cite{harten1983upstream} and \cite[Section 10.4]{toro2001shock}) for $h$ and $m$. We observe a similar numerical artifact and the wrong convergence in Figure \ref{fig:db-nwb}. Note that LxF and HLLC schemes converge to slightly different limit functions, especially in the momentum $m$.  

Next we test this problem with several well-balanced schemes, including the wbLxF scheme \eqref{eq:wbLxF}, the HR scheme \eqref{eq:hr}, and the first-order version of XS scheme \eqref{eq:fsr}. %\zs{Has the first-order version of XS scheme appeared in earlier literature?} 
In these schemes, different ingredients are added to preserve the well-balanced property. In Section \ref{sec:schemes}, we will see that all these methods can essentially be viewed as modifications of the LxF scheme.

\begin{figure}[!h]
	\centering
	\subfloat[wbLxF, $h$.]{%
		\includegraphics[width=0.3\textwidth,trim={2cm 0 1cm 0}]{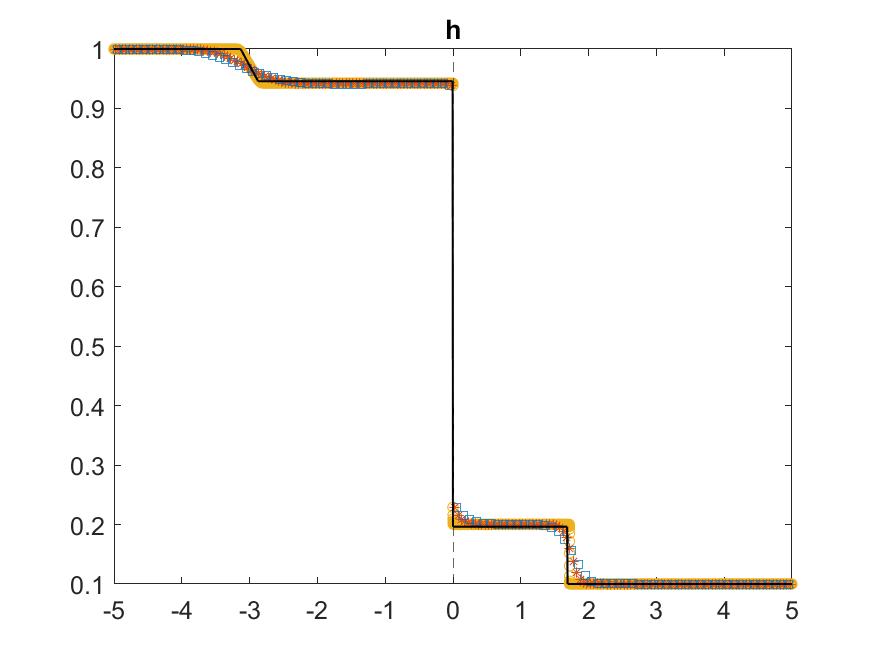}
	}
	\subfloat[HR, $h$.]{%
		\includegraphics[width=0.3\textwidth,trim={2cm 0 1cm 0}]{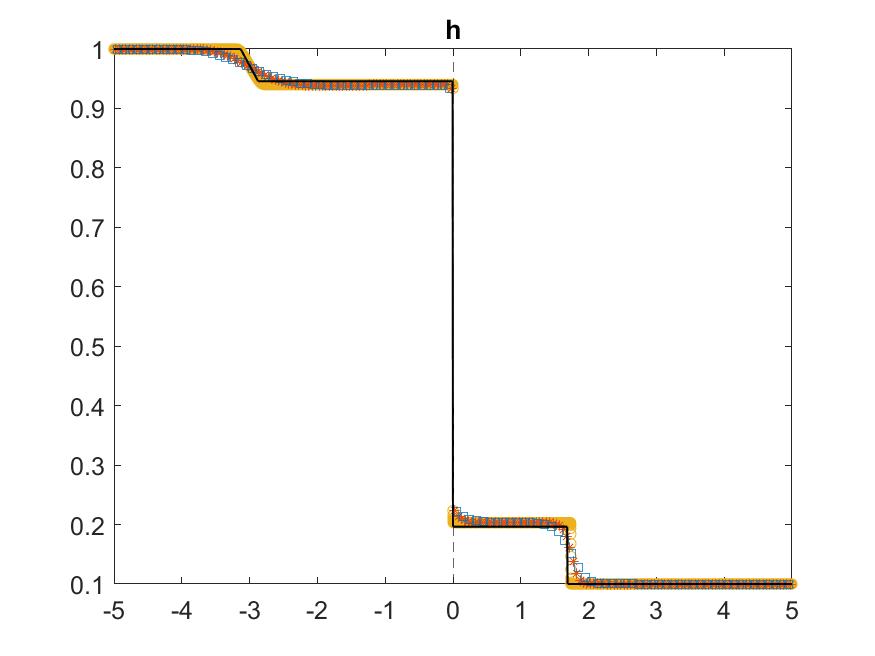}
	}
	\subfloat[XS, $h$.]{%
		\includegraphics[width=0.3\textwidth,trim={2cm 0 1cm 0}]{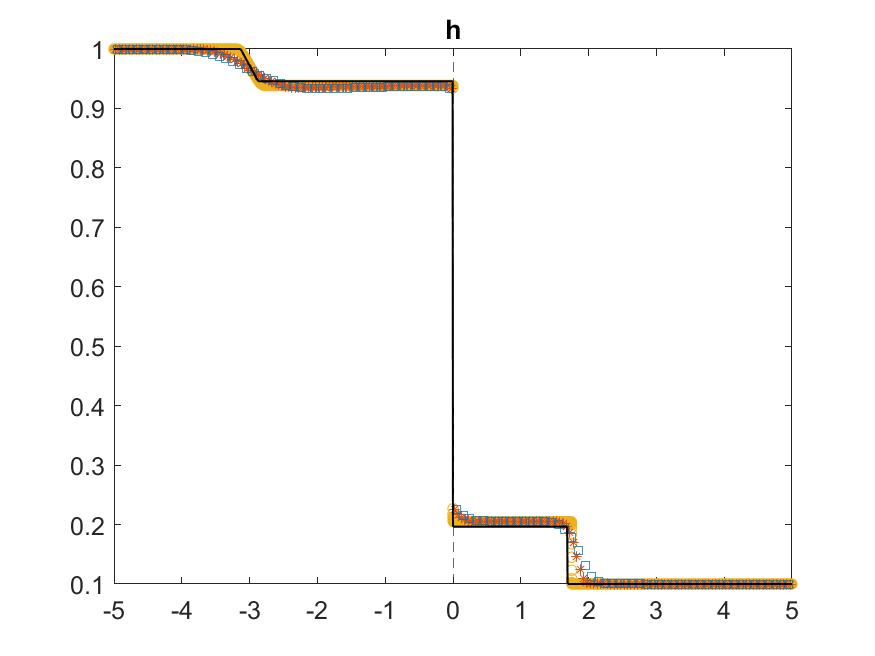}
	}\\
	\subfloat[wbLxF, $m= hu$.]{%
		\includegraphics[width=0.3\textwidth,trim={2cm 0 1cm 0}]{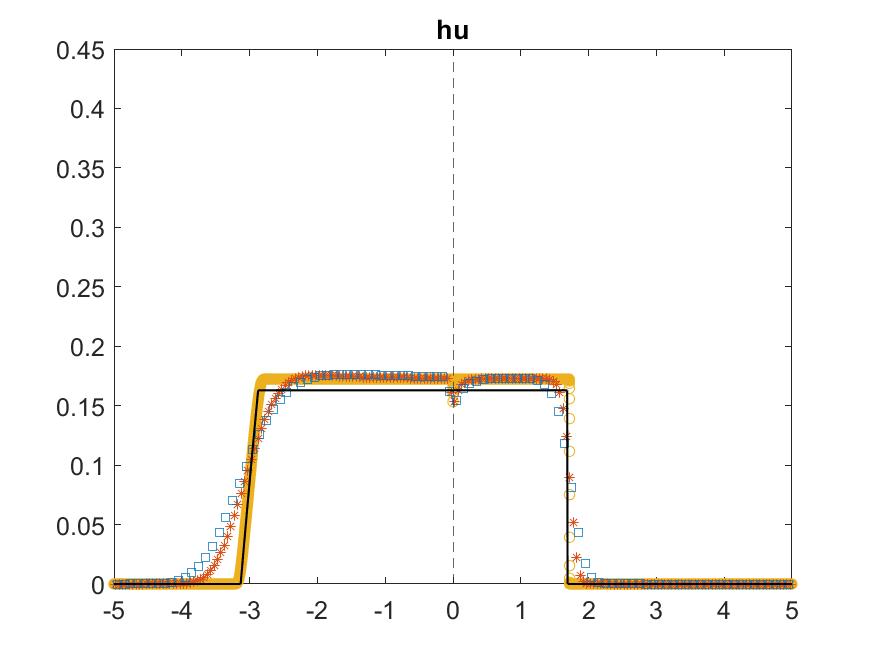}
	}
	\subfloat[HR, $m= hu$.]{%
		\includegraphics[width=0.3\textwidth,trim={2cm 0 1cm 0}]{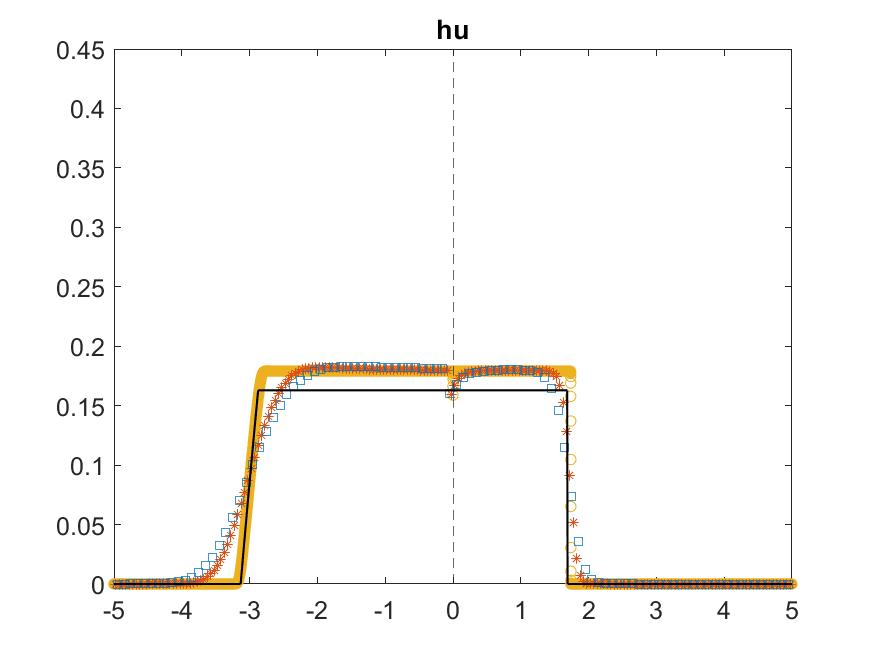}
	}
	\subfloat[XS, $m= hu$.]{%
		\includegraphics[width=0.3\textwidth,trim={2cm 0 1cm 0}]{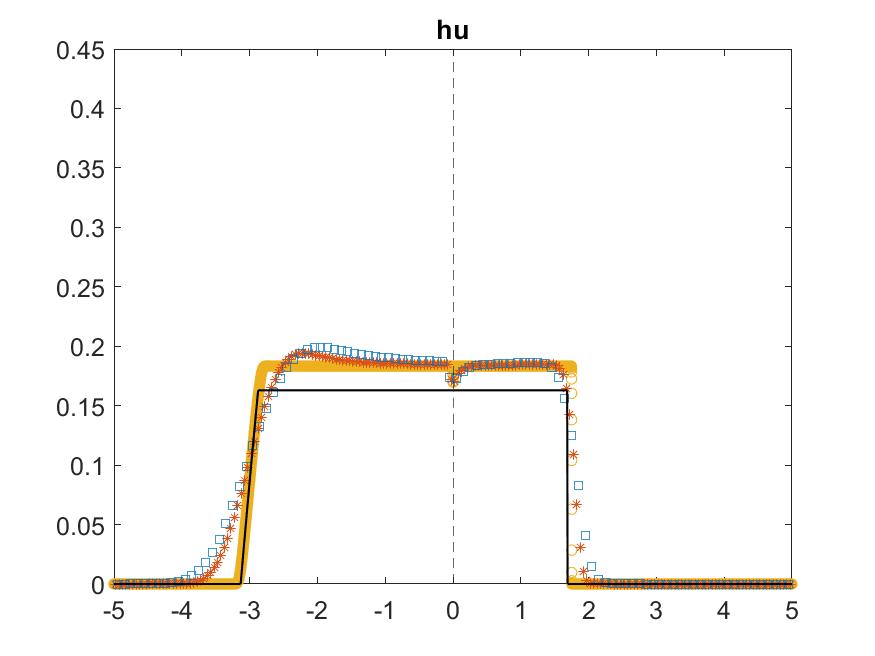}
	}\\
	\caption{Solution of well-balanced numerical schemes to the dam-break problem \eqref{eq:db}. Black solid line: exact solution; blue squares: $N = 100$; red stars: $N = 200$, yellow circles: $N = 25600$.}\label{fig:db-wb}
\end{figure}

\begin{figure}[!h]
	\centering
	\subfloat[$L^1$ error of $h$.]{%
		\includegraphics[width=0.35\textwidth,trim={2cm 0 1cm 0}]{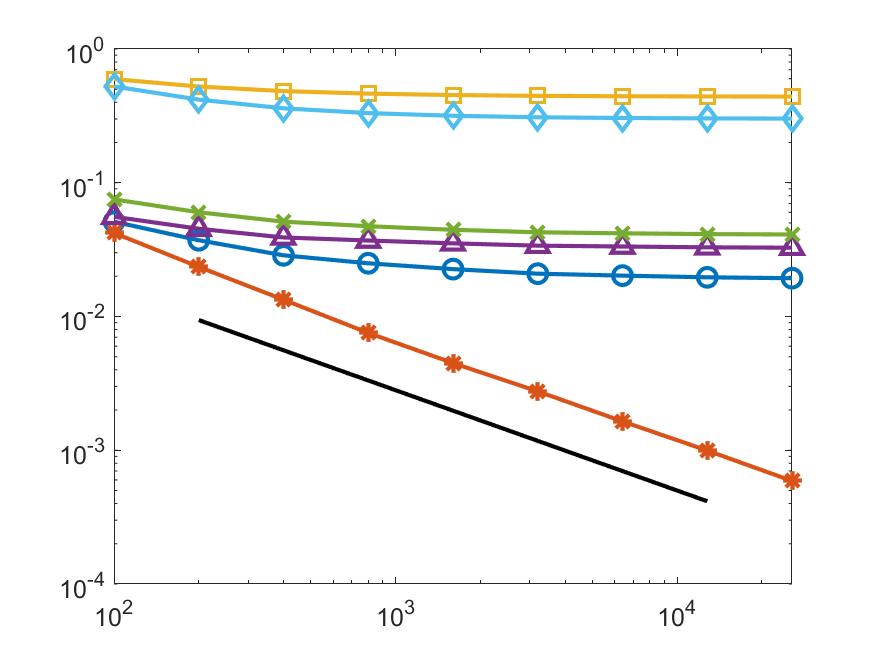}
		\label{fig:db_error_h}}
	\qquad
	\subfloat[$L^1$ error of $m = hu$.]{%
		\includegraphics[width=0.35\textwidth,trim={2cm 0 1cm 0}]{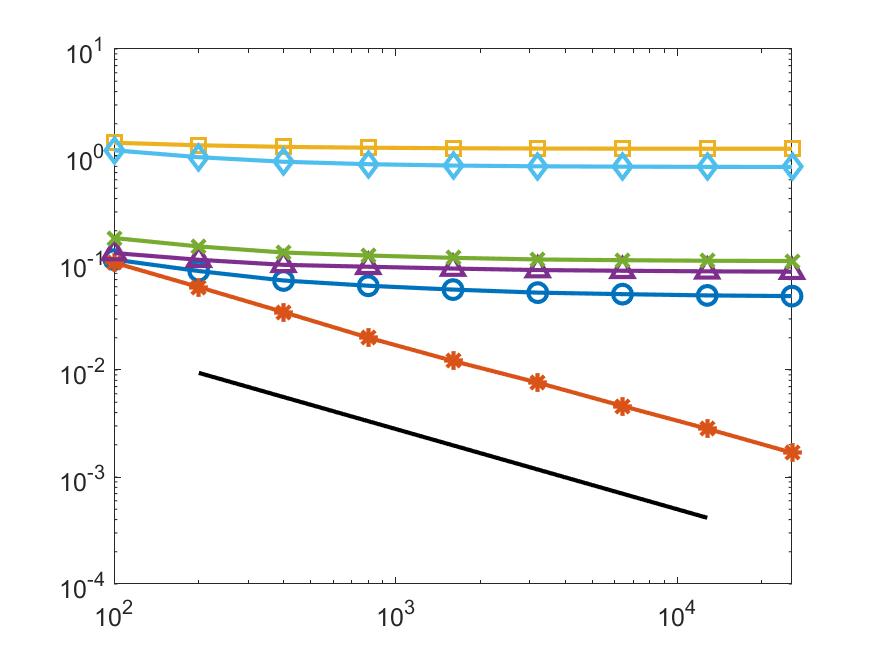}
	}
	\caption{$L^1$ error plot of the dam-break problem \eqref{eq:db}. The $x$ axis is the number of mesh points $N$, and the $y$ axis is the $L^1$ error. Yellow squares: LxF; cyan diamonds: HLLC;  blue circles: wbLxF; purple triangles: HR; green crosses: XS; red stars: cLxF; black solid line: reference slope with $\mathcal{O}(h^{3/4})$ convergence rate. }\label{fig:db-error}
\end{figure}

The numerical results of these well-balanced schemes are given in Figure \ref{fig:db-wb}. As we can see, these schemes seem to converge to solutions that are ``closer" to the exact solution compared with the non-well-balanced schemes. However, there is still a clear mismatch between the numerical and the exact momentum $m$. Furthermore, these schemes also suffer the numerical artifact of having a spurious spike at the bottom discontinuity (while for this test, it is less evident compared with the non-well-balanced schemes, and we will explain the reason in Remark \ref{rmk:wbnwbspike}). Even though it is hard to visually distinguish the numerical and exact water heights, we can see from the error plot Figure \ref{fig:db-error} that the convergence rates decrease to 0 and hence eventually all these schemes converge to wrong solutions.

\begin{figure}[!h]
	\centering
	\subfloat[cLxF, $h$.]{%
		\includegraphics[width=0.35\textwidth,trim={2cm 0 1cm 0}]{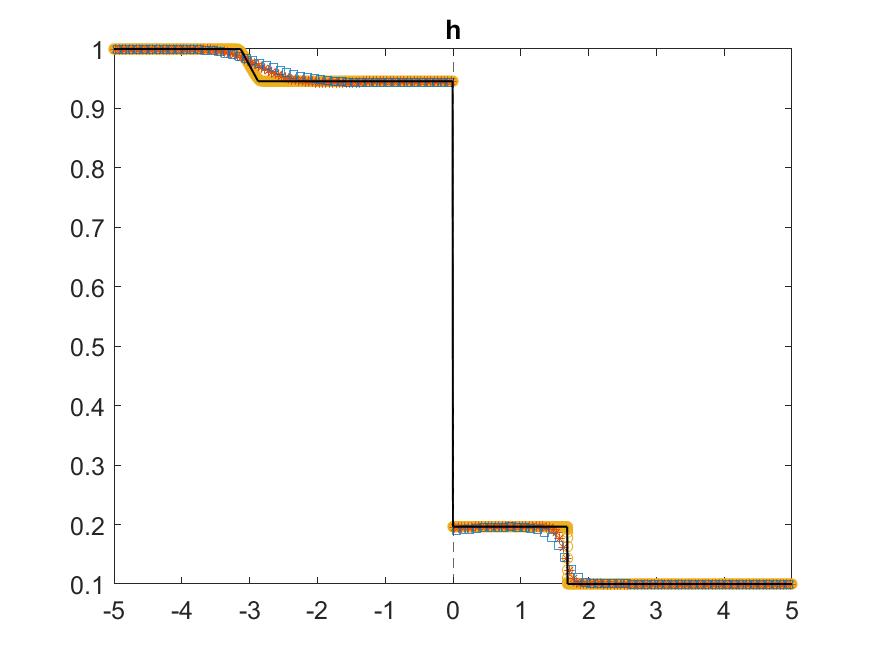}
	}
	\qquad
	\subfloat[cLxF, $m = hu$.]{%
		\includegraphics[width=0.35\textwidth,trim={2cm 0 1cm 0}]{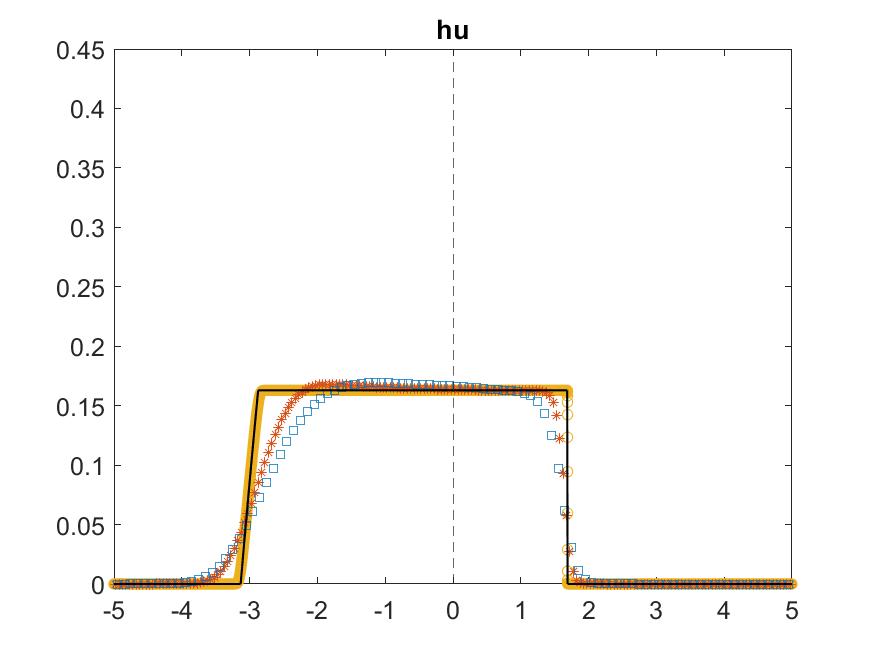}
	}
	\caption{Solution to the dam-break problem \eqref{eq:db} with cLxF scheme. Black solid line: exact solution; blue squares: $N = 100$; red stars: $N = 200$, yellow circles: $N = 25600$.}\label{fig:db-cLxF}
\end{figure}

However, if we make a very minor change in the LxF scheme, in which we remove the numerical viscosity at $x = 0$ by applying the central flux in the equation of mass conservation (which will be called the cLxF scheme), the numerical solution will converge to the exact Riemann solution. See Figure \ref{fig:db-cLxF}. From the error plot Figure \ref{fig:db-error}, it can be seen that the numerical solution converges at a rate of around 0.75 in $L^1$. %\YX{Just curious, is the expected $L^\infty$ norm 0-th order, and the $L^2$ norm $1/2$ order?} and 
Furthermore, there is no spurious spike in the numerical momentum. In the following sections, we will explain the formation of the spurious spike and the mechanism of why removing the numerical viscosity helps with the convergence. 

\section{Weak solutions and Riemann solutions}\label{sec:ws}
\setcounter{equation}{0}

In this section, we briefly review the definitions of the weak solutions and the Riemann solutions of the SWEs \eqref{eq-swe} with a discontinuous bottom topography. 
The main points of this section are summarized in below. 
\begin{enumerate}
	\item The weak solutions to the SWEs \eqref{eq-swe} are defined through a specific choice of the path. Physically, this choice of the path relates to the definition of the hydrostatic pressure at the bottom step. In this paper, we define the weak solutions in \eqref{eq-weak} with the path-related terms specified in \eqref{eq-S_phi} and \eqref{eq:sgnb}.
	\item There are controversial discussions on how the Riemann solutions should be defined at the bottom discontinuity. In the numerical tests of this paper, the exact solutions are generated based on the generalized Rankine--Hugoniot condition \cite{bernetti2008exact,rosatti2010riemann}.
\end{enumerate}

\subsection{Weak solutions}

\textbf{Weak solution and the path.} As is mentioned in Section \ref{sec:introduction}, the augmented nonconservative hyperbolic system \eqref{eq-ncsv} is considered. When defining weak solutions of \eqref{eq-ncsv}, complications arise at a point where $\tU$ is discontinuous. Due to the nonconservative product, the formal integration by parts no longer applies. In \cite{maso1995definition}, the authors propose to consider a smooth regularization in the neighborhood of discontinuity $x=x_\star$ to connect the two states
\begin{equation}
\tU^\veps = \left\{
\begin{array}{ll}
\tU^-,& x-x_\star \leq -\veps,\\
\Phi\left(\frac{x - x_\star + \veps}{2\veps};\tU^-,\tU^+\right),& - \veps < x - x_\star  \leq \veps,\\
\tU^+,& x-x_\star > \veps,
\end{array}
\right.
\end{equation}
and then consider the limit $\veps \to 0^+$. Here $x_\star$ is the point separating the left state $\tU^-$ and the right state $\tU^+$. %\YX{Do we have $x_\star=0$ based on the assumption in Section 1?} 
The Lipschitz function $\Phi: [0,1]\times \tOmega \times \tOmega \to \tOmega^3$, satisfying
\begin{equation}
 \Phi\left(0;\tU^-,\tU^+\right) = \tU^- \quand\Phi\left(1;\tU^-,\tU^+\right) = \tU^+,
\end{equation}
is called a \emph{path} connecting the two states. It turns out that, once the path $\Phi$ is given, one can interpret $A(\tU)\tU_x$ as a Borel measure, formally defined as
\begin{equation}
[A(\tU)\tU_x]_\Phi = A(\tU)\tU_x \dd x + \sum_{x_\star} \left(\int_0^1 A\left(\Phi(s;\tU^-,\tU^+)\right) \frac{\partial\Phi}{\partial s}\left(s;\tU^-,\tU^+\right) \dd s\right) \delta(x_{\star}), 
\end{equation}
where the summation is taken over all discontinuous points $x_\star$ and $\delta$ is the Dirac measure. Then we can consider \eqref{eq-ncsv} as a measure-valued equation, 
\begin{equation}\label{eq-measeq}
\tU_t + [A(\tU)\tU_x]_\Phi = 0.
\end{equation} 
Across the discontinuity, the following generalized Rankine--Hugoniot condition holds, 
\begin{equation}\label{eq-gRH}
\xi \left(\tU^+ - \tU^-\right) = \int_0^1 A\left(\Phi(s;\tU^-,\tU^+)\right)\frac{\partial \Phi}{\partial s}\left(s,\tU^-,\tU^+\right)\dd s,
\end{equation}
where $\xi$ is the shock speed. In the conservative case, $A$ is the Jacobian of a flux function and one can apply the chain rule and the Newton--Leibniz formula to retrieve the classical Rankine--Hugoniot condition.

For $\tU=(h,m,b)^T$, let us denote the path by $\Phi = (\Phi^h, \Phi^\mm, \Phi^b)^T$ and choose the test function $\Psi = (\zeta, \eta, \theta)^T \in C_c^1((-\infty,\infty)\times [0,\infty))$. Using these component-wise notations and the definitions of $A$, $S$ and $\tU$, \eqref{eq-measeq} becomes
\begin{equation}\label{eq-weak-1}
\begin{aligned}
&- \int_0^\infty\int_{-\infty}^\infty U\cdot \Psi_t  \dd x \dd t -\int_0^\infty U^0\cdot \Psi \dd t- \int_0^\infty\int_{-\infty}^\infty F(U) \cdot \Psi_x \dd x \dd t\\
&\hskip2cm  
=\int_0^\infty\int_{-\infty}^\infty S(U) b_x \cdot \Psi \dd x \dd t +  \int_0^\infty \sum_{x_\star} S_\Phi^h(x_\star,t)\eta(x_\star,t) \dd t,
\end{aligned}
\end{equation}
where
\begin{equation}\label{eq-sphi}
S_\Phi^h = -  g\int_0^1 \Phi^h\left(s;\tU^-,\tU^+\right) \frac{\partial \Phi^b}{\partial s}  \left(s;\tU^-,\tU^+\right)\dd s. 
\end{equation}
Clearly, the weak solution to \eqref{eq-weak-1} varies with the choice of the path $\Phi$ and one needs to specify $\Phi^h$ and $\Phi^b$ to close the definition of $U$. Note $b$ is assumed to be continuous except at $x = 0$. One typically takes $\Phi^b \equiv b$ to be independent of $s$ at all other discontinuities of $h$ and $m$ (where $b$ is continuous). Therefore, according to \eqref{eq-sphi}, $S_\Phi^h$ is nonzero only at $x = 0$. Hence one can drop the summation on the right-hand side of \eqref{eq-weak-1} to rewrite it as  
\begin{equation}\label{eq-weak}
	\begin{aligned}
		&- \int_0^\infty\int_{-\infty}^\infty U\cdot \Psi_t \dd x \dd t -\int_0^\infty U^0\cdot \Psi \dd t- \int_0^\infty\int_{-\infty}^\infty F(U) \cdot \Psi_x\dd x \dd t \\
		&\hskip2cm  
		=\int_0^\infty\int_{-\infty}^\infty S(U) b_x \cdot \Psi \dd x \dd t +  \int_0^\infty  S_\Phi^h(0,t)\eta(0,t) \dd t.
	\end{aligned}
\end{equation}
In other words, the path $\Phi$ only needs to be introduced at the discontinuity of $b$ (at $x = 0$ in our settings).  

\textbf{The path and the hydrostatic pressure.} To give an appropriate definition of $S_\Phi^h$ at $x = 0$, we look into another stream of study on solutions of \eqref{eq-ncsv}, which concerns the Riemann problem of \eqref{eq-swe} and uses the controlled volume across the discontinuity to derive the Rankine--Hugoniot condition. See, for example, \cite{bernetti2008exact}. In this section, we give a heuristic derivation of the weak solution following the idea in \cite{bernetti2008exact}. We will only consider the momentum equation to address the effect of the source term with a discontinuous bottom. 

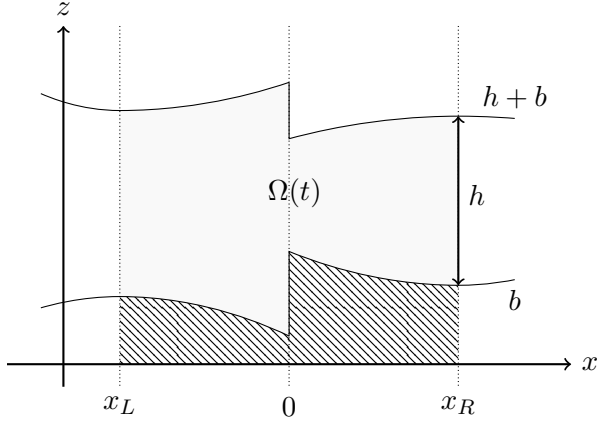
\begin{figure}[h!]
	\centering
	\begin{tikzpicture}[scale=1.5]
		
		\draw[->, thick] (-0.5,0) -- (4.5,0) node[right] {$x$};
		\draw[->, thick] (0,-0.2) -- (0,3) node[above] {$z$};
		
		% fill left
		\filldraw[gray!5] (.5,0.6) parabola (2,0.25) -- (2,2.5);
		\filldraw[gray!5] (.5,2.25) parabola (2,2.5) -- (2,0.25);
		\filldraw[gray!5] (.5,2.25) parabola (0.5,0.6) -- (2,1);
		% fill right
		\filldraw[gray!5] (3.5,0.7) parabola (2,1) -- (2,2);
		\filldraw[gray!5] (3.5,2.2) parabola (2,2) -- (2,1);
		\filldraw[gray!5] (3.5,0.7) parabola (3.5,2.2) -- (2,1.5);
		\draw (2.05,1.3) node[above] {$\Omega(t)$};
		
		% Riverbed
		\fill[pattern=north west lines] (.5,0.6) parabola (2,0.25) -- (2,0) -- (0.5,0);
		\fill[pattern=north west lines] (3.5,0.7) parabola (2,1) -- (2,0) -- (3.5,0);
		% lower
		\draw (-.2,0.5) parabola bend (.5,0.6) (2,0.25);
		\draw (4,0.75) parabola bend (3.5,0.7) (2,1) -- (2,0.25);
		\draw (4,0.75) node[below] {$b$};
		
		%upper
		\draw (-.2,2.4) parabola bend (.5,2.25) (2,2.5);
		\draw (4,2.18) parabola bend (3.5,2.2) (2,2) -- (2,2.5); 
		\draw (4,2.18) node[above] {$h+b$};
		
		%h
		\draw[<->, thick] (3.5,0.7) -- (3.5,2.2);
		\draw (3.5,1.5)  node[right] {$h$};
		
		\draw[densely dotted] (0.5,3) -- (0.5,-0.2)  node[below] {$x_L$};
		\draw[densely dotted] (3.5,3) -- (3.5,-0.2)  node[below] {$x_R$};
		\draw[densely dotted] (2,3) -- (2,-0.2)  node[below] {$0$};
		
	\end{tikzpicture}
	\caption{Control volume near the discontinuity of $b$.}
 \label{fig1}
\end{figure}

Assume the velocity $u = u(x)$ to be independent of the water depth $z$. From the conservation of momentum over a control volume $\Omega(t)$ surrounding the discontinuity at $x = 0$, we have 
\begin{equation}\label{eq-mom}
	\frac{\dd}{\dd t}\int_{\Omega(t)} \rho u \dd x \dd z + \int_{\partial \Omega(t)} p \nux\dd s = 0.
\end{equation}
Here $\rho$ is the water density (constant), $p$ is the hydrostatic pressure to be specified later, and $\nux$ is the horizontal component of the unit outer normal along the boundary $\partial \Omega(t)$. See Figure \ref{fig1} for an illustration of the control volume $\Omega(t)$.

By substituting in the definitions of $\Omega$, $\partial\Omega$ and $\nu^x$, one can obtain (see Appendix \ref{app:weaksoln}) 
\begin{equation}\label{eq-mom-cont}
	\begin{aligned}
		&\int_{x_L}^{x_R} (h u)_t  \dd x + \left(hu^2+\hf g h^2\right)\bigg|_{x=x_R} -   \left(hu^2+\hf g h^2\right) \bigg|_{x=x_L}= - \int_{x_L}^{x_R} g h b_x \dd x - \int_{b^-}^{b^+} \frac{ p(0,z) }{\rho}\dd z. 
	\end{aligned}
\end{equation}
Then using a Lax--Wendroff type argument, we get (details are omitted)
\begin{equation}\label{eq:weak-mom}
	\begin{aligned}
	&-\int_0^\infty \int_{-\infty}^\infty h u \eta_t  \dd x\dd t  - \int_0^\infty hu(x,0) \eta(x,0) dx - \int_0^\infty \int_{-\infty}^\infty \left(hu^2 + \frac{1}{2}gh^2\right) \eta_x \dd x \dd t\\
	&\hskip2cm
	=-\int_0^\infty \int_{-\infty}^\infty  gh b_x \eta \dd x \dd t  - \int_0^\infty \left(\int_{b^-}^{b^+} \frac{p(0,z)}{\rho} \dd z\right) \eta(0,t)  \dd x \dd t.
		\end{aligned}
\end{equation}
Comparing \eqref{eq-weak} with \eqref{eq:weak-mom}, it can be seen that different choices of the path in \eqref{eq-weak} may correspond to different definitions of the pressure $p(0,z)$ in \eqref{eq:weak-mom}. The two definitions of the weak solution coincide if 
\begin{equation}\label{eq:S_phi-p}
	S_\Phi^h(0,t) = - \int_{b^-}^{b^+} \frac{p(0,z)}{\rho} \dd z.
\end{equation}
A similar relation has also been derived in \cite[equation (22)]{cozzolino2011numerical}.

One usually assumes the hydrostatic pressure $p = \rho g(h + b - z)$ in SWEs. But at the bottom discontinuity, $h+b$ is double-valued. Hence we consider the weighted average of $h+b$ and define 
\begin{equation}\label{eq:p}
p(0, z) = \rho g(\hb{h+b} -z),
\end{equation}
where we have used the notations 
\begin{equation}\label{eq-hb}
\hb{w} = \{w\} - \frac{\gamma}{2}[w], \quad \{w\} = \hf\left(w^+ + w^-\right), \quad \jp{w} = w^+ - w^-,
\end{equation}
and $w^\pm = \lim_{\delta \to 0^\pm}w(x+\delta)$ are the left and right limits of $w$ at $x = 0$. 
In the literature \cite{fraccarollo2002riemann,rosatti2006well,bernetti2008exact,cozzolino2011numerical}, the pressure is usually taken from the lower side of the bottom step, corresponding to 
\begin{equation}\label{eq:sgnb}
	\gamma = \sgnb.
\end{equation}
Substituting \eqref{eq:p} into \eqref{eq:S_phi-p} and using the fact $[b^2] = 2\{b\}[b]$, one can get
\begin{equation}\label{eq-S_phi}
	S_\Phi^h(0,t) = -  g\check{h}_\gamma[b], \quad \text{where } \check{h}_\gamma := \hb{h+b} -\{b\},
\end{equation}
which corresponds to the choice of the path:
\begin{subequations}\label{eq-path}
\begin{align}
	\Phi^h\left(s;\tU^-,\tU^+\right) =& h^- + \left([h]-3\gamma[h+b]\right)s + 3\gamma[h+b]s^2,\\
	\Phi^b\left(s;\tU^-,\tU^+\right) =& b^- + [b]s.
\end{align}	
\end{subequations}
for $S_\Phi^h$ defined in \eqref{eq-sphi}. In \eqref{eq-path}, $\Phi^h$ and $\Phi^b$ are Lipschitz continuous in $\tU$. 

In summary, the weak solution to \eqref{eq-swe} is defined by \eqref{eq-weak} with $S_\Phi^h$ specified in \eqref{eq-S_phi}, or the path specified in \eqref{eq-path}, and $\gamma$ specified in \eqref{eq:sgnb}.

\subsection{Riemann solutions}\label{sec:riemann}

There have been several works on studying Riemann solutions to the SWEs over a bottom step. For the nonconservative augmented system \eqref{eq-ncsv}, it can be seen that its eigenstructure shares some similarities with the homogeneous SWEs (without the source term). Indeed, two of the eigenvalues in \eqref{eq-ncsv} are equivalent to those in the homogeneous SWEs, leading to two genuinely nonlinear waves that can develop either shock or rarefaction waves. It also includes a third eigenvalue that is identically equal to zero, and it corresponds to a contact wave appearing only in presence of the bottom discontinuity. In the literature, there has been some discussions regarding how the flow variables should be connected at this contact wave. In particular, the following two approaches have been presented.

\textbf{The first approach: Riemann invariants.} The first approach is based on the mass and energy conservation principle. They consider a reformulation of \eqref{eq-swe} in the form
\begin{equation}\label{eq:prim}
	\begin{pmatrix}
		h\\
		u
	\end{pmatrix}_t + 
\begin{pmatrix}
	hu\\
	\frac{u^2}{2}+g(h+b)
\end{pmatrix}_x= \begin{pmatrix}
0\\0
\end{pmatrix}. 
\end{equation}
According to this formulation, the stationary wave at the bottom discontinuity should satisfy
\begin{subequations}\label{eq:ri}
\begin{align}
	[hu] =& 0,\\
	\left[\frac{u^2}{2}+g(h+b)\right] =&
	0.\label{eq:ri-2}
\end{align}
\end{subequations}
This condition ensures the constancy of the Riemann invariants, namely the mass and the energy. Many works in the literature have been adopting this approach, see  \cite{alcrudo2001exact,bukreev2003breakdown,andrianov2005performance,chinnayya2004well,gallouet2003some,lefloch2007riemann,han2014exact,aleksyuk2019uniqueness} for an incomplete list of references. 

\textbf{The second approach: the generalized Rankine--Hugoniot condition.} In contrast, the second approach is based on the mass and momentum conservation, which is used by Bernetti et al. in \cite{bernetti2008exact} and then further investigated by Rosatti and Begnudelli in \cite{rosatti2010riemann}. This approach is based on integral form of the SWEs as have been derived in the previous section
\begin{subequations}\label{eq:swe-conserve}
	\begin{align}
	\int_{x_L}^{x_R} h_t  \dd x + hu\big|^{x=x_R}_{x=x_L}=& 0,\\ 
			\int_{x_L}^{x_R} (h u)_t  \dd x + \left(hu^2+\hf g h^2\right)\bigg|^{x=x_R}_{x=x_L}=& - \int_{x_L}^{x_R} g h b_x \dd x - \int_{b^-}^{b^+} \frac{ p(0,z) }{\rho}\dd z. 
\end{align}
\end{subequations}
Recall that $-\int_{b^-}^{b^+}{p(0,z)}/{\rho}\dd z = S_\Phi^h(0,t) = -g\check{h}_\gamma[b]$. After taking $x_R\to 0^+$ and $x_L\to 0^-$, one can see that the stationary wave at the bottom discontinuity admits the following conditions
\begin{subequations}\label{eq:GRH}
	\begin{align}
	[hu] =& 0,\\
	\left[hu^2+\hf gh^2\right] =& -  g\check{h}_\gamma[b].\label{eq:GRH-2}
\end{align}
\end{subequations}
\eqref{eq:GRH} is consistent with \eqref{eq-gRH} and is also referred to as the generalized Rankine--Hugoniot condition. 

\textbf{Exact Riemann solver in this paper.} One can see that \eqref{eq:ri} and \eqref{eq:GRH} are not mathematically equivalent, at the place where the bottom topography is discontinuous. Hence the Riemann solutions derived through these two approaches will be different. Detailed discussion of these two approaches is presented in \cite{rosatti2010riemann}. Note that it is demonstrated in \cite[Theorem 1]{rosatti2010riemann} that the two solutions will coincide with each other if the integral curve (IC) coincides with the generalized Hugoniot locus (HL). Although we have the fact that ``The HL curves and the IC curves coincide in case of contact waves in conservative systems. In other words, both Riemann invariants and Rankine--Hugoniot relations hold at the same time." \cite[Section 3.4]{rosatti2010riemann}, unfortunately, in general, ``In contact waves of nonconservative systems where the relevant eigenvalue does not depend on $U$, IC and generalized HL may not coincide."\cite[Corollary 1]{rosatti2010riemann}. We refer to \cite{rosatti2010riemann} for further details  on the difference and comparison of the two approaches. %\YX{cite few comments/summary from that paper?}

In our numerical tests, we adopt the second approach using the generalized Rankine-Hugoniot condition to generate the exact solutions. The reason is that we will use Lax--Wendroff type argument as a tool for the convergence analysis, and it agrees more naturally with the second approach. To be more specific, the numerical schemes we consider in this paper are designed by discretizing the unknowns $h$ and $m$ as that in \eqref{eq:swe-conserve}, rather than $h$ and $u$ as that in \eqref{eq:prim}. In the case $b$ is constant, the numerical schemes will retrieve the classical Lax--Friedrichs scheme. It preserves the local conservation of mass and momentum, but it doesn't have a mechanism to preserve the local energy conservation automatically. Therefore, the classical Lax--Wendroff theorem will show that the corresponding solution limit satisfies \eqref{eq-swe} with $b$ being constant, instead of \eqref{eq:prim} with $b$ being constant, in the weak sense. In the case that $b$ has a jump discontinuity at $x = 0$, the solution limit by the numerical schemes may admit a relationship in the similar form as \eqref{eq:GRH-2}, but recovering \eqref{eq:ri-2} could be difficult. 
In addition, in our numerical tests, we do observe the correct convergence to the exact solution generated by the second approach in $L^1$ norm for the cLxF scheme for nontransonic tests and the LxF scheme for negative supercritical tests. %Hence, we believe that using the second approach to generate the reference solution is, at the very least, not wrong. 
Moreover, regardless of how the exact solution is defined, the spurious spike we observe seems to be unphysical. The definition of the exact Riemann solution should not affect the main results presented in our paper.

%\YX{consider how to rewrite this, emphasizing that we will use the Lax-Wendroff as a tool and adopt the second approach, while not saying any bad about the first approach?}

%In this special case, When analyzing the convergence by proving a Lax--Wendroff theorem, we expect the solution limit should satisfy an equation more similar to \eqref{eq:GRH} rather than \eqref{eq:prim}. 

%Despite the wide use of the first approach in the literature, Rosatti and Begnudelli argued in \cite{rosatti2010riemann} that one of the Riemann invariant, namely the energy, should not be constant wave across the bottom discontinuity due to the nonconservative nature of the problem. They further demonstrated that, in this type of system, Riemann invariants do not generally hold in contact waves. We refer to \cite{rosatti2010riemann} for details. 

\section{Numerical schemes}\label{sec:schemes} 
\setcounter{equation}{0}

In this section, we detail the numerical schemes which have been tested in Section \ref{sec:artifact}. We will omit the HLLC scheme and focus on LxF type schemes only.  We will derive a unified formulation and rewrite three well-balanced schemes as the LxF scheme added with some source modification terms. This will facilitate the convergence analysis in Section \ref{sec:ca}. 

Let us consider a mesh partition of the spacetime domain $(-\infty, \infty)\times[0,\infty) = \cup_{j,n} \left(I_j\times I^n\right)$ with $I_j = [x_{j-\hf}, x_{j+\hf})$, and $I^n = [t^n,t^{n+1})$. We define $x_{j}= (x_{j+\hf}+x_{j-\hf})/2$, $\dx_j = x_{j+\hf}- x_{j-\hf}$, $\dt^n = t^{n+1}-t^n$, $\dx = \max_j \dx_j$ and $\dt = \max_n \dt^n$. The numerical solution is denoted as $U_j^n$, which approximates the exact solution $U(x_j, t^n)$. Without loss of generality, we assume that the jump discontinuity of $b$ locates at $x_\hf = 0$. 

In this paper, we consider the LxF type schemes in the following form
\begin{equation}\label{eq:pcs}
\frac{U_j^{n+1} -U_j^n}{\dt^n} + \frac{\wh{F}_{j+\hf}^n-\wh{F}_{j-\hf}^n}{\dx_j}= \frac{\wh{S}_{j+\hf}^{n}+\wh{S}^{n}_{j-\hf}}{\dx_j} + \frac{\wh{M}_{j+\hf}^{n}-\wh{M}^{n}_{j-\hf}}{\dx_j},
\end{equation}
where we have (after dropping all superscripts $n$)
%\begin{align}
%\widehat{{F}}_{j+\hf}^n &= \wh{F}(U_j^n,U_{j+1}^n) \\
%&=  \hf\left(F(U_j^n) + F(U_{j+1}^n)\right) - \frac{\alpha}{2}\left(U_{j+1}^n-U_j^n\right) = \{F(U^n)\}_{j+\hf} - \frac{\alpha}{2}[U^n]_{j+\hf},\\
%{S}_{j+\hf}^{n,\pm} =& \left(
%\begin{array}{c}
%0\\
% \frac{g}{2} \left(\frac{1+\gamma}{2}w_{j+1}^n+\frac{1-\gamma}{2}w_j^n - \hf(b_{j+1} + b_j))\right)
%\end{array}\right) = \left(
%\begin{array}{c}
%0\\
%\frac{g}{2} \left(\{w^n\}+ \frac{\gamma}{2}[w^n] - \{b\}\right)_{j+\hf}
%\end{array}\right),\\
%\wh{M}_{j+\hf}^n =& \left(
%\begin{array}{c}
%\alpha/2\\
%m 
%\end{array}\right)(b_{j+1}-b_j)= \left(
%\begin{array}{c}
%\alpha/2\\
%m 
%\end{array}\right)[b]_{j+\hf}.
%\end{align}
\begin{subequations}\label{eq:FSM}
\begin{align}
%\widehat{{F}}_{j+\hf} =& \wh{F}\left(U_j,U_{j+1}\right)= \{F(U)\}_{j+\hf} - \frac{\alpha_{j+\hf}}{2}[U]_{j+\hf},\quad \alpha_{j+\hf} = \max_{l = j,j+1}\left(\left|\frac{\mm_{l}}{h_{l}}\right| + \sqrt{gh_{l}}\right),\label{eq:lxfflux}\\
\widehat{{F}}_{j+\hf} =& \wh{F}\left(U_j,U_{j+1}\right)= \{F(U)\}_{j+\hf} - \frac{1}{2}A_{j+\hf}[U]_{j+\hf},\quad A_{j+\hf} = \begin{pmatrix}
	\alpha_{1,j+\hf}&0\\0&\alpha_{2,j+\hf}
\end{pmatrix}, \label{eq:lxfflux}\\
\wh{S}_{j+\hf} =& -\frac{g}{2}\left(
\begin{array}{c}
0\\
\check{h}_{\gamma,j+\hf}
\end{array}\right)[b]_{j+\hf}, \qquad \check{h}_{\gamma, j+\hf} = \check{h}_\gamma\left(h_j,h_{j+1}\right) = \hbn{h+b}_{\gamma,j+\hf} - \{b\}_{j+\hf},\label{eq:lxfsource}\\
\wh{M}_{j+\hf} =& \hf \widehat{N}_{j+\hf}[b]_{j+\hf}, \qquad \qquad \qquad   \widehat{N}_{j+\hf} = \widehat{N}\left(U_j,U_{j+1}\right).\label{eq:lxfM}
\end{align}
\end{subequations}
Here if the local LxF flux is used, we usually choose 
	\begin{equation}\label{eq:alpha}
		A_{j+\hf} = A\left(U_j,U_{j+1}\right) = 
		\begin{pmatrix}
			\alpha_1\left(U_j,U_{j+1}\right)&0\\0&\alpha_2\left(U_j,U_{j+1}\right)
		\end{pmatrix}= \alpha\left(U_j,U_{j+1}\right)
	\begin{pmatrix}
	1&0\\0&1
\end{pmatrix}
	\end{equation} with
	\begin{equation}
		\alpha(U_j,U_{j+1}) = \max_{k = j,j+1}\left(\left|\frac{\mm_{k}}{h_{k}}\right| + \sqrt{gh_{k}}\right) := \alpha_{j+\hf}.
	\end{equation}
$\widehat{N}$ is an extra term to be specified later to accommodate different numerical methods. We remark here that the term $\check{h}_{\gamma,j+\hf}$ is an approximation of $\check{h}_\gamma$ in \eqref{eq-S_phi} at $x_{j+\hf}$, which relates to the definition of the pressure in the SWEs. In \eqref{eq:FSM}, we have adopted the notations in \eqref{eq-hb}, namely, 
\begin{align}
	v_{j+\hf}^- = v_j, \quad v_{j+\hf}^+ = v_{j+1}, \quad \{v\}_{j+\hf} = \frac{1}{2}\left(v_j+v_{j+1}\right), \quad [v]_{j+\hf} = v_{j+1} - v_{j},
\end{align}
and 
\begin{equation}
	\hbn{v}_{\gamma,j+\hf} = \{v\}_{j+\hf} - \frac{\gamma}{2}[v]_{j+\hf}.
\end{equation}

Next, we present the numerical schemes discussed in Section \ref{sec:artifact}, and show that they can be reformulated as \eqref{eq:pcs} with specially defined $\widehat{N}$. For ease of notation, we may omit the superscript $n$ and the subscript $j\pm \hf$ when it does not cause confusion. 

\textbf{(1) Lax--Friedrichs scheme (LxF scheme).} The simplest case is to take $\widehat{N}= 0$, which gives

\begin{equation}\label{eq:LxF}
\frac{U_j^{n+1}-U_j^n}{\Delta t^n} + \frac{\wh{F}_{j+\hf}^n-\wh{F}_{j-\hf}^n}{\dx_j}= \frac{\wh{S}_{j+\hf}^n+\wh{S}_{j-\hf}^n}{\dx_j}.
\end{equation}

\textbf{(2) The well-balanced Lax--Friedrichs  scheme (wbLxF scheme).} One can enforce the well-balanced property of \eqref{eq:LxF} by modifying the numerical flux, which takes the form
\begin{equation}\label{eq:wbLxF}
\frac{U_j^{n+1}-U_j^n}{\Delta t^n} + \frac{\wh{F}_{j+\hf}^{b,n}-\wh{F}_{j-\hf}^{b,n}}{\dx_j}= \frac{\wh{S}_{j+\hf}^n+\wh{S}_{j-\hf}^n}{\dx_j},
\end{equation}
where
\begin{equation}\label{eq:fsrFhat}
\wh{F}^b = \wh{F} - \frac{1}{2}A\begin{pmatrix}
[b]\\
0
\end{pmatrix}.
\end{equation}
Clearly, this numerical scheme can be written in the form of \eqref{eq:pcs}. 
\begin{THM}\label{thm:wblxf}
	The wbLxF scheme \eqref{eq:wbLxF} and \eqref{eq:fsrFhat} can be reformulated as \eqref{eq:pcs} with 
	\begin{equation}
	\gamma = \sgnb\quand \widehat{N} =\begin{pmatrix}
	\alpha_1\\
	0
	\end{pmatrix}.
	\end{equation}
\end{THM}

\textbf{(3) The well-balanced hydrostatic reconstruction scheme (HR scheme).} This scheme is initially proposed by Audusse et al. in \cite{audusse2004fast} and the first-order version is given as follows:
\begin{equation}\label{eq:hr}
\frac{U_j^{n+1}-U_j^n}{\Delta t^n} + \frac{\wh{F}_{j+\hf}^{n,*,-} - \wh{F}_{j-\hf}^{n,*,+}}{\dx_j} = 0,
\end{equation}
where
\begin{align}\label{eq:hrFhat}
\wh{F}^{*,\pm} =&  
\wh{F}^* + \left(
\begin{array}{c}
0\\
\frac{g}{2}(h^\pm)^2 - \frac{g}{2}(h^{*,\pm})^2
\end{array}\right),
\end{align}
and 
\begin{align}\label{eq:hrFhat2}
\wh{{F}}^* :=  \wh{F}\left(U^{*,-},U^{*,+}\right),\quad U^{*,\pm} = \begin{pmatrix}
h^{*,\pm}\\
\mm^\pm
\end{pmatrix}, \quad h^{*,\pm} = h^\pm+b^\pm-\max(b^-,b^+).
\end{align}
\begin{THM}\label{thm-HR}
	The HR scheme \eqref{eq:hr}, \eqref{eq:hrFhat} and \eqref{eq:hrFhat2} can be reformulated as \eqref{eq:pcs} with 
	\begin{equation}
		\gamma = \sgnb\quand \widehat{N} =\begin{pmatrix}
			\alpha_1\\
			 \hb{\frac{\mm^2}{h}}\frac{\mathrm{sgn}([b])}{|[b]|-\hb{h}}
		\end{pmatrix}.
	\end{equation}
\end{THM}

The proof of Theorem \ref{thm-HR} is given in Appendix \ref{app-HR}. 

\textbf{(4) The well-balanced scheme with flux and source modification (XS scheme).} Here we consider the first-order version of the well-balanced scheme proposed in \cite{xing2006high}:
\begin{equation}\label{eq:fsr}
\frac{U_j^{n+1}-U_j^n}{\Delta t^n}+ \frac{\wh{F}_{j+\hf}^{b,n} - \wh{F}_{j-\hf}^{b,n}}{\dx_j} = \frac{1}{\dx_j}\begin{pmatrix}
0\\
s_j^n
\end{pmatrix},
\end{equation}
where $\wh{F}^b$ is defined in \eqref{eq:fsrFhat} and
\begin{equation}\label{eq:fsrS}
s_j^n = \frac{g}{2}\left(\{b^2\}_{j+\hf} - \{b^2\}_{j-\hf}\right) - g(h_j^n + b_j)\left(\{b\}_{j+\hf} - \{b\}_{j-\hf}\right).
\end{equation}

\begin{THM}\label{thm-XS}
The XS scheme \eqref{eq:fsr}, \eqref{eq:fsrFhat} and \eqref{eq:fsrS} can be reformulated as  \eqref{eq:pcs} with 
\begin{equation}
	\gamma = 0\quand\widehat{N} = \begin{pmatrix}
		\alpha_1\\
		\frac{g}{2}[h+b]
	\end{pmatrix}. 
\end{equation}
\end{THM}

The proof of Theorem \ref{thm-XS} is given in Appendix \ref{app-XS}. 

\section{Convergence analysis}\label{sec:ca}
\setcounter{equation}{0}

In this section, we analyze the limit of the numerical solution of the generic scheme \eqref{eq:pcs} as $\dx, \dt \to 0$. With an abuse of notation, this limit will be denoted by $U = (h, m)^T$ in this section (note it may be different from the exact weak solution in Section \ref{sec:ws}). In Section \ref{sec:glw}, we prove a Lax--Wendroff type theorem on two sides of the domain separated by the bottom discontinuity (i.e., $x=0$), which characterizes the solution limit in the weak formulation. In Section \ref{sec:pwlimit}, we study the point limits at $x = 0$ to explain the formation of the spurious spike, which is caused by the numerical viscosity. The analysis indicates that removing the numerical viscosity may help eliminate the artifact. In Section \ref{sec:cLxF}, we present and analyze a so-called cLxF scheme, which removes the numerical viscosity at $x = 0$ by replacing the LxF flux with the central flux. 

\subsection{A Lax--Wendroff type theorem}\label{sec:glw}

In this subsection, we use $\cc$ for a generic constant, whose value may vary at different places.

Recall that $I_j = [x_{j-\hf},x_{j+\hf})$ and $I^n = [t^n,t^{n+1})$. For a given set $A$, we denote $\chi_A$ the characteristic function on $A$. Then the discrete solution $(U_j^n)$ can be interpreted as a piecewise constant function $U_\Delta(x,t)$ on $(-\infty,\infty)\times [0,\infty)$
\begin{equation}
U_\Delta(x,t) = \sum_{n = 0}^\infty\sum_{j = -\infty}^\infty U_j^n \chi_{I_j}(x)\chi_{I^n}(t). 
\end{equation}

To state our convergence theorem, we need to make several assumptions on the numerical flux $\widehat{F}$, the numerical solution $U_\Delta$, and the limit of $U_\Delta$. In below, Assumption \ref{assp:LW} collects standard assumptions of Lax--Wendroff theorems that are also required for classical homogeneous conservation laws; Assumption \ref{assp:spike} states additional assumptions needed for the convergence analysis at $x = 0$. 
\begin{ASSP}\label{assp:LW}\;
	\begin{enumerate}
		\item $\widehat{F}(\cdot,\cdot)$ is consistent and Lipschitz continuous, i.e.
		\begin{enumerate}
			\item (consistency) 
			\begin{equation}
			\widehat{F}(V,V) = F(V) \qquad \forall\, V;
			\end{equation}
			\item (Lipschitz continuity)
			\begin{equation}\label{eq:Lipschitz}
			\left|\widehat{F}(V,W) - \widehat{F}(\overline{V},\overline{W})\right| \leq \cc \left(\left|V - \overline{V}\right|+ \left|W - \overline{W}\right|\right)\qquad \forall \, V, \overline{V}, W, \overline{W}.
			\end{equation}
		\end{enumerate}
		\item For all $\Delta x$ and $\Delta t$, the function value and the total variation of $U_\Delta$ are uniformly bounded.
		\begin{enumerate}
			\item (uniform boundedness)
			\begin{equation}
				\sup_{\dx_j, \dt^n} \nm{U_\Delta}_{L^\infty((-\infty,\infty)\times[0,\infty))} \leq \cc\qquad \text{for some constant }\cc;
			\end{equation}
			\item (uniformly bounded total variation)
			\begin{equation}
				\sup_{\dx_j, \dt^n, t} \mathrm{TV}(U_\Delta(\cdot,t)) \leq \cc\qquad \text{for some constant }\cc.
			\end{equation}
			Here $\mathrm{TV}$ denotes the total variation function $TV(V) = \sup\sum_{j = 1}^N \left|V(\xi_j) - V(\xi_{j-1})\right|$ with the supremum taken over all subdivisions of the real line $-\infty = \xi_0 <\xi_1 <\cdots < \xi_N = \infty$. 
		\end{enumerate}
	\item $U_\Delta$ converges to a function $U$ in $L_{\mathrm{loc}}^1$:
	\begin{equation} \label{eq5.6}
	\lim_{\Delta x, \Delta t \to 0} \|U - U_\Delta\|_{L^1_{\mathrm{loc}}((-\infty,\infty)\times[0,\infty))} = 0.
	\end{equation}
	\end{enumerate}
\end{ASSP}
\begin{ASSP}\label{assp:spike}\;
	\begin{enumerate}
		\item $\widehat{N}(\cdot,\cdot)$ is Lipschitz continuous in the sense of  \eqref{eq:Lipschitz}. 
		\item Let
		\begin{equation}
			U_\Delta^\pm(t) = \lim_{x\to 0^\pm} U_\Delta (x,t).
		\end{equation} 
		$U_\Delta^\pm(t)$ converge to some functions $U_\star^\pm(t)$ in $L_{\mathrm{loc}}^1$:
		\begin{equation}
			\lim_{\Delta t \to 0} \|U_\star^\pm - U_\Delta^\pm\|_{L^1_{\mathrm{loc}}([0,\infty))} = 0.
		\end{equation} 
	\end{enumerate}
\end{ASSP}

\begin{REM}
	Note it may yield $U_\star^\pm(t) \neq \lim_{x\to 0^\pm} U(x,t)$. In other words, the limit of traces $U_\Delta^\pm(t)$ may not be the traces of the limit $U(0^\pm,t)$. Let $U_\star^\pm(t) = \left(h_\star^\pm(t), m_\star^\pm(t)\right)$ and $U(0^\pm,t) = (h_\hf^\pm(t), m_\hf^\pm(t))$. One can interpret $m_\star^\pm$ as the peak of the spurious spike and $m^\pm$ as the foot of the spike in Figure \ref{fig:db-nwb}.
\end{REM}
Next we introduce a few more notations. Let $I_\Psi\times I^\Psi$ be a rectangle on the spacetime domain that contains the support of the test function $\Psi \in C_c^1((-\infty,\infty)\times[0,\infty))$. Furthermore,  we define
\begin{equation}\label{eq-psinotations}
	\Psi_\hf(t) = \Psi\left(x_\hf,t\right),\quad  \Psi_{\Delta,\hf}(t) = \sum_{n = 0}^\infty \Psi_\hf(t^n)\chi_{I^n}(t),  \quand \Psi^n(x) = \Psi(x,t^n).
\end{equation}
Similar notations will also be used for other variables including $U$ and $[b]$. 

Now we are ready to state a Lax--Wendroff type theorem, which gives the weak forms of the solution limit on the first and the second quadrants of the spacetime plane, respectively. These weak forms will facilitate our analysis in Subsection \ref{sec:pwlimit} to explain the numerical artifact at the bottom discontinuity. It should be mentioned that in the classical Lax--Wendroff theorem, one usually considers the entire upper half plane, as is stated in Corollary \ref{cor-lw}. Here the argument for separately considering the solution limit on both sides of the bottom discontinuity shares some similar flavors as that in the deriviation of the Rankine--Hugoniot jump condition. 

\begin{THM}\label{thm:lw}
	Suppose $b(x) \in C^1((-\infty,0)\cup(0,\infty))$ has a jump discontinuity at $x_\hf = 0$. 
	Then under Assumptions \ref{assp:LW} and \ref{assp:spike} and for any $\Psi = (\eta,\theta)^T\in C_c^1\left((-\infty,\infty)\times [0,\infty)\right)$, the limit solution $U$ in \eqref{eq5.6} satisfies  %\zs{check how to define $U^0$}
	\begin{equation}\label{eq:lwleft}
		\begin{aligned}
			&-\int_0^\infty \int_{-\infty}^0 U\cdot \Psi _t \dd x \dd t - \int_{-\infty}^0 U^0\cdot\Psi^0 \dd x
			- \int_0^\infty\int_{-\infty}^0 F(U)\cdot \Psi_x \dd x\dd t + \int_0^\infty \wh{F}\left(U_\star^-,U_\star^+\right)\cdot \Psi_\hf \dd t\\
			&\quad =
			-\int_0^\infty\int_{-\infty}^0g h b_x \theta \dd x \dd t - \int_0^\infty\frac{g}{2}\check{h}_\gamma\left(h_\star^-,h_\star^+\right)[b]_\hf\theta_\hf \dd t + \int_0^\infty \widehat{M}\left(U_\star^-,U_\star^+\right)\cdot \Psi_\hf \dd t,
		\end{aligned}
	\end{equation}
	and
	\begin{equation}\label{eq:lwright}
		\begin{aligned}
			&-\int_0^\infty \int_0^\infty U\cdot \Psi _t \dd x \dd t - \int_{0}^\infty U^0\cdot\Psi^0 \dd x
			- \int_0^\infty\int_0^\infty F(U)\cdot \Psi_x \dd x\dd t - \int_0^\infty \wh{F}\left(U_\star^-,U_\star^+\right)\cdot \Psi_\hf \dd t\\
			&\quad =
			-\int_0^\infty\int_0^\infty g h b_x \theta \dd x \dd t - \int_0^\infty\frac{g}{2}\check{h}_\gamma\left(h_\star^-,h_\star^+\right)[b]_\hf\theta_\hf \dd t - \int_0^\infty \widehat{M}\left(U_\star^-,U_\star^+\right)\cdot \Psi_\hf \dd t.  
		\end{aligned}
	\end{equation}
\end{THM}

Before proving Theorem \ref{thm:lw}, let us remark that one can combine two equations \eqref{eq:lwleft} and \eqref{eq:lwright} to obtain a Lax--Wendroff theorem on the entire domain, as stated in Corollary \ref{cor-lw}. This corollary is a special case of \cite[Claim 1]{castro2008many} for the SWEs \eqref{eq-ncsv}. The additional assumption $\hbn{h_\star + b}_\gamma = \hbn{h + b}_\gamma$ is related with the assumption of ``convergence in the sense of graph" in \cite[Claim 1]{castro2008many}. 

\begin{COR}\label{cor-lw}
	Under the assumption of Theorem \ref{thm:lw}, we have  
	\begin{equation}\label{eq-lwsum}
		\begin{aligned}
			&-\int_0^\infty \int_{-\infty}^\infty U\cdot \Psi _t \dd x \dd t - \int_{-\infty}^\infty U^0\cdot\Psi^0 \dd x
			- \int_{0}^\infty\int_{-\infty}^\infty F(U)\cdot \Psi_x \dd x\dd t\\
			&\qquad =
			-\int_0^\infty\int_{-\infty}^\infty g h b_x \theta \dd x \dd t - \int_0^\infty g\check{h}_\gamma\left(h_\star^-,h_\star^+\right)[b]_\hf\theta_\hf \dd t.  
		\end{aligned}
	\end{equation}
	Furthermore, if $\hbn{h_\star + b}_\gamma = \hbn{h + b}_\gamma$, then the limit of the numerical solution given by \eqref{eq:pcs} and \eqref{eq:FSM} is a weak solution of \eqref{eq-swe}. 
\end{COR}

The rest of this subsection is dedicated to the proof of Theorem \ref{thm:lw}. 
\begin{PROP}\label{prop:sbp} Let $\left(w_{j+\hf}\right)$ and $\left(v_j\right)$ be two sequences with only finite nonzero elements in $(v_j)$. Then we have 
	\begin{equation}
		\sum_{j  = - \infty}^0 \left(w_{j+\hf}-w_{j-\hf}\right)v_j = - \sum_{j = -\infty}^{-1} w_{j+\hf}\left(v_{j+1}-v_j\right) + w_\hf v_0. \label{eq:sbp-}
	\end{equation}
\end{PROP}

\begin{LEM}\label{lem:tracelim}
	Suppose $\widehat{Q}(\cdot,\cdot)$ is Lipschitz continuous. Then for any $\Psi \in C_c^1((-\infty,\infty)\times [0,\infty))$,
	\begin{equation}
		\lim_{\Delta t\to 0} \sum_{n = 0}^\infty \widehat{Q}\left(U_0^n,U_1^n\right)\cdot\Psi_\hf^n \Delta t^n = \int_0^\infty \widehat{Q}\left(U_\star^-,U_\star^+\right)\cdot \Psi_\hf \dd t.
	\end{equation}
\end{LEM}
\begin{proof}
	Using the triangle inequality, we have 
	\begin{equation}\label{eq-T1T2}
	\left| \int_0^\infty \widehat{Q}\left(U_\star^-,U_\star^+\right)\cdot \Psi_\hf \dd t - \sum_{n = 0}^\infty \widehat{Q}\left(U_0^n,U_1^n\right)\cdot\Psi_\hf^n \Delta t^n  \right|\leq T_1 + T_2,
	\end{equation}
	where 
	\begin{align}
		T_1 =&  \left|\int_0^\infty \widehat{Q}\left(U_\star^-,U_\star^+\right)\cdot\Psi_\hf \dd t - \int_0^\infty  \widehat{Q}\left(U_\Delta^-,U_\Delta^+\right)\cdot\Psi_\hf \dd t\right|,\\
		T_2 =& \left|\int_0^\infty  \widehat{Q}\left(U_\Delta^-,U_\Delta^+\right)\cdot\Psi_\hf \dd t - \sum_{n = 0}^\infty \widehat{Q}\left(U_0^n,U_1^n\right)\cdot\Psi_\hf^n \Delta t^n \right|.\label{eq-T2}
	\end{align}
	For $T_1$, recall that $\widehat{Q}(\cdot,\cdot)$ is Lipschitz continuous, $\Psi$ is bounded, and $U_\Delta^\pm \to U_\star^\pm$ in $L_{\mathrm{loc}}^1$. Therefore,
	\begin{equation}\label{eq-T2est}
	\begin{aligned}
	T_1	\leq \cc\int_{I^\Psi}  \left|U_\Delta^- - U_\star^-\right| + \left|U_\Delta^+ - U_\star^+\right| \dd t \to 0, \qquad \text{ as } \Delta t \to 0.
	\end{aligned}
	\end{equation}
	For $T_2$, we use \eqref{eq-psinotations} and the fact that $U_0^n=U_\Delta^-$, $U_1^n=U_\Delta^+$ on a fixed mesh to obtain
	\begin{equation}\label{eq-x0sumt}
	\sum_{n = 0}^\infty \widehat{Q}\left(U_0^n,U_1^n\right)\cdot\Psi_\hf^n \Delta t^n = \int_0^\infty \widehat{Q}\left(U_\Delta^-,U_\Delta^+\right)\cdot \Psi_{\Delta, \hf} \dd t.
	\end{equation}
	As a two-variable function $\widehat{Q}(\cdot,\cdot)$, there must be some input $V, W$ independent of $\Delta x, \Delta t$, such that $|\widehat{Q}(V,W)|$ is finite (for example, consider $|\widehat{Q}(0,0)|$ if it is well-defined). Therefore, using the Lipschitz continuity of $\widehat{Q}\left(\cdot, \cdot\right)$, and the uniform boundedness of $U_\Delta$, we have 
	\begin{equation}\label{eq:boundedQ}
		\left|\widehat{Q}\left(U_\Delta^-,U_\Delta^+\right)\right|  \leq \left|\widehat{Q}\left(U_\Delta^-,U_\Delta^+\right) - \widehat{Q}(V,W)\right| + \left|\widehat{Q}(V,W)\right| \leq \cc\left(|U_\Delta^--V|+|U_\Delta^+-W|\right)  + \left|\widehat{Q}(V,W)\right|\leq \cc.
	\end{equation}
	 Therefore, $\nm{\widehat{Q}\left(U_\Delta^-,U_\Delta^+\right)}_{L^\infty}\leq \cc$ is uniformly bounded in terms of $\Delta x,\Delta t$. Using this fact and substituting \eqref{eq-x0sumt} into \eqref{eq-T2} yield
	\begin{equation}\label{eq-T1est}
	\begin{aligned}
	T_2=&\; \left|\int_0^\infty \widehat{Q}\left(U_\Delta^-,U_\Delta^+\right) \cdot\left( \Psi_{\Delta,\hf} - \Psi_\hf\right) \dd t\right| 
		\leq \cc\int_{I^\Psi} \left|\Psi_{\Delta,\hf} - \Psi_\hf\right|\dd t \\
	\leq&\;  \cc \Delta t \sup_{t}\left|\Psi_t\left(x_{\hf}, t\right)\right| \left|I^\Psi\right| \to 0, \quad \text{ as } \Delta t \to 0.
	\end{aligned}
	\end{equation}
	 The proof is completed after substituting \eqref{eq-T1est} and \eqref{eq-T2est} into \eqref{eq-T1T2}.
\end{proof}

\begin{proof}[Proof of Theorem \ref{thm:lw}]
	We will only prove \eqref{eq:lwleft} and one can follow similar lines to deduce \eqref{eq:lwright}. 
	Taking a dot product of \eqref{eq:pcs} with $\Psi_j^n\Delta x_j\Delta t^n$ and summing over $-\infty \leq j\leq 0$ and $0 \leq n \leq \infty$ lead to
	\begin{equation}\label{eq:lwsum}
	\begin{aligned}
	&\sum_{n = 0}^{\infty}\sum_{j = -\infty}^0 \left(U_{j}^{n+1} - U_{j}^n\right)\cdot \Psi_j^n\Delta x_j + \sum_{n = 0}^{\infty}\sum_{j = -\infty}^0 \left(\wh{F}_{j+\hf}^n - \wh{F}_{j-\hf}^n\right)\cdot \Psi_j^n \Delta t^n\\
	 =& \sum_{n = 0}^{\infty}\sum_{j = -\infty}^0 \left(\wh{S}_{j+\hf}^n + \wh{S}_{j-\hf}^n\right)\cdot \Psi_j^n\Delta t^n+ \sum_{n = 0}^{\infty}\sum_{j = -\infty}^0 \left(\wh{M}_{j+\hf}^n - \wh{M}_{j-\hf}^n\right)\cdot \Psi_j^n \Delta t^n.
	\end{aligned}
	\end{equation}
	We will then take summation by parts in Proposition \ref{prop:sbp} and send $\dx, \dt \to 0$. 
	
	For the left hand side of \eqref{eq:lwsum}, one only needs to pay attention to the right boundary $x_\hf = 0$. We can use Assumption \ref{assp:LW} and follow similar lines as those in \cite[Theorem 2.3]{shi2018local} to derive
	\begin{equation}\label{eq:LW-U}
	\begin{aligned}
	\sum_{n = 0}^{\infty}\sum_{j = -\infty}^0 \left(U_{j}^{n+1} - U_{j}^n\right)\cdot \Psi_j^n\Delta x_j =& - \sum_{n = 0}^{\infty}\sum_{j = -\infty}^0 U_{j}^n \cdot \left(\Psi_j^n - \Psi_j^{n-1}\right)\Delta x_j - \sum_{j = -\infty}^0 U_{j}^0\Psi_j^0\Delta x_j\\
	\to& -\int_{0}^\infty \int_{-\infty}^0 U \cdot \Psi_t \dd x \dd t	- \int_{-\infty}^0 U^0\cdot \Psi^0 \dd x,
	\end{aligned}
	\end{equation}
	\begin{equation}\label{eq:LW-F}
	\begin{aligned}
	\sum_{n = 0}^{\infty}\sum_{j = -\infty}^0 \left(\wh{F}_{j+\hf}^n - \wh{F}_{j-\hf}^n\right)\cdot \Psi_j^n\Delta t^n =&  - \sum_{n = 0}^{\infty}\sum_{j = -\infty}^{-1} \wh{F}_{j+\hf}^n\cdot \left(\Psi_{j+1}^n -\Psi_j^n\right)\Delta t^n  + \sum_{n = 0}^{\infty}\wh{F}_{\hf}^n\cdot \Psi_\hf^n\Delta t^n\\
	\to& -\int_0^\infty \int_{-\infty}^0 F\cdot \Psi_x \dd x\dd t + \int_0^\infty \wh{F}(U_\star^-,U_\star^+) \cdot \Psi_\hf \dd t,
	\end{aligned}
	\end{equation}
	as $\Delta x$, $\Delta t \to 0$, where Lemma \ref{lem:tracelim} was used to show the convergence of the second term in \eqref{eq:LW-F}. 
	
	For the right hand side of \eqref{eq:lwsum}, one can apply Proposition \ref{prop:sbp} to get
	\begin{equation}
	\begin{aligned}
	&\sum_{n = 0}^{\infty}\sum_{j = -\infty}^0 \left(\widehat{M}^n_{j+\hf} - \widehat{M}^n_{j-\hf}\right)\cdot \Psi_j^n\Delta t^n \\
	=&  - \sum_{n = 0}^{\infty}\sum_{j = -\infty}^{-1} \widehat{M}^n_{j+\hf}\cdot \left(\Psi_{j+1}^n -\Psi_j^n\right)\Delta t^n  + \sum_{n = 0}^{\infty}\widehat{M}_{\hf}^n\cdot \Psi_\hf^n\Delta t^n:= T_1 + T_2.
	\end{aligned}
	\end{equation}
	Note that $U_j^n$ is uniformly bounded and $\widehat{N}(\cdot,\cdot)$ is Lipschitz continuous. Using the same argument as that in \eqref{eq:boundedQ}), one can show that $\widehat{N}_{j+\hf}^n$ is bounded. Since $b$ is differentiable on $(-\infty,0)$, we have $\left|\widehat{M}^n_{j+\hf}\right| = \left|\hf \widehat{N}_{j+\hf}^n[b]_{j+\hf}\right|\leq \cc\nm{b_x}_{L^\infty} \Delta x_j$. Also note that $\left|\Psi_{j+1}^n - \Psi_j^n\right|\leq  \nm{\Psi_x}_{L^\infty} \Delta x_j$. Hence we have 
	\begin{equation}\label{eq:MT1}
	T_1 \leq \cc \nm{b_x}_{L^\infty}\nm{\Psi_x}_{L^\infty} \sum_{n = 0}^\infty \sum_{j = -\infty}^{-1} \chi_{\{j:x_{j+1} \in I_\Psi\}}\Delta x_j^2 \Delta t^n \to 0,  \quad \text{ as } \Delta x, \Delta t \to 0.
	\end{equation} 
	Consider $b$ to be given. $\widehat{N}(\cdot,\cdot)$ is Lipschitz implies that $\widehat{M}(\cdot,\cdot)$ is Lipschitz. The limit of $T_2$ can be obtained by applying Lemma \ref{lem:tracelim} with $\wh{Q} = \wh{M}$. Combining with the bound of $T_1$ in \eqref{eq:MT1}, we have
	\begin{equation}\label{eq:LW-M}
	\sum_{n = 0}^{\infty}\sum_{j = -\infty}^0 \left(\widehat{M}_{j+\hf} ^n- \widehat{M}_{j-\hf}^n\right)\cdot \Psi_j^n\Delta t^n \to \int_0^\infty \widehat{M}(U_\star^-,U_\star^+)\cdot \Psi_\hf dt,  \quad \text{ as } \Delta x, \Delta t \to 0.
	\end{equation}
	For the summation of $\wh{S}$ terms, it can be shown that
	\begin{equation}
	\begin{aligned}
	&\sum_{n = 0}^{\infty}\sum_{j = -\infty}^0 \left(\wh{S}_{j+\hf}^n + \wh{S}_{j-\hf}^n\right)\cdot \Psi_j^n\Delta t^n \\
	=& 
	-\frac{g}{2}\sum_{n = 0}^{\infty}\sum_{j = -\infty}^0\left(\check{h}_{\gamma,j+\hf}^n[b]_{j+\hf} + \check{h}_{\gamma,j-\hf}^n[b]_{j-\hf}\right) \theta_j^n\Delta t^n \\
	=& -\frac{g}{2}\sum_{n = 0}^{\infty}\sum_{j = -\infty}^{-1} \check{h}_{\gamma, j+\hf}^n[b]_{j+\hf} (\theta^n_j+\theta^n_{j+1}) \Delta t^n -\frac{g}{2}\sum_{n = 0}^{\infty}{\check{h}_{\gamma,\hf}^n}[b]_\hf \theta_0^n\Delta t^n.
	\end{aligned}
	\end{equation}
	Note that $\theta$ is compactly supported, $\theta_j + \theta_{j+1} = 2\theta_{j+\hf} + o(1)$, $\theta_0^n = \theta_\hf^n + o(1)$, and %$|\theta_j + \theta_{j+1} - 2\theta_{j+\hf}|\leq c \Delta x$, and 
	\begin{equation}
		[b]_{j+\hf} = \frac{b_{j+1}-b_j}{(\Delta x_j + \Delta x_{j+1})/2} \cdot \frac{\Delta x_j + \Delta x_{j+1}}{2}= \left((b_x)_{j+\hf}+ o(1)\right) \frac{\Delta x_j + \Delta x_{j+1}}{2} .
	\end{equation} One can deduce that 
	\begin{equation}\label{eq:LW-S}
		\begin{aligned}
				&\sum_{n = 0}^{\infty}\sum_{j = -\infty}^0 \left(\wh{S}_{j+\hf}^n + \wh{S}_{j-\hf}^n\right)\cdot \Psi_j^n\Delta t^n \\
		%	=& -g\sum_{n = 0}^{\infty}\sum_{j = -\infty}^{-1} \check{h}_{\gamma, j+\hf}^n\frac{b_{j+1}-b_j}{\Delta x_j}\theta_{j+\hf}^n\Delta x_j \Delta t^n  -\frac{g}{2}\sum_{n = 0}^{\infty}{\check{h}_{\gamma,\hf}^n}[b]_\hf \theta_0^n\Delta t^n + o(1)\\
	=& -g\sum_{n = 0}^{\infty}\sum_{j = -\infty}^{-1} \check{h}_{\gamma,j+\hf}^n(b_x)_{j+\hf}\theta_{j+\hf}^n\frac{\Delta x_j+\Delta x_{j+1}}{2} \Delta t^n  -\frac{g}{2}\sum_{n = 0}^{\infty}{\check{h}_{\gamma,\hf}^n}[b]_\hf \theta_\hf^n\Delta t^n + o(1)\\
	\to & -\int_0^\infty \int_{-\infty}^0 g h b_x \theta \dd x \dd t- \int_0^\infty \frac{g}{2} \check{h}_\gamma\left(h_\star^-,h_\star^+\right) [b]_\hf \theta_\hf \dd t,  \qquad \text{ as } \Delta t, \Delta x \to 0.
	\end{aligned}
	\end{equation}
	Here the convergence of the first term can be obtained through the standard estimate, and the convergence of the second term can be shown by applying Lemma \ref{lem:tracelim} with $\widehat{Q} = \left(0, \check{h}\right)^T[b]$.
	
	The equation \eqref{eq:lwleft} follows after substituting \eqref{eq:LW-U}, \eqref{eq:LW-F}, \eqref{eq:LW-M}, and \eqref{eq:LW-S} into \eqref{eq:lwsum}, and the proof is completed.
\end{proof}

\subsection{Solution limit at the discontinuous bottom}\label{sec:pwlimit}

In this subsection, we provide a characterization of $U_\star^\pm$, the point limits of $U$ at the bottom discontinuity $x = 0$. We assume that the spatial support of $\Psi$ is concentrated near $x = 0$ and $U \in C_c^1$ except along $x = 0$ on the region of interest. We denote 
\begin{equation}
	U^\pm(t) = \lim_{x\to 0^\pm} U(x,t)= \lim_{x\to 0^\pm}\lim_{\Delta x, \Delta t\to 0} U_\Delta(x,t),
\end{equation} 
for the limit of $U$ at $x = 0$. Note it should be distinguished from the limit of the numerical solution \begin{equation}
	U_\star^\pm(t) = \lim_{\Delta x, \Delta t\to 0} U_\Delta^\pm(t) = \lim_{\Delta x, \Delta t\to 0}\lim_{x\to 0^\pm} U_\Delta(x,t).
\end{equation} 
We use the shorthand notations $\wh{F}_\star = \wh{F}\left(U_\star^-,U_\star^+\right)$. Furthermore, we denote 
\begin{equation}
\{F\}_\star = \hf \left(F\left(U_\star^-\right)+F\left(U_\star^+\right)\right), \quad \{m\}_\star = \hf\left(m_\star^- + m_\star^+\right) \quand [m]_\star = m_\star^+ - m_\star^-. 
\end{equation} 
Similar notations will also be used for other unknowns, including $\wh{A}_\star$, $\alpha_{1,\star}$, $\wh{M}_\star$, $\wh{N}_\star$, and $\check{h}_{\gamma,\star}$. 

Applying integration by parts to \eqref{eq:lwleft} yields
\begin{equation}
\begin{aligned}
&\int_0^\infty\int_{-\infty}^0 \left(U_t + F(U)_x\right)\cdot \Psi \dd x \dd t + \int_0^\infty \left(\wh{F}_\star-F\left(U_\hf^-\right)\right)\cdot \Psi_\hf \dd t \\
=& -\int_0^\infty\int_{-\infty}^0ghb_x \theta \dd x \dd t - \int_0^\infty \frac{g}{2}\check{h}_{\gamma,\star}[b]_{\hf}\theta_\hf \dd t  + \int_0^\infty \widehat{M}_\star\cdot \Psi_\hf \dd t. 
\end{aligned}
\end{equation}
Since the equation holds strongly in the interior, we have 
\begin{equation}
U_t + F(U)_x = \begin{pmatrix}
0\\
-ghb_x
\end{pmatrix}.
\end{equation}
As a result, one can get
\begin{equation}
\begin{aligned}
\int_0^\infty \left(\left(\wh{F}_\star - \widehat{M}_\star \right) -F\left(U_\hf^-\right) \right)\cdot \Psi_\hf dt =-  \int_0^\infty \frac{g}{2}\check{h}_{\gamma,\star}[b]_\hf\theta_\hf
dt. 
\end{aligned}
\end{equation}
The arbitrariness of $\Psi$ indicates that
\begin{equation}\label{eq:ibplflux}
\left(\wh{F}_\star- \widehat{M}_\star\right)-F\left(U_\hf^-\right) = \begin{pmatrix}
0\\
- \frac{g}{2}\check{h}_{\gamma,\star}[b]_\hf
\end{pmatrix}.
\end{equation}
Similarly, applying integration by parts to \eqref{eq:lwright} would yield
\begin{equation}\label{eq:ibprflux}
F\left(U_\hf^+\right) - \left(\wh{F}_\star - \widehat{M}_\star\right) = \begin{pmatrix}
0\\
- \frac{g}{2}\check{h}_{\gamma,\star}[b]_\hf
\end{pmatrix}.
\end{equation}
We can add \eqref{eq:ibplflux} to \eqref{eq:ibprflux} to obtain (recall the notations in \eqref{eq-hb})
\begin{equation}\label{eq:add}
[F]_\hf = \begin{pmatrix}
0\\
- g\check{h}_{\gamma,\star}[b]_\hf
\end{pmatrix}\Rightarrow
\begin{pmatrix}
[\mm]_\hf\\
\left[\frac{\mm^2}{h}+\frac{g}{2}h^2\right]_\hf
\end{pmatrix} = \begin{pmatrix}
0\\
- g\check{h}_{\gamma,\star}[b]_\hf
\end{pmatrix}.
\end{equation}
Note that \eqref{eq:add} can be interpreted as a Rankine--Hugoniot jump condition along $x = 0$ for the limit of the numerical solution in \eqref{eq:lwleft} and \eqref{eq:lwright}. 

The first line in \eqref{eq:add} implies the following theorem, which roughly says that the solution limit of $m$ is continuous at the bottom discontinuity. 
\begin{THM}\label{thm:flatm}
	$\mm^- = \mm^+$ at $x = 0$. 
\end{THM}
We now try to explain the spurious spike formed in $m$ at $x_\hf = 0$. 
Subtracting \eqref{eq:ibprflux} from \eqref{eq:ibplflux} and then dividing by $2$ yield
\begin{equation}\label{eq:sub}
\wh{F}_\star -\widehat{M}_\star= \{F\}_\hf.
\end{equation}
Substituting in the definitions of $\wh{F}$ in \eqref{eq:FSM} and $\wh{M}$ in \eqref{eq:lxfM} leads to
\begin{equation}
\{F\}_\star - \{F\}_\hf = \hf A_{\star}[U]_{\star}  + \frac12\widehat{N}_\star [b]_\hf,
\end{equation}
where
	\begin{equation}
		A_\star = A(U_\star^-,U_\star^+) = \begin{pmatrix}
			\alpha_{1,\star}&0\\0&\alpha_{2,\star}
		\end{pmatrix}.
\end{equation}
Or equivalently, it can be written as
\begin{equation}
\begin{aligned}
\begin{pmatrix}
\{\mm\}_\star - \{\mm\}_\hf\\
\{\frac{\mm^2}{h}\}_\star - \{\frac{\mm^2}{h}\}_\hf + \frac{g}{2}\left(\{h^2\}_\star - \{h^2\}_\hf\right)
\end{pmatrix} = 
\frac{1}{2} \begin{pmatrix}
\alpha_{1,\star}[h]_\star\\
\alpha_{2,\star}[\mm]_\star
\end{pmatrix} + \frac12\widehat{N}_\star [b]_\hf.
\end{aligned}
\end{equation}
Note that $\{\mm\}_\hf=\mm_\hf$ from Theorem \ref{thm:flatm}, and the first component of $\wh{N}_\star$ is $0$ for the LxF scheme and is $ \alpha_{1,\star}$ for the well-balanced schemes. Hence we can obtain the following theorem. 
\begin{THM}\label{thm:spike}
	\;
	\begin{enumerate}
		\item $\{\mm\}_\star - \mm_\hf = \hf \alpha_{1,\star}[h]_\star$ for the LxF scheme \eqref{eq:LxF};
		\item $\{\mm\}_\star - \mm_\hf = \hf \alpha_{1,\star}[h+b]_\star$ for the well-balanced schemes, including the wbLxF scheme \eqref{eq:wbLxF}, the HR scheme \eqref{eq:hr}, and the XS scheme \eqref{eq:fsr}.
	\end{enumerate}
	In other words, the averaged height of the spurious spike in the momentum $\mm$ is proportional to the viscosity constant $\alpha_{1,\star}$ and the jump of $h_\star$ or $h_\star + b$. 
\end{THM}
\begin{REM}\label{rmk:wbnwbspike}
 	From Theorem \ref{thm:spike}, we see that for a fixed value of $\alpha_{1,\star}$, the height of the spurious spike is proportional to $[h]$ for the LxF scheme and is proportional to $[h+b]$ for the well-balanced schemes. For the dam-break problem \eqref{eq:db}, we actually have $[h]>[h+b]$, as seen in Figure \ref{fig:db-nwb}. Therefore, well-balanced schemes form a shorter spike and suffer less from the wrong convergence compared with non-well-balanced schemes, as observed in Section \ref{sec:artifact}. In general, for subcritical Riemann problems, we usually have the water to be shallower on the higher side of the step, hence $[h]>[h+b]$. For these problems, well-balanced schemes typically perform better than non-well-balanced schemes, although the spurious spike still appears following the result in Theorem \ref{thm:spike}.
\end{REM}

\subsection{Avoiding numerical artifacts by removing viscosity}\label{sec:cLxF}

In this subsection, we investigate a possible approach to remove the numerical artifact. We focus on the modification of the LxF scheme \eqref{eq:LxF}, and briefly summarize the generalization to well-balanced methods in Remark \ref{remark5.4}.

From Theorem \ref{thm:spike}, we see that the averaged height of the spurious spike $\{\mm\}_\star - \mm_\hf$ is zero if and only if $\alpha_{1,\star}[h]_\star = 0$. Since the exact solution $h$ may have a discontinuity at $x = 0$, in general, we cannot expect $[h]_\star$ to be zero. The remaining choice is to set $\alpha_{1,\star}= \alpha_1(U_\star^-,U_\star^+) = 0$, which essentially depends on $\alpha_{1,\hf}$ in the computation. This motivates us to introduce the following scheme. 

\textbf{Central--Lax--Friedrichs scheme (cLxF scheme).} The scheme is similar to the original LxF scheme, except for using the central flux in the definition of the numerical flux $\widehat{m}_\hf$. In other words, we change the definition of $\alpha_{1,\hf}=\alpha_\hf$ in \eqref{eq:alpha} to $\alpha_{1,\hf}= 0$, then the new scheme is given by \eqref{eq:pcs}, \eqref{eq:FSM}, and
%\zss{\begin{align}\label{eq:malpha}
%A_{j+\hf}^n = \left\{
%\begin{array}{ll}
%\diag(\alpha_{j+\hf}^n,\alpha_{j+\hf}^n),& j\neq 0\\
%\diag(0,\alpha_{j+\hf}^n),& j = 0
%\end{array}\right. \quad \text{with} \quad \alpha_{j+\hf}^n = \max_{l = j, j+1}\left(\left|\frac{\mm_{l}^n}{h_{l}^n}\right| + \sqrt{gh_{l}^n}\right).
%\end{align}}
%\zss{
%\begin{equation}
%	\alpha_{1,\hf}^n = 0.
%\end{equation}}
\begin{equation}\label{eq:malpha}
	A_{\hf} = \begin{pmatrix}
		0&0\\0&\alpha_\hf
	\end{pmatrix} \quand
	A_{j+\hf} = \alpha_{j+\hf} \begin{pmatrix}
		1&0\\0&1
	\end{pmatrix}\text{ for } j\neq 0. 
%	\alpha_{1,j+\hf} = \left\{\begin{array}{ll}
%		0,&j = 0\\
%		\alpha_{j+\hf}, & j\neq 0
%%	\end{array}\right.\quad \text{and}\quad  \alpha_{2,j+\hf} = \alpha_{j+\hf}.
\end{equation}

Although numerically we see that the cLxF scheme converges to exact Riemann solutions in various tests, we are not able to mathematically prove the guaranteed convergence. However, we do obtain partial results to explain why the cLxF could avoid certain numerical artifacts. 

Theorem \ref{thm:cLxFm} is a direct consequence of Theorem \ref{thm:flatm} and Theorem \ref{thm:spike}. It implies that the cLxF scheme cannot form a one-sided spurious spike, for which $(\mm_\star - \mm)^\pm$ are of the same sign. 

\begin{THM}\label{thm:cLxFm}
	For the cLxF scheme, we have $\{\mm_\star\} = \mm^- = \mm^+$ at $x = 0$.
\end{THM}

\begin{REM}
	However, it is still possible for the cLxF scheme to generate a spurious spike by having $(\mm_\star - \mm)^\pm$ admit different signs. %See Example \ref{ex-coincide} which is a transonic test case. \YX{Is the new Ex 6.8 still this case? No?}
\end{REM}
Theorem \ref{thm:transpts} partially explains why the correct convergence of $\mm$ may lead to the correct convergence of $h$ in the cLxF scheme. Let $l_0$ be any fixed positive integer, it says that if $\mm_\Delta$ converges to a constant on $p\Delta x\leq 0 \leq q\Delta x$ (with fixed $p$, $q$), a shrinking neighborhood of 0, as $\Delta x \to 0$, then $h_\Delta$ will converge to a piecewise constant on the same region without any transition points. Roughly speaking, by avoiding the spurious spike in $m$, we may also avoid transition points of the limit function $h$ at the discontinuity of the bottom -- if one consider the first, the second, the third, $\cdots$ points on the left (right) of $x = 0$, they all converge to the same value, instead of spreading out along $x = 0$. 

\begin{THM}\label{thm:transpts}
	Let $j$, $n\geq 0$, $p\leq 0$ and $q\geq 1$ be fixed integers that are independent of $\Delta x$ and $\Delta t$. For the cLxF scheme, suppose
	\begin{equation}
	\lim_{\Delta x \to 0}U_\Delta(x_j,t^n)= 
	\begin{pmatrix}
	h_{\star,j}^n\\
	\mm_{\star,j}^n
	\end{pmatrix},
	\end{equation}
	and
	\begin{equation}\label{eq5.42}
		\mm_{\star,j}^n \equiv \mm_{\star,0}^n \quad \forall\, p\leq j\leq q
	\end{equation} 
	is a constant.  Then we have 
	\begin{equation}
		h_{\star,j}^n =
		\left\{\begin{array}{cc}
			h_{\star,0}^n& \forall\, p\leq j\leq 0\\
			h_{\star,1}^n& \forall\, 1\leq j\leq q
		\end{array} \right.
	\end{equation} 
	to be a piecewise constant. In other words, the limit function $h$ does not contain any transition points near the bottom discontinuity $x=0$. 
\end{THM}
\begin{proof}
	The first equation in \eqref{eq:LxF} is
	\begin{equation}\label{eq:scheme-eq1}
	\frac{h_{j}^{n+1} - h_{j}^{n} }{\Delta t^n} + \frac{\wh{\mm}_{j+\hf}^n - \wh{\mm}_{j-\hf}^n }{\Delta x_j} = 0. 
	\end{equation}	
	Let $\eta \in C_c^1((-\infty, \infty)\times[0,\infty))$. As before, we multiply \eqref{eq:scheme-eq1} with the test function $\eta_j^n\Delta t^n \Delta x_j$, and then sum from $n=0$ to $+\infty$ and from $j = -\infty$ to $l\leq 0$. It gives
	\begin{equation}
	\sum_{n = -\infty}^\infty \sum_{j = -\infty}^l \left(h_{j}^{n+1} - h_{j}^{n}\right)\eta_j^n \Delta x_j + \sum_{n = -\infty}^\infty \sum_{j = -\infty}^l\left(\wh{\mm}_{j+\hf}^n - \wh{\mm}_{j-\hf}^n \right)\eta_j^n\Delta t^n = 0 \quad \forall\, l \leq 0.
	\end{equation}
	Note here, the summation from $j = -\infty$ to $l$ corresponds to $x\in (-\infty, x_l)$ and we have $x_l \to 0$ as $\Delta x \to 0$. We take $\Delta x, \Delta t \to 0$ and use similar argument as those in Subsection \ref{sec:glw}. Then it yields
	\begin{equation}\label{eq:mstarl}
	- \int_{0}^\infty\int_{-\infty}^0 h \eta_t \dd x \dd t - \int_0^\infty \int_{-\infty}^0 \mm \eta_x \dd x \dd t + \int_0^\infty \wh{\mm}_{\star,l+\hf} \eta_{\hf}\dd t = 0 \quad \forall\, l\leq 0.
	\end{equation}
	Here 
	\begin{equation}
	\wh{\mm}_{\star,l+\hf} = \wh{\mm} \left(U_{\star,l}, U_{\star,l+1}\right)\quand U_{\star,l}(t) = 
	\begin{pmatrix}
	h_{\star,l}(t)\\
	m_{\star,l}(t)
	\end{pmatrix}
	= \lim_{\dx,\dt \to 0} U_{\Delta}(x_l,t). 
	\end{equation}
	Comparing \eqref{eq:mstarl} with the special case $l = 0$ in \eqref{eq:mstarl}, we obtain
	\begin{equation}
	\int_0^\infty \wh{\mm}_{\star,l+\hf} \eta_{\hf}\dd t = \int_0^\infty \wh{\mm}_{\star,\hf} \eta_{\hf}\dd t \quad \forall\, l\leq -1.
	\end{equation}
	Then by the arbitrariness of $\eta$, we have
	\begin{equation}
	\wh{\mm}_{\star,l+\hf} = \wh{\mm}_{\star,\hf} \quad \forall\, l\leq-1.
	\end{equation}
	Recall the definition of the numerical flux in \eqref{eq:lxfflux} and \eqref{eq:malpha}. It can be seen that
	\begin{equation}
	\frac{\mm_{\star,l} + \mm_{\star,l+1}}{2} - \frac{\alpha_{1,\star,l+\hf}}{2}\left(h_{\star,l+1} - h_{\star,l}\right) = \frac{\mm_{\star,0} + \mm_{\star,1}}{2} \quad \forall\, l \leq -1,
	\end{equation}
	where $\alpha_{1,\star,l+\hf} = \max_{k =l,l+1}\left(\left|\frac{m_{\star,k}}{h_{\star,k}}\right|+\sqrt{gh_{{\star,k}}}\right)$.

	According to the assumption \eqref{eq5.42} that $\mm_{\star,l}\equiv \mm_{\star,0}$ for all $p\leq l\leq q$, one can deduce that
	\begin{equation}
	\frac{\alpha_{1,\star,l+\hf}}{2}\left(h_{\star,l+1} - h_{\star,l}\right) = 0 \quad \forall\, p\leq l \leq -1.
	\end{equation}
	Note that $\alpha_{1,\star,l+\hf} \geq \max_{k=l,l+1} \sqrt{gh_{\star,k}}$ is uniformly positive for $l\leq-1$ (since we avoid the dry bed in the discussion). We have
	\begin{equation}
	h_{\star,l}=h_{\star,l+1}, \quad p\leq l\leq-1 \quad \Rightarrow \quad h_{\star,l} \equiv h_{\star,0}, \quad \forall\, p\leq l \leq 0.
	\end{equation}
	Using a similar argument with summation of $j$ from $1\leq l$ to $+\infty$, we can show that $h_{\star,l} \equiv h_{\star,1}$ for all $1\leq l\leq q$. 
\end{proof}

\begin{REM}\label{remark5.4}
			We have also examined the following modified numerical schemes:
			\begin{itemize}
				\item $\alpha_{1,\hf} = \alpha_{2,\hf} = 0$ in the LxF scheme (central fluxes for both equations). 
				\item $\alpha_{1,\hf} = 0$ in the wbLxF scheme (which implies $\wh{M}_\hf = 0$ from \eqref{eq:lxfM} and Theorem \ref{thm:wblxf}). 
				\item $\alpha_{1,\hf} = 0$, $\widehat{M}_\hf = 0$ in the HR scheme.
				\item $\alpha_{1,\hf} = 0$, $\widehat{M}_\hf = 0$ in the XS scheme.
			\end{itemize}
			For the first three schemes, we observe similar numerical results as those of the cLxF scheme for nontransonic tests in Examples \ref{ex-1s2r}--\ref{ex-1r2r-supern}. For the last scheme, it may converge to a different solution since we have $\gamma = 0$ in the numerical scheme while the reference exact solution is computed with $\gamma = \sgnb$. But other than that, the last method performs as expected for nontransonic tests and avoids the formation of the spurious spike. %For transonic tests, all these methods may not converge to the correct solution, similar to the cLxF scheme. Detailed numerical tests of these methods are omitted. 
\end{REM}

\section{Further numerical tests}\label{sec:num}
\setcounter{equation}{0}

In this section, we perform further numerical tests to examine the convergence of these first-order schemes. Our main interest is in the (local) LxF scheme and the cLxF scheme. 

To distinct subcritical, supercritical and transonic tests, we introduce the Froude number $Fr = u/\sqrt{gh}$ and define $W = (h, Fr)^T$. The subscripts $L$ and $R$ refer to the left and right states, respectively. The exact solutions for Examples \ref{ex-1s2r}--\ref{ex-1r2r-supern} are obtained through the exact solver in \cite{bernetti2008exact} and those for Examples \ref{ex-resonant} and \ref{ex-coincide} are obtained through the exact solver in \cite{han2014exact}. Values for exact Riemann solutions in the text are rounded to the fourth digit after the decimal point. 

In all the numerical tests, we set the CFL number to be $0.5$ and $g = 9.81$. %Unless otherwise noted, the meshes with $N =100$, $200$, $\cdots$, $25600$ grid points are used in the spatial discretization to produce the error tables. 
All tests concern Riemann problems with the initial data set as 
\begin{equation}
	W(x,0) = \left\{
	\begin{array}{cc}
		W_L\quad x<0,\\
		W_R\quad x\geq 0,
	\end{array}\right.
\end{equation}
where $W_L$ and $W_R$ are defined individually in each test. $b_L$ is always set as $0$, and $b_R$ is given for each test. $U(x,0)$ can be computed accordingly from $W(x,0)$. The settings and numerical observations are summarized in Table \ref{tab-settings}.

\begin{table}[h!]
	\centering
	\footnotesize
	\begin{tabular}{c|c|c|c|c|c|c|c|c|c|c|c|c}
		\hline
		\multirow{2}{*}{Test}&\multirow{2}{*}{Domain}&\multirow{2}{*}{$T$}
		&\multirow{2}{*}{$h_\mathrm{L}$}&\multirow{2}{*}{$h_\mathrm{R}$}
		&\multirow{2}{*}{$Fr_\mathrm{L}$}&\multirow{2}{*}{$Fr_\mathrm{R}$}
		&\multirow{2}{*}{$b_R$}&\multirow{2}{*}{Wave}&\multirow{2}{*}{$|Fr^\pm(0)|$}&\multirow{2}{*}{$\gamma$}&\multicolumn{2}{c}{$L^1$ convergence}\\
		\cline{12-13}
		&&&&&&&&&&&cLxF&LxF\\
		\hline
		\eqref{eq:db}&$[-5,5]$&$1$&$1$&$0.1$&$0$&$0$&$0.7$&1R-0-2S&Sub&$\star$&\cmark&\xmark\\
		\ref{ex-1s2r}&$[-5,5]$&$0.7$&$0.95$&$0.7$&$0.55$&$0.85$&$0.5$&1S-0-2R&Sub&$\star$&\cmark&\xmark\\
		\ref{ex-1r2r}&$[-5,5]$&$0.5$&$1$&$1.2$&$0.3$&$0.95$&$0.2$&1R-0-2R&Sub&$\star$&\cmark&\xmark\\
		\ref{ex-1s2s}&$[-5,5]$&$1$&$0.7$&$0.2$&$0.2$&$0.2$&$0.5$&1S-0-2S&Sub&$\star$&\cmark&\xmark\\
		\ref{ex-1s2s-super}&$[-1,5]$&$1$&$0.5$&$0.3$&$1.5$&$0.0$&$0.2$&1S-0-2S& Sub&$\star$&\cmark&\xmark\\
		\ref{ex-1s2s-supern}&$[-8,2]$&$1$&$0.5$&$0.7$&$-1.5$&$-1.05$&$0.2$&1S-2S-0& nSup&$\star$&\cmark&\halfcmark\\
		\ref{ex-1r2r-supern}&$[-8,2]$&$1$&$0.5$&$0.7$&$-2$&$-1.05$&$0.2$&1R-2R-0& nSup&$\star$&\cmark&\halfcmark\\
		\ref{ex-db-gamma0}&$[-5,5]$&$1$&$1$&$0.1$&$0$&$0$&$0.7$&1R-0-2S&Sub&$0$&\xmark&\xmark\\
		\ref{ex-coincide}&$[-5,15]$&$0.8$&$4$&$1.0299$&$1.1175$&$2.2428$&$1$&1S-$0$(R)-2S&Tran&$\star$&\xmark\xmark&\xmark\\
		%\ref{ex-coincide}&\multirow{2}{*}{$[-0.5,0.5]$}&\multirow{2}{*}{$0.03$}&\multirow{2}{*}{$0.5674$}&\multirow{2}{*}{$0.558$}&\multirow{2}{*}{$0.8283$}&\multirow{2}{*}{$-1.2822$}&\multirow{2}{*}{$0.8$}&\multirow{2}{*}{1S-$0'$-$0$-2S}&\multirow{2}{*}{Tran}&\halfcmark&\halfcmark\\
		\ref{ex-resonant}&$[-25,15]$&$0.8$&$6$&$8$&$-2.0855$&$0$&$1$&1R-2R-$0$(R)&Tran&$\star$&\xmark\xmark& \halfcmark \\
		\hline
	\end{tabular}\bigskip

\begin{tabular}{c|l}
	\hline
	Sub&subcritical test: $|Fr^-(0)|<1$ and $|Fr^+(0)|<1$\\
	nSup&negative supercritical test: $Fr^-(0)<-1$ and $Fr^+(0)<-1$\\
	Tran&transonic test: $|Fr^-(0)|<1<|Fr^+(0)|$ or $|Fr^+(0)|<1<|Fr^-(0)|$\\
	\hline
	$\star$&$\sgnb$\\
	\hline
	1S/1R& 1-shock wave / 1-rarefaction wave\\
	2S/2R& 2-shock wave / 2-rarefaction wave\\
	\hline
	-& unstable\\
	\cmark& converge to the exact solution without the numerical artifact\\
	\halfcmark& converge to the exact solution with the numerical artifact\\
	\xmark& converge to a wrong solution close to the exact solution with a similar wave pattern\\
	\xmark\xmark&converge to a wrong solution very different from the exact solution, \\
	&potentially with a wrong wave pattern\\
	\hline
\end{tabular}
\caption{Summary of numerical tests. We refer to the wave at $x = 0$ as the (stationary) 0-wave, the left-most wave as the 1-wave, and the other wave to be the 2-wave. The convergence is in the sense of $L^1$, and the ``correct convergence" is determined by whether the numerical solution exhibits a non-vanishing convergence rate in the error tables.}\label{tab-settings}
\end{table}

Based on Table \ref{tab-settings}, we summarize our numerical observations for problems tested in this paper. 
\begin{enumerate}
	\item For nontransonic problems (both the subcritical and the negative supercritical tests), the cLxF scheme with $\gamma = \sgnb$ will converge to the correct solution in $L^1$ without the numerical artifact.
 \item For transonic problems, the cLxF scheme with $\gamma = \sgnb$ will converge to a wrong solution with the wrong wave pattern. This behavior seems to be similar to the convergence to the non-entropy solution of certain schemes for hyperbolic conservation laws. 
 \item For both the transonic and nontransonic tests, the LxF scheme will suffer the aforementioned numerical artifact. Moreover, if the 1- and 2-waves are developed on both sides of the bottom step (which includes all subcritical tests and certain transonic tests), then the LxF scheme with $\gamma = \sgnb$ will converge to a wrong solution close to the exact solution; if the 1- and 2-waves are developed on the same side of the bottom step (which includes all negative supercritical tests and certain transonic tests), then the LxF scheme with $\gamma = \sgnb$ will converge to the correct exact solution (although with the numerical artifact at the bottom discontinuity).
%	\item For nontransonic subcritical tests, the LxF scheme captures the correct wave patterns, but converges to a wrong solution close to the exact solution and suffers the numerical artifact at the bottom discontinuity. 
%	\item For nontransonic negative supercritical tests, the LxF scheme converges to the correct solution in $L^1$. But it suffers the numerical artifact at the bottom discontinuity. 
	\item $\gamma = \sgnb$ is needed for the correct convergence. The $\gamma$ used to define the exact solution should be consistent with the $\gamma$ used in the numerical scheme. 
%	\item For the transonic tests Example \ref{ex-coincide}, the cLxF scheme may converge to a wrong solution with the wrong wave pattern. This seems to be similar to the convergence to a non-entropy solution for homogeneous hyperbolic conservation laws in the resonant case. The LxF scheme will suffer the numerical fact, while it may either converge to the correct solution or converge to a wrong solution close to the exact solution.
%It seems that the LxF solution is driven towards an entropy solution, but the numerical artifact caused by the discontinuous bottom rune the $L^1$ convergence. 
%	\item The cLxF scheme may still suffer the numerical artifact for certain transonic problems.
\end{enumerate}
In the rest of this section, we present detailed numerical results of these tests.

\begin{ex}\label{ex-1s2r}
	The exact solution of this Riemann test consists of a 1-shock and a 2-rarefaction. The exact solution admits the following states in $W$:
	\begin{equation}
		W_L = \begin{pmatrix}
			0.95\\
			0.55
		\end{pmatrix}\xrightarrow{\text{1-shock}}\begin{pmatrix}
		1.2295\\
		0.24
	\end{pmatrix}\xrightarrow{\text{0-wave}}
	\begin{pmatrix}
		0.5814\\
		0.7381
	\end{pmatrix}\xrightarrow{\text{2-rarefaction}}
	\begin{pmatrix}
		0.7\\0.85
	\end{pmatrix} = W_R.
	\end{equation}
	
From Figure \ref{fig-1s2r}, we clearly observe that the numerical solution of the LxF scheme forms a downward spike in the momentum at $x = 0$. In the mean time, there is a clear mismatch between the exact solution and the numerical solution even with an extremely refined mesh. The error table (Table \ref{tab-1s2r}) also confirms that the numerical solution does not converge to the exact solution. But from Figure \ref{fig-1s2r} and Table \ref{tab-1s2r}, we can see that the cLxF scheme converges to the exact solution in $L^1$ without the spurious spike in $m$ and without any transition points in $h$. 
	\begin{figure}[!ht]
		\centering
		\subfloat[cLxF, $h$.]{%
			\includegraphics[width=0.25\textwidth,trim={2cm 0 1cm 0}]{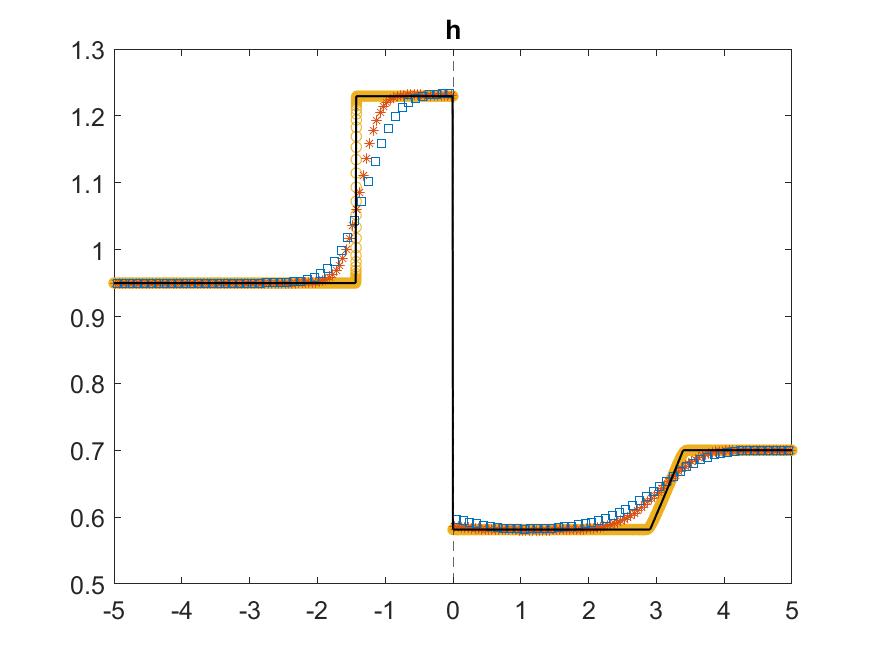}
		}
		\subfloat[cLxF, $m = hu$.]{%
			\includegraphics[width=0.25\textwidth,trim={2cm 0 1cm 0}]{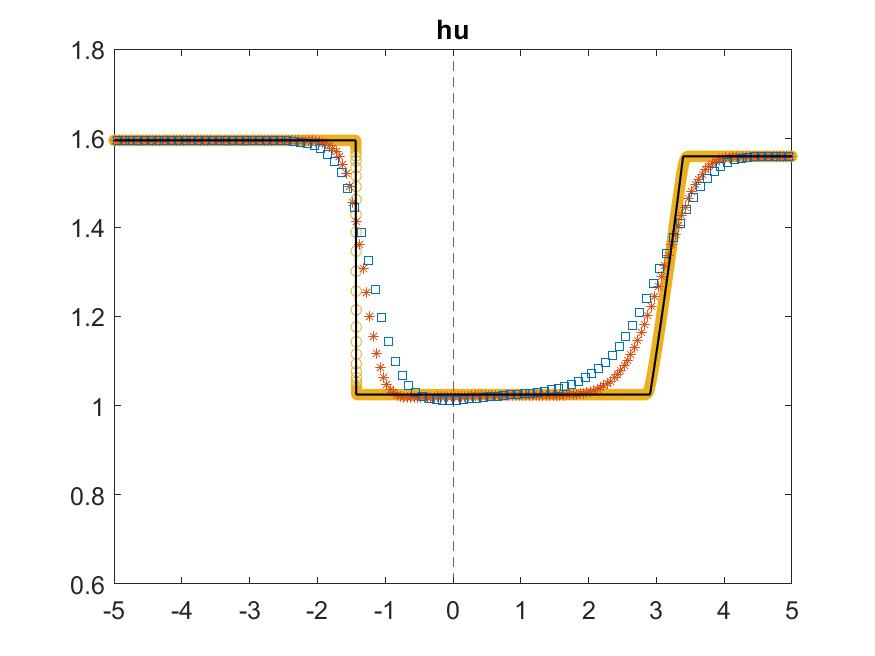}
		}
		\subfloat[LxF, $h$.]{%
			\includegraphics[width=0.25\textwidth,trim={2cm 0 1cm 0}]{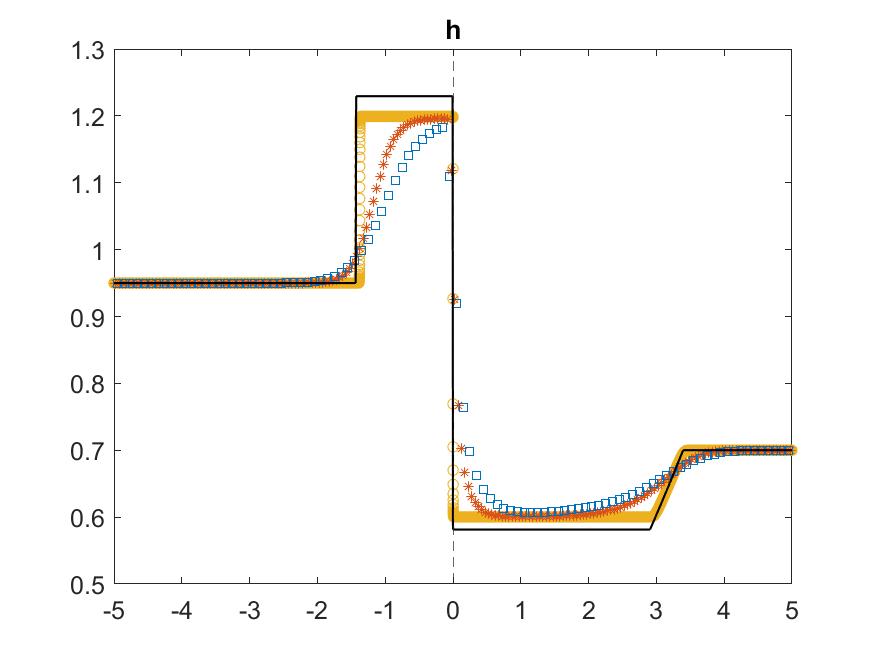}
		}
		\subfloat[LxF, $m = hu$.]{%
			\includegraphics[width=0.25\textwidth,trim={2cm 0 1cm 0}]{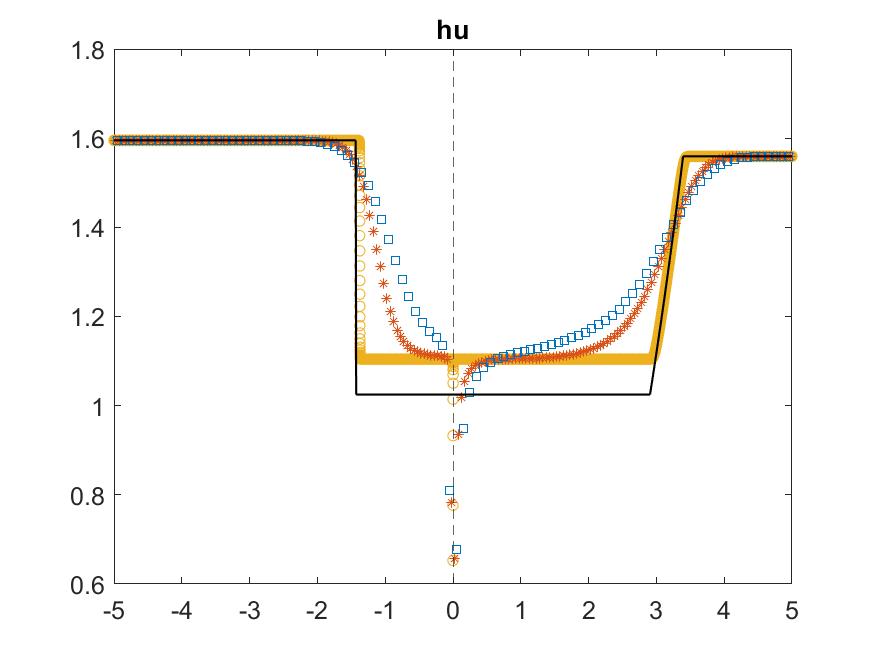}
		}
		\caption{Solutions to Example \ref{ex-1s2r}. Black solid line: exact solution; blue squares: $N = 100$; red stars: $N = 200$, yellow circles: $N = 25600$.}\label{fig-1s2r}
	\end{figure}
	\begin{table}[h!]
		\centering
		\begin{tabular}{c|c|c|c|c|c|c|c|c}
			\hline
			$N$&cLxF $e_{h}$&order&cLxF $e_{\mm}$&order&LxF $e_{h}$&order&LxF $e_{\mm}$&order\\
			\hline
			100&  1.42e-01&  -&  4.05e-01&  -&  3.65e-01&  -&  8.88e-01&  -\\ 
			200&  8.53e-02&  0.74&  2.52e-01&  0.69&  2.60e-01&  0.49&  6.87e-01&  0.37\\ 
			400&  4.83e-02&  0.82&  1.45e-01&  0.80&  1.87e-01&  0.47&  5.39e-01&  0.35\\ 
			800&  2.73e-02&  0.82&  8.35e-02&  0.80&  1.50e-01&  0.32&  4.63e-01&  0.22\\ 
			1600&  1.53e-02&  0.83&  4.84e-02&  0.79&  1.32e-01&  0.18&  4.24e-01&  0.12\\ 
			3200&  8.34e-03&  0.88&  2.72e-02&  0.83&  1.22e-01&  0.12&  4.00e-01&  0.08\\ 
			6400&  4.57e-03&  0.87&  1.54e-02&  0.82&  1.17e-01&  0.06&  3.88e-01&  0.04\\ 
			12800&  2.49e-03&  0.87&  8.62e-03&  0.84&  1.14e-01&  0.03&  3.81e-01&  0.03\\ 
			25600&  1.37e-03&  0.87&  4.84e-03&  0.83&  1.13e-01&  0.02&  3.78e-01&  0.01\\ 
			\hline
		\end{tabular}
	\caption{$L^1$ error table for Example \ref{ex-1s2r}.}\label{tab-1s2r}
	\end{table}
\end{ex}
\begin{ex}\label{ex-1r2r}
		The exact solution of this Riemann test consists of a 1-rarefaction and a 2-rarefaction. The exact solution admits the following states in $W$:
		\begin{equation}
			W_L = \begin{pmatrix}
				1\\
				0.3
			\end{pmatrix}\xrightarrow{\text{1-rarefaction}}\begin{pmatrix}
				0.9443\\
				0.3669
			\end{pmatrix}\xrightarrow{\text{0-wave}}
			\begin{pmatrix}
				0.6780\\
				0.6031
			\end{pmatrix}\xrightarrow{\text{2-rarefaction}}
			\begin{pmatrix}
				1.2\\0.95
			\end{pmatrix} = W_R.
		\end{equation}
		
		Again, from Figure \ref{fig-1r2r}, the LxF solution forms a spurious spike in the numerical momentum $m$ and admits transition points in the numerical water height $h$, while the cLxF solution does not suffer similar numerical artifacts. From Table \ref{tab-1r2r}, we see that the order of accuracy for LxF scheme keeps decreasing as we refine the mesh. Indeed, we expect that the convergence rate will decrease to 0 and the numerical error will be trapped at some nonzero value as we keep refining the mesh. In contrast, the cLxF scheme converges to the exact solution in $L^1$ at a rate between $0.5$ and $1$. 
		\begin{figure}[!h]
			\centering
		\subfloat[cLxF, $h$.]{%
			\includegraphics[width=0.25\textwidth,trim={2cm 0 1cm 0}]{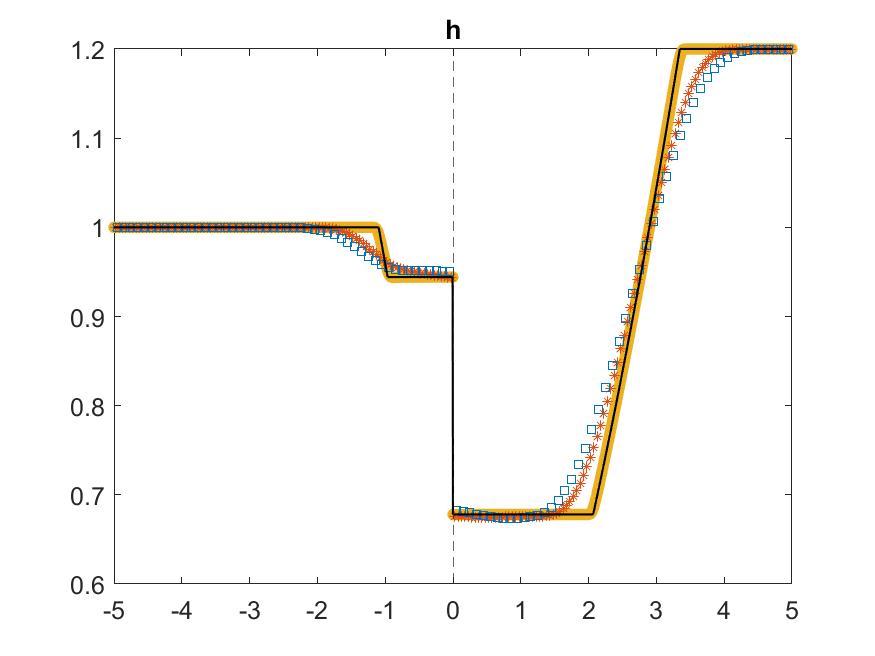}
		}
		\subfloat[cLxF, $m = hu$.]{%
			\includegraphics[width=0.25\textwidth,trim={2cm 0 1cm 0}]{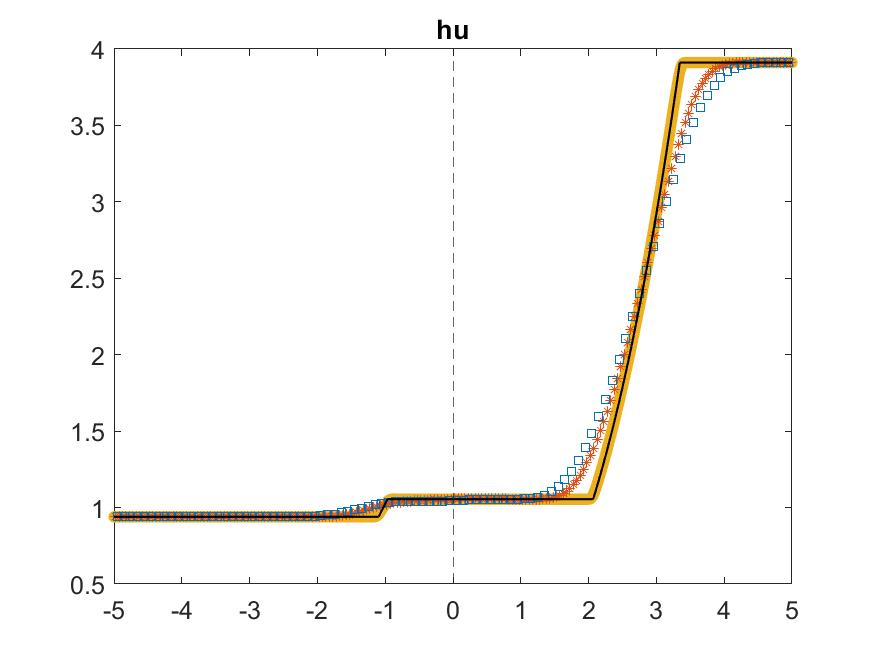}
		}
		\subfloat[LxF, $h$.]{%
			\includegraphics[width=0.25\textwidth,trim={2cm 0 1cm 0}]{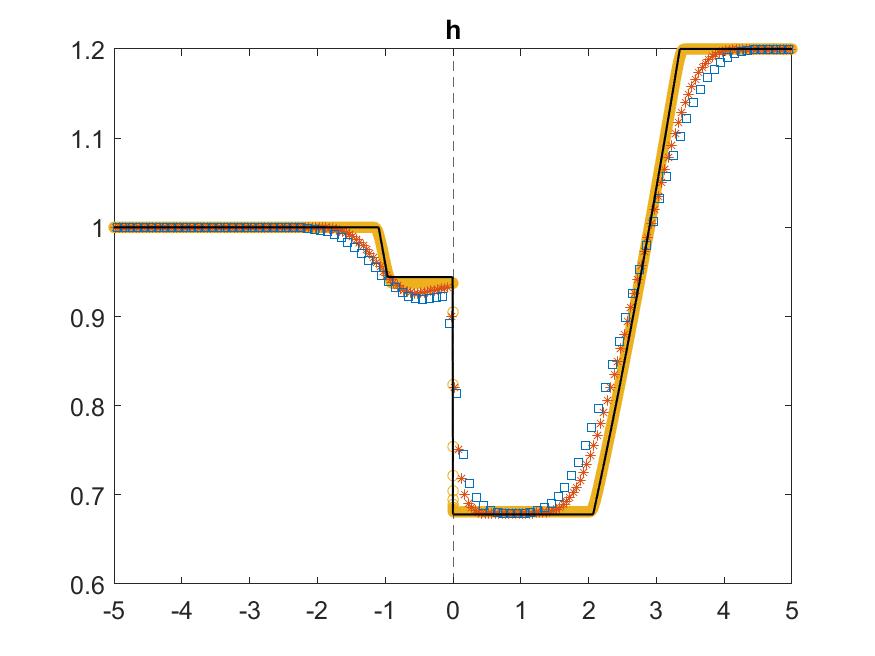}
		}
		\subfloat[LxF, $m = hu$.]{%
			\includegraphics[width=0.25\textwidth,trim={2cm 0 1cm 0}]{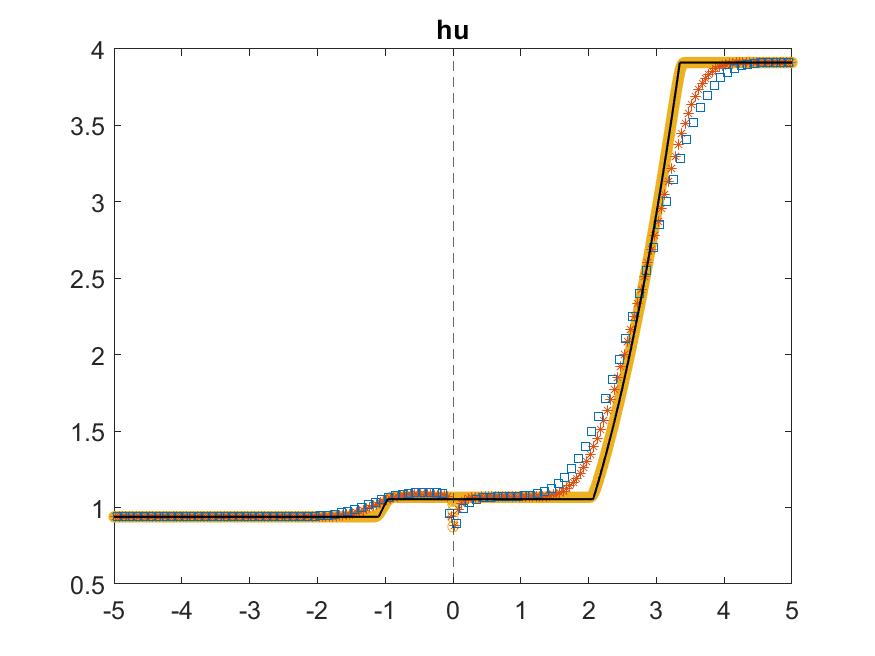}
		}
		\caption{Solutions to Example \ref{ex-1r2r}. Black solid line: exact solution; blue squares: $N = 100$; red stars: $N = 200$, yellow circles: $N = 25600$.}\label{fig-1r2r}
	\end{figure}
\begin{table}[h!]
	\centering
	\begin{tabular}{c|c|c|c|c|c|c|c|c}
		\hline
		$N$&cLxF $e_{h}$&order&cLxF $e_{\mm}$&order&LxF $e_{h}$&order&LxF $e_{\mm}$&order\\
		\hline
		100&  1.53e-01&  -&  7.45e-01&  -&  2.00e-01&  -&  8.31e-01&  -\\ 
		200&  9.89e-02&  0.63&  4.72e-01&  0.66&  1.27e-01&  0.66&  5.40e-01&  0.62\\ 
		400&  6.26e-02&  0.66&  2.96e-01&  0.68&  8.06e-02&  0.65&  3.47e-01&  0.64\\ 
		800&  3.89e-02&  0.69&  1.82e-01&  0.70&  5.28e-02&  0.61&  2.25e-01&  0.63\\ 
		1600&  2.38e-02&  0.71&  1.09e-01&  0.73&  3.62e-02&  0.54&  1.49e-01&  0.59\\ 
		3200&  1.43e-02&  0.74&  6.46e-02&  0.76&  2.64e-02&  0.46&  1.03e-01&  0.53\\ 
		6400&  8.45e-03&  0.76&  3.76e-02&  0.78&  2.06e-02&  0.36&  7.57e-02&  0.45\\ 
		12800&  4.94e-03&  0.77&  2.16e-02&  0.80&  1.72e-02&  0.26&  5.97e-02&  0.34\\ 
		25600&  2.85e-03&  0.79&  1.23e-02&  0.82&  1.53e-02&  0.17&  5.04e-02&  0.24\\ 
		\hline
	\end{tabular}
	\caption{$L^1$ error table for Example \ref{ex-1r2r}.}\label{tab-1r2r}
\end{table}

\end{ex}

\begin{ex}\label{ex-1s2s}
		The exact solution of this Riemann test consists of a 1-shock and a 2-shock. The exact solution admits the following states in $W$:
		\begin{equation}
			W_L = \begin{pmatrix}
				0.7\\
				0.2
			\end{pmatrix}\xrightarrow{\text{1-shock}}\begin{pmatrix}
				0.7849\\
				0.0774
			\end{pmatrix}\xrightarrow{\text{0-wave}}
			\begin{pmatrix}
				0.2569\\
				0.4133
			\end{pmatrix}\xrightarrow{\text{2-shock}}
			\begin{pmatrix}
				0.2\\0.2
			\end{pmatrix} = W_R.
		\end{equation}
		Again, from Figure \ref{fig-1s2s} and Table \ref{tab-1s2s}, we see that the LxF scheme converges to a wrong solution, with a spurious spike in $m$ and transition points in $h$ at $x = 0$. The convergence rate decays to below $0.1$ at $N=800$. On the other hand, the cLxF scheme converges to the exact solution without such numerical artifacts, and the convergence rate approaches $1$.
		\begin{figure}[h!]
			\centering
		\subfloat[cLxF, $h$.]{%
			\includegraphics[width=0.25\textwidth,trim={2cm 0 1cm 0}]{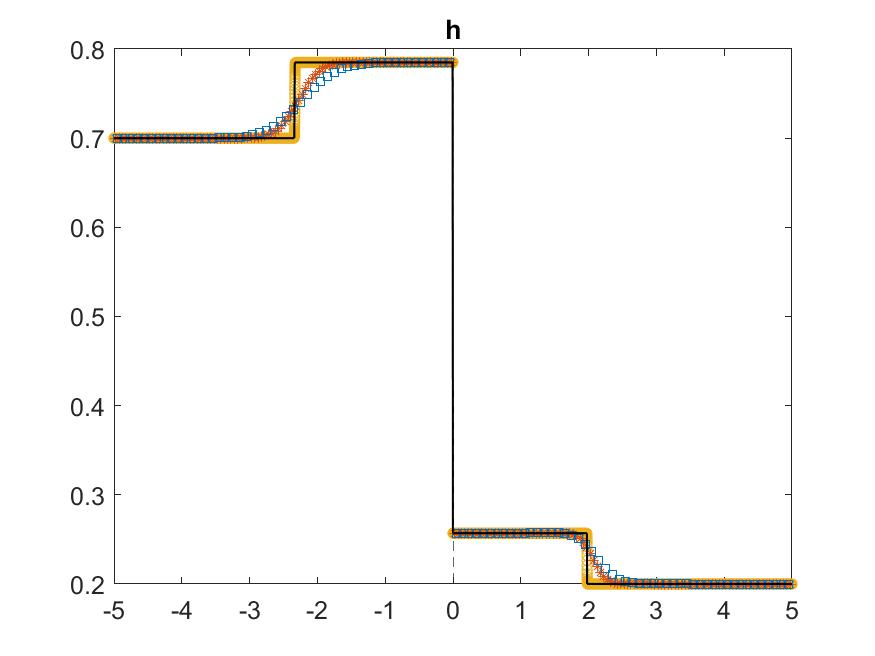}
		}
%		\hfill
		\subfloat[cLxF, $m = hu$.]{%
			\includegraphics[width=0.25\textwidth,trim={2cm 0 1cm 0}]{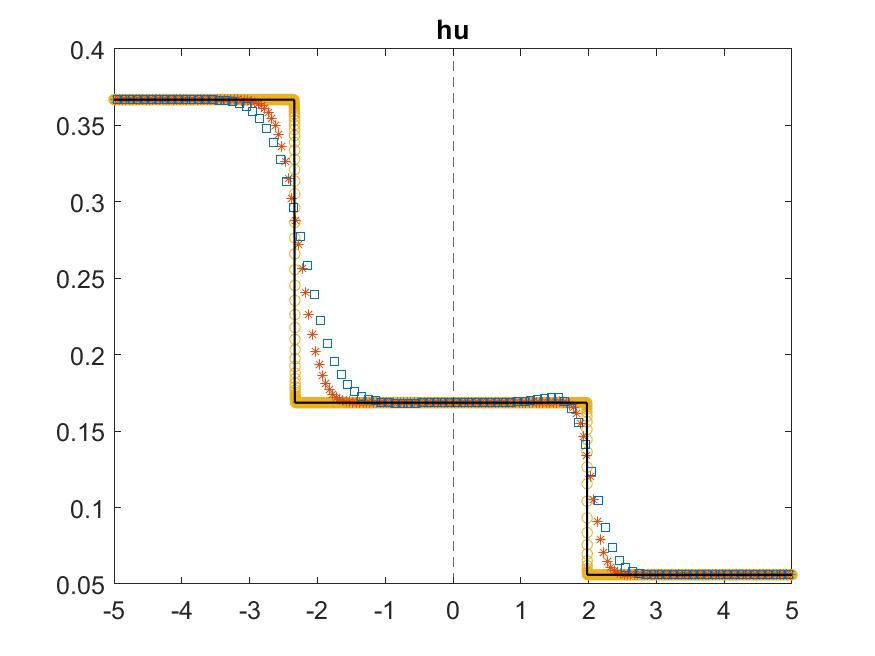}
		}
		\subfloat[LxF, $h$.]{%
			\includegraphics[width=0.25\textwidth,trim={2cm 0 1cm 0}]{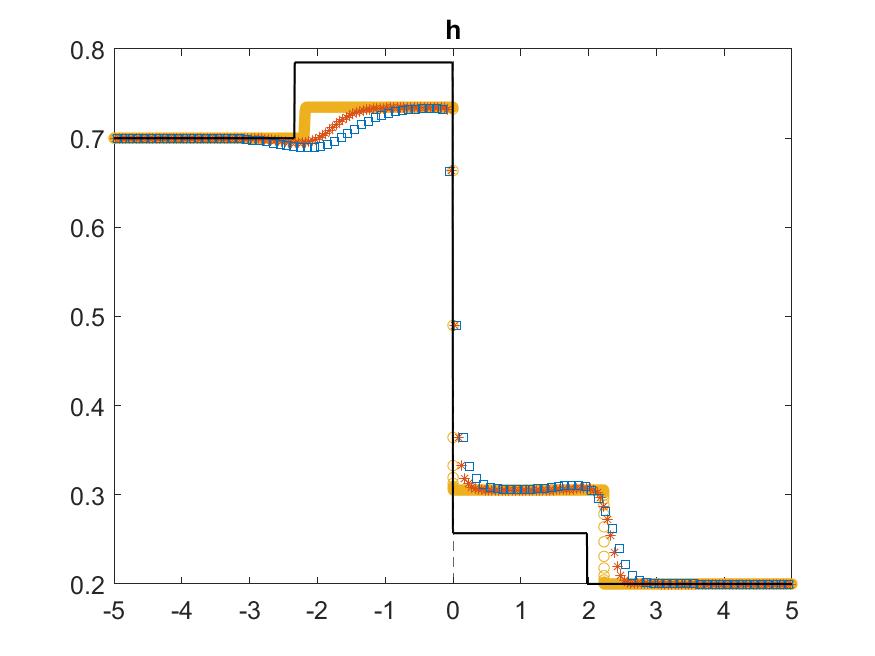}
		}
%		\hfill
		\subfloat[LxF, $m = hu$.]{%
			\includegraphics[width=0.25\textwidth,trim={2cm 0 1cm 0}]{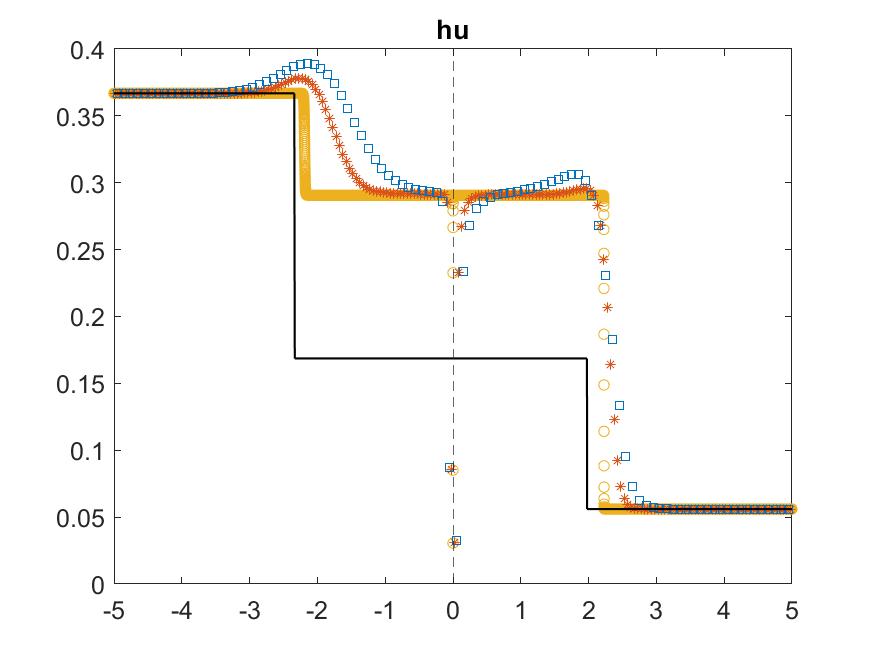}
		}
		\caption{Solutions to Example \ref{ex-1s2s}. Black solid line: exact solution; blue squares: $N = 100$; red stars: $N = 200$, yellow circles: $N = 25600$.}\label{fig-1s2s}
	\end{figure}
	\begin{table}[h!]
		\centering
		\begin{tabular}{c|c|c|c|c|c|c|c|c}
			\hline
			$N$&cLxF $e_{h}$&order&cLxF $e_{\mm}$&order&LxF $e_{h}$&order&LxF $e_{\mm}$&order\\
			\hline
			100&  4.18e-02&  -&  9.20e-02&  -&  3.42e-01&  -&  7.15e-01&  -\\ 
			200&  2.63e-02&  0.67&  5.78e-02&  0.67&  2.96e-01&  0.21&  6.61e-01&  0.11\\ 
			400&  1.60e-02&  0.72&  3.48e-02&  0.73&  2.70e-01&  0.13&  6.27e-01&  0.08\\ 
			800&  9.23e-03&  0.79&  2.03e-02&  0.77&  2.58e-01&  0.07&  6.13e-01&  0.03\\ 
			1600&  4.91e-03&  0.91&  1.08e-02&  0.91&  2.51e-01&  0.04&  6.03e-01&  0.02\\ 
			3200&  2.47e-03&  0.99&  5.52e-03&  0.97&  2.47e-01&  0.02&  5.99e-01&  0.01\\ 
			6400&  1.25e-03&  0.98&  2.80e-03&  0.98&  2.46e-01&  0.01&  5.97e-01&  0.00\\ 
			12800&  6.07e-04&  1.05&  1.36e-03&  1.05&  2.45e-01&  0.01&  5.96e-01&  0.00\\ 
			25600&  3.11e-04&  0.96&  6.95e-04&  0.96&  2.44e-01&  0.00&  5.95e-01&  0.00\\ 
			\hline
		\end{tabular}
	\caption{$L^1$ error table for Example \ref{ex-1s2s}.}\label{tab-1s2s}
	\end{table}
	
\end{ex}

\begin{ex}\label{ex-1s2s-super}
	This Riemann test consists of a 1-shock and 2-shock on two sides of the bottom step, admitting the following states in $W$:
	\begin{equation}
		W_L = \begin{pmatrix}
			0.5\\
			1.5
		\end{pmatrix}\xrightarrow{\text{1-shock}}\begin{pmatrix}
			1.0141\\
			0.4295
		\end{pmatrix}\xrightarrow{\text{0-wave}}
		\begin{pmatrix}
			0.7041\\
			0.7424
		\end{pmatrix}\xrightarrow{\text{2-shock}}
		\begin{pmatrix}
			0.3\\0
		\end{pmatrix} = W_R.
	\end{equation}
	
	This differs from Example \ref{ex-1s2s} with a supercritical velocity for the initial left states. Besides the numerical artifact and the wrong convergence of the LxF scheme at $x = 0$ (see Figure \ref{fig-1s2s-super}), we want to address the solution profile at the 1-shock. See Figure \ref{fig-1s2s-super-zoomin} for the zoomed-in pictures. We note that the LxF scheme admits an overshoot at the 1-shock and captures the wrong shock speed. For the cLxF scheme, the correct speed for 1-shock is captured. But its numerical solution still suffers the overshoot. This is reasonable since we only remove the numerical viscosity at $x = 0$ to eliminate the spurious spike there, while the overshoot or undershoot elsewhere has not been taken care of. 
	\begin{figure}[h!]
		\centering
	\subfloat[cLxF, $h$.]{%
		\includegraphics[width=0.25\textwidth,trim={2cm 0 1cm 0}]{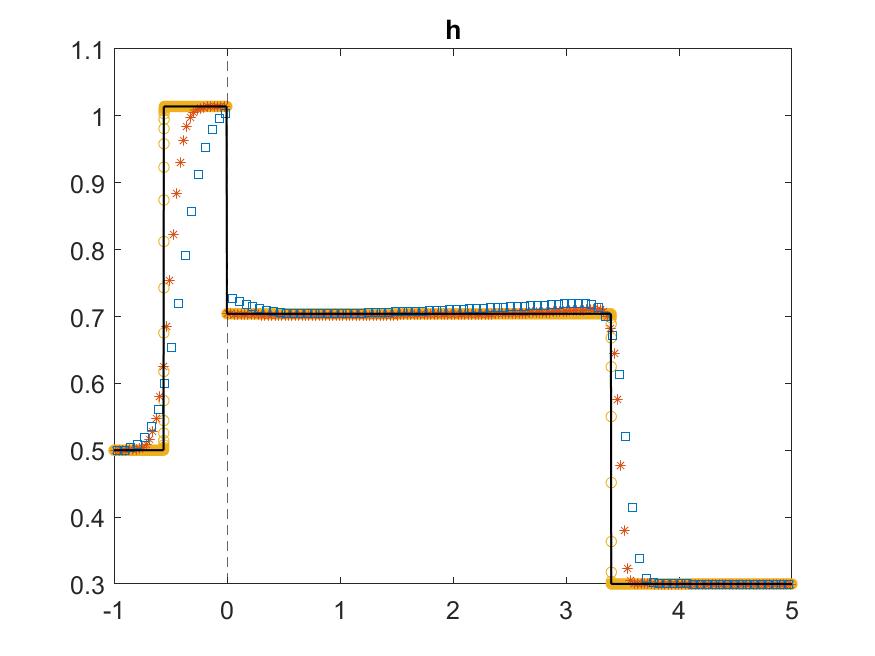}
	}
	\subfloat[cLxF, $m = hu$.]{%
		\includegraphics[width=0.25\textwidth,trim={2cm 0 1cm 0}]{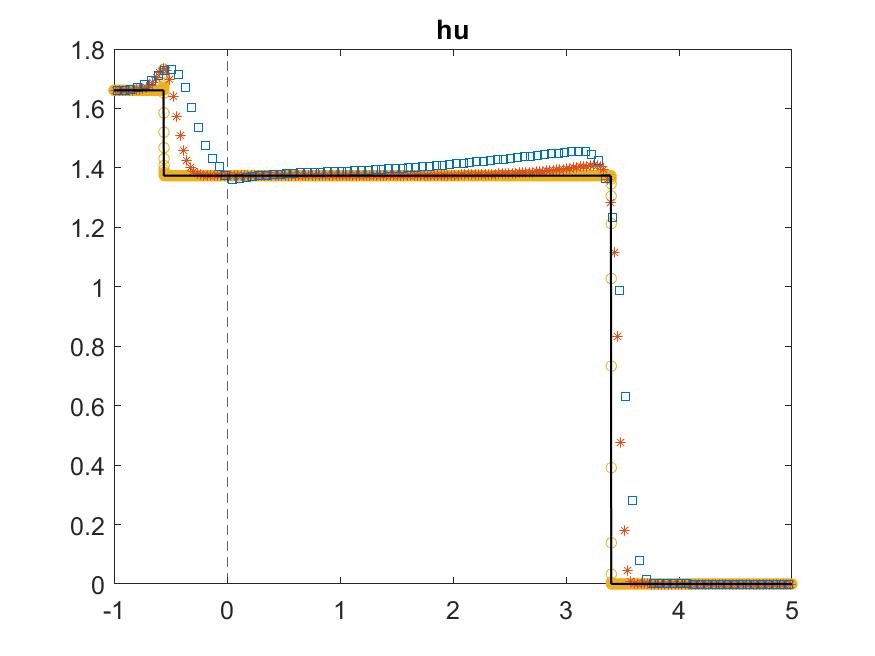}
	}
	\subfloat[LxF, $h$.]{%
		\includegraphics[width=0.25\textwidth,trim={2cm 0 1cm 0}]{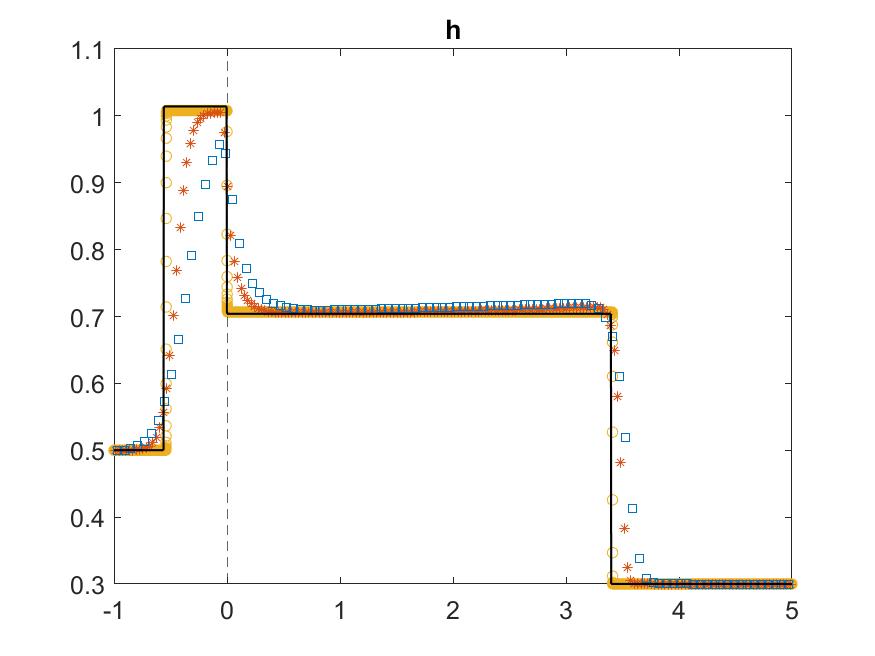}
	}
	\subfloat[LxF, $m = hu$.]{%
		\includegraphics[width=0.25\textwidth,trim={2cm 0 1cm 0}]{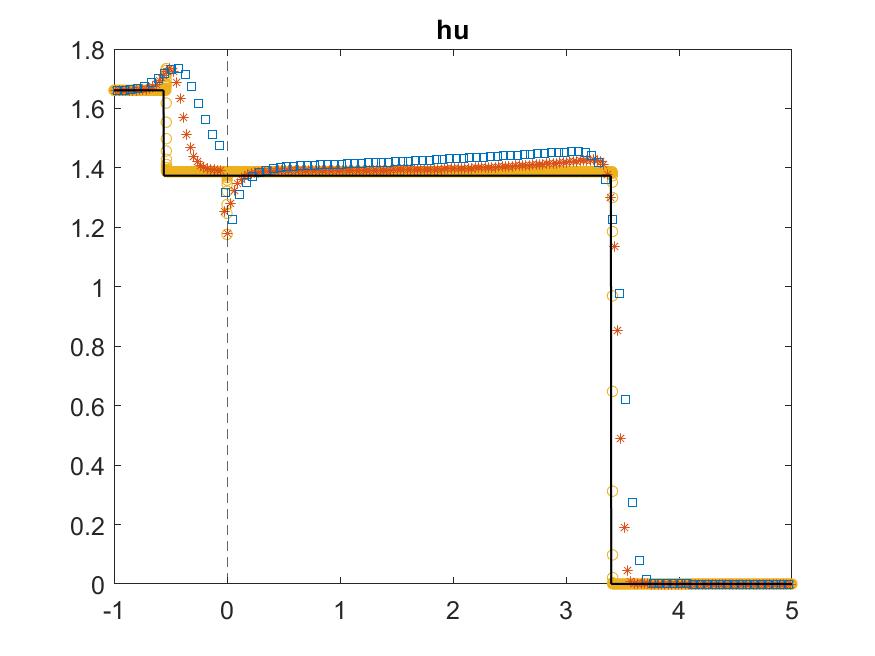}
	}
	\caption{Solutions to Example \ref{ex-1s2s-super}. Black solid line: exact solution; blue squares: $N = 100$; red stars: $N = 200$, yellow circles: $N = 25600$.}\label{fig-1s2s-super}
\end{figure}
\begin{figure}[h!]
	\centering
	\subfloat[cLxF, $m = hu$.]{%
		\includegraphics[width=0.4\textwidth,trim={1cm 0 1cm 0}]{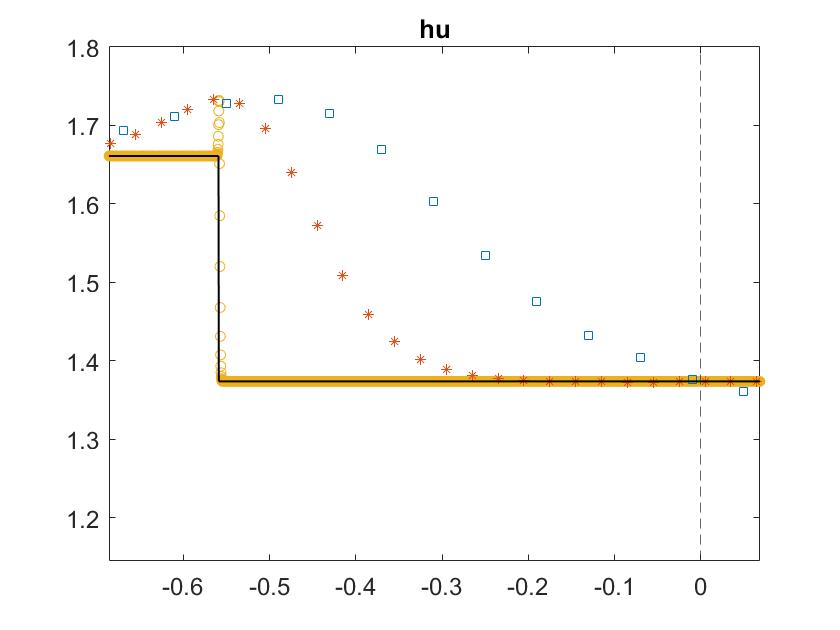}
	}
	\subfloat[LxF, $m = hu$.]{%
		\includegraphics[width=0.4\textwidth,trim={1cm 0 1cm 0}]{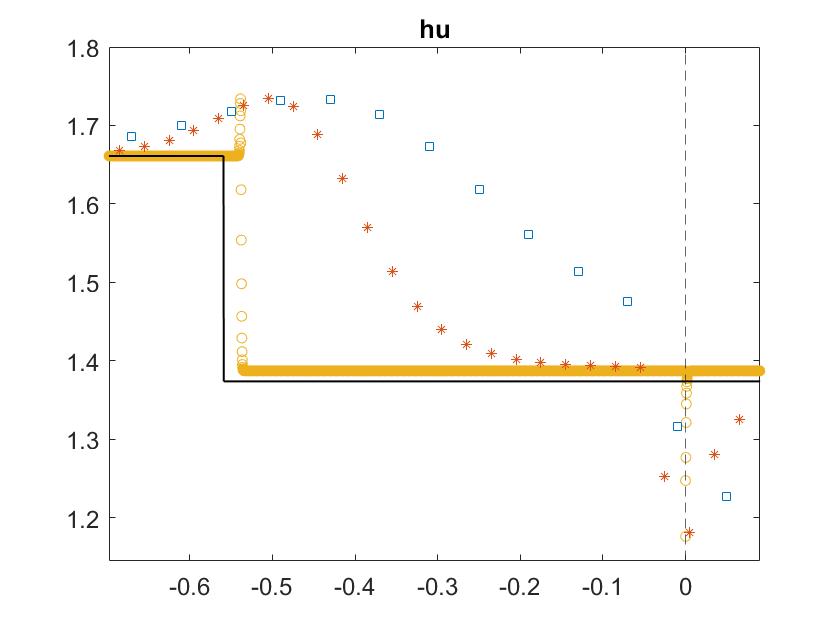}
	}
	\caption{Zoomed-in solutions to Example \ref{ex-1s2s-super}. Black solid line: exact solution; blue squares: $N = 100$; red stars: $N = 200$, yellow circles: $N = 25600$.}\label{fig-1s2s-super-zoomin}
\end{figure}
	\begin{table}[h!]
		\centering
		\begin{tabular}{c|c|c|c|c|c|c|c|c}
			\hline
			$N$&cLxF $e_{h}$&order&cLxF $e_{\mm}$&order&LxF $e_{h}$&order&LxF $e_{\mm}$&order\\
			\hline
			100&  1.95e-01&  -&  4.32e-01&  -&  2.57e-01&  -&  5.19e-01&  -\\ 
			200&  7.73e-02&  1.34&  1.55e-01&  1.48&  1.21e-01&  1.09&  2.55e-01&  1.02\\ 
			400&  5.47e-02&  0.50&  1.10e-01&  0.50&  9.13e-02&  0.41&  1.97e-01&  0.38\\ 
			800&  2.06e-02&  1.41&  4.04e-02&  1.44&  5.12e-02&  0.83&  1.19e-01&  0.72\\ 
			1600&  1.25e-02&  0.72&  2.68e-02&  0.59&  4.24e-02&  0.27&  1.04e-01&  0.20\\ 
			3200&  5.26e-03&  1.25&  1.02e-02&  1.40&  3.38e-02&  0.33&  8.53e-02&  0.29\\ 
			6400&  3.33e-03&  0.66&  6.64e-03&  0.61&  3.12e-02&  0.11&  8.13e-02&  0.07\\ 
			12800&  1.36e-03&  1.29&  2.78e-03&  1.26&  2.87e-02&  0.12&  7.70e-02&  0.08\\ 
			25600&  7.83e-04&  0.80&  1.64e-03&  0.76&  2.79e-02&  0.04&  7.57e-02&  0.03\\ 
			\hline
		\end{tabular}
	\caption{$L^1$ error table for Example \ref{ex-1s2s-super}.}
	\end{table}
	
\end{ex}

\begin{ex}\label{ex-1s2s-supern}
	This is a negative supercritical Riemann test taken from \cite[Section 4.5]{bernetti2008exact}. The 2-wave moves downstream to the left across the bottom step. Both 1-shock and 2-shock are on the left of the 0-wave. The exact solution admits the following states in $W$:
	\begin{equation}
		W_L = \begin{pmatrix}
			0.5\\
			-1.5
		\end{pmatrix}\xrightarrow{\text{1-shock}}\begin{pmatrix}
			0.5565\\
			-1.5262
		\end{pmatrix}\xrightarrow{\text{2-shock}}
		\begin{pmatrix}
			0.5138\\
			-1.6697
		\end{pmatrix}\xrightarrow{\text{0-wave}}
		\begin{pmatrix}
			0.7\\-1.05
		\end{pmatrix} = W_R.
	\end{equation}
	
	In this test, although the LxF method still suffers the previously stated numerical artifacts, its numerical solution indeed converges to the exact solution in $L^1$. See Figure \ref{fig-1s2s-supern} and Table \ref{tab-1s2s-supern}. The cLxF scheme can converge to the exact solution without a spike in $m$ or transition points in $h$ at the bottom step. 
	\begin{figure}[h!]
		\centering
		\subfloat[cLxF, $h$.]{%
			\includegraphics[width=0.25\textwidth,trim={2cm 0 1cm 0}]{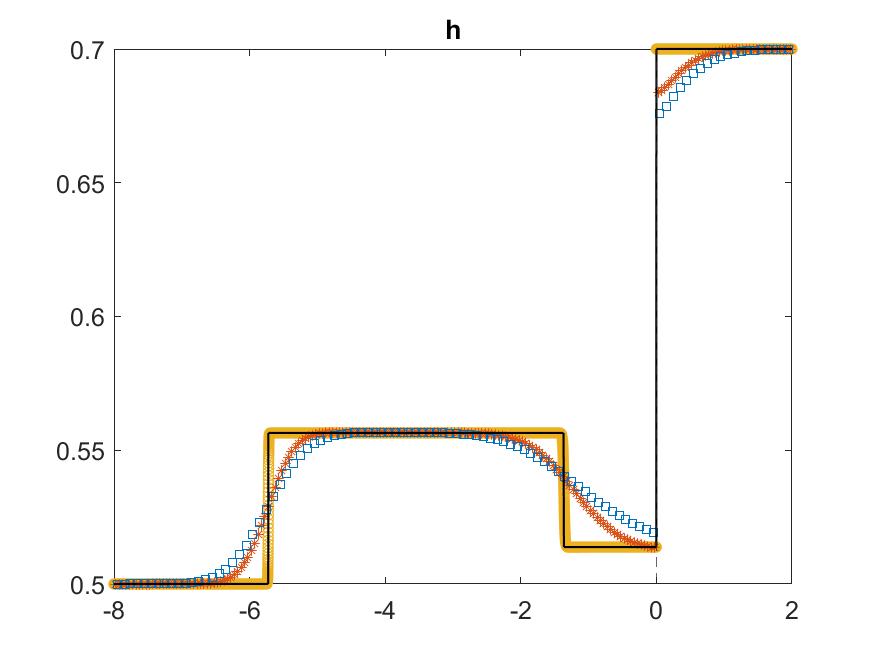}
		}
		\subfloat[cLxF, $m = hu$.]{%
			\includegraphics[width=0.25\textwidth,trim={2cm 0 1cm 0}]{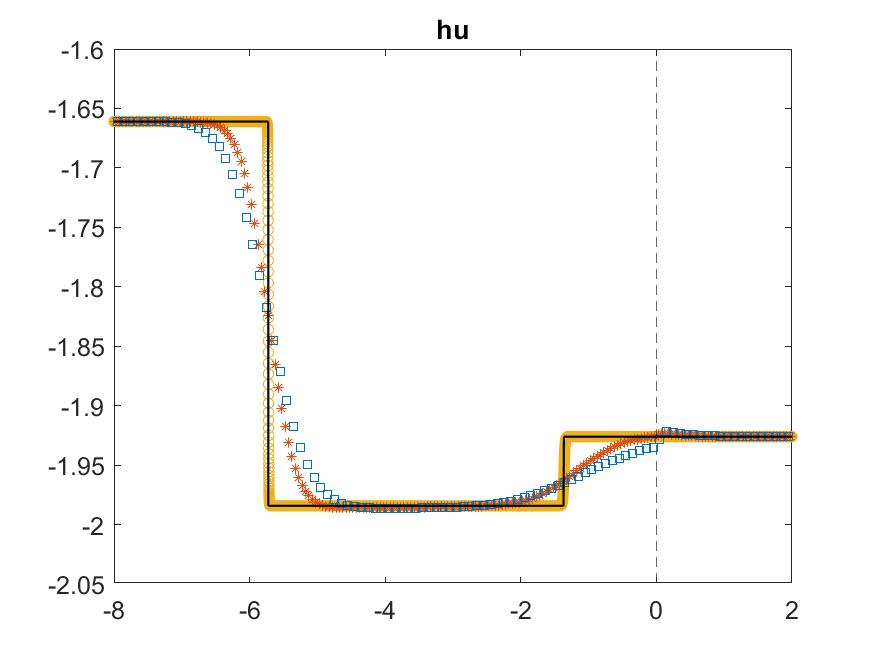}
		}
		\subfloat[LxF, $h$.]{%
			\includegraphics[width=0.25\textwidth,trim={2cm 0 1cm 0}]{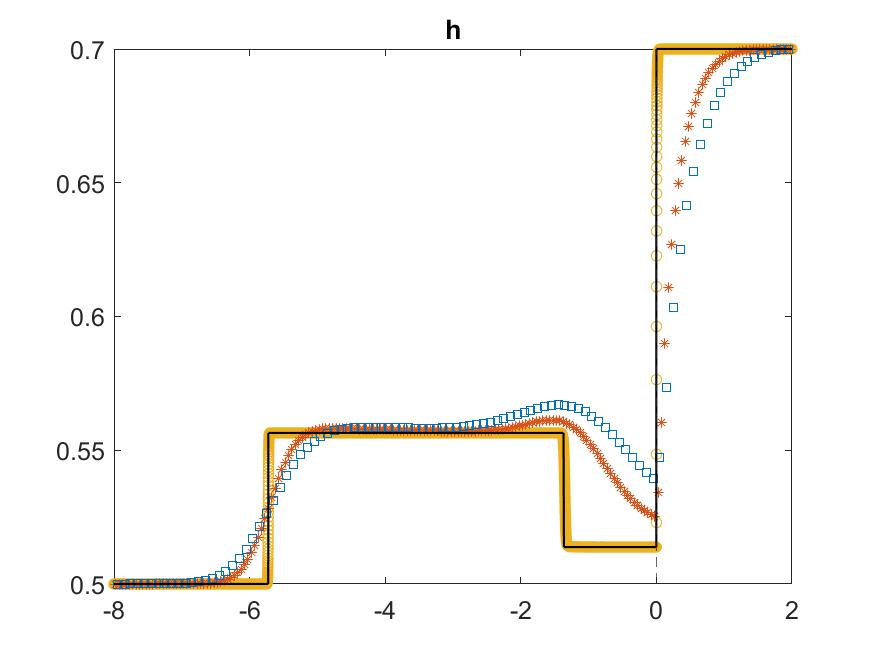}
		}
		\subfloat[LxF, $m = hu$.]{%
			\includegraphics[width=0.25\textwidth,trim={2cm 0 1cm 0}]{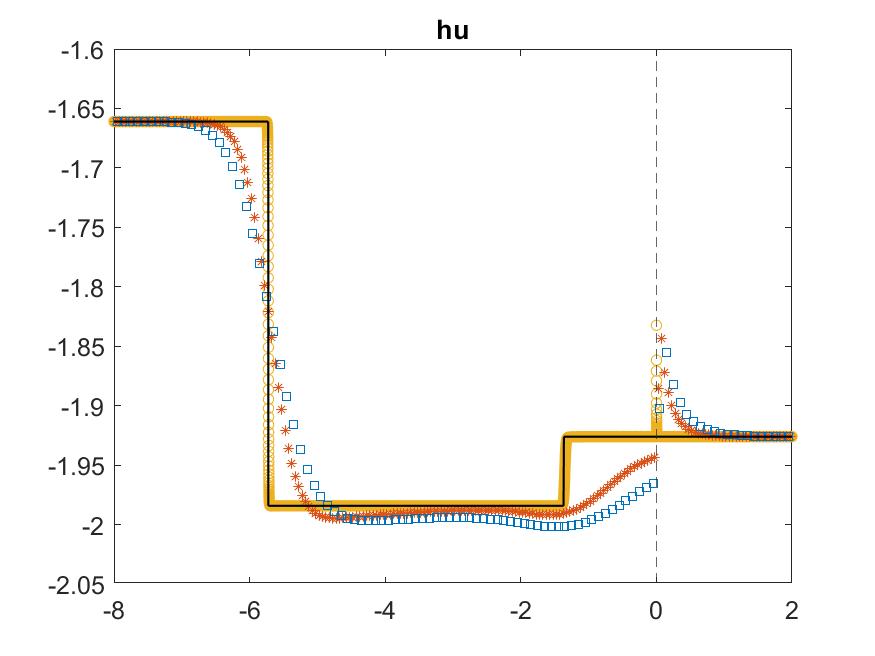}
		}
		\caption{Solutions to Example \ref{ex-1s2s-supern}. Black solid line: exact solution; blue squares: $N = 100$; red stars: $N = 200$, yellow circles: $N = 25600$.}\label{fig-1s2s-supern}
	\end{figure}
\begin{table}[h!]
	\centering
	\begin{tabular}{c|c|c|c|c|c|c|c|c}
		\hline
		$N$&cLxF $e_{h}$&order&cLxF $e_{\mm}$&order&LxF $e_{h}$&order&LxF $e_{\mm}$&order\\
		\hline
		100&  6.30e-02&  -&  1.66e-01&  -&  1.61e-01&  -&  2.61e-01&  -\\ 
		200&  3.98e-02&  0.66&  1.10e-01&  0.59&  1.04e-01&  0.63&  1.70e-01&  0.62\\ 
		400&  2.42e-02&  0.72&  7.08e-02&  0.63&  6.44e-02&  0.69&  1.06e-01&  0.68\\ 
		800&  1.46e-02&  0.73&  4.44e-02&  0.67&  3.78e-02&  0.77&  6.26e-02&  0.76\\ 
		1600&  8.76e-03&  0.74&  2.68e-02&  0.73&  2.17e-02&  0.80&  3.62e-02&  0.79\\ 
		3200&  5.09e-03&  0.78&  1.54e-02&  0.80&  1.21e-02&  0.85&  2.02e-02&  0.84\\ 
		6400&  2.82e-03&  0.85&  8.30e-03&  0.89&  6.35e-03&  0.92&  1.06e-02&  0.92\\ 
		12800&  1.49e-03&  0.92&  4.28e-03&  0.96&  3.23e-03&  0.97&  5.40e-03&  0.98\\ 
		25600&  7.52e-04&  0.98&  2.15e-03&  0.99&  1.61e-03&  1.01&  2.69e-03&  1.00\\ 
		\hline
	\end{tabular}
\caption{$L^1$ error table for Example \ref{ex-1s2s-supern}.}\label{tab-1s2s-supern}
\end{table}

\end{ex}
\begin{ex}\label{ex-1r2r-supern}
		This is another negative supercritical Riemann test taken from \cite[Section 4.5]{bernetti2008exact}. By increasing the flow speed of the left state in Example \ref{ex-1s2s-supern}, both the 1-wave and the 2-wave become rarefaction and are located on the left of the 0-wave. The exact solution admits the following states in $W$:
		\begin{equation}
			W_L = \begin{pmatrix}
				0.5\\
				-2
			\end{pmatrix}\xrightarrow{\text{1-shock}}\begin{pmatrix}
				0.4325\\
				-2
			\end{pmatrix}\xrightarrow{\text{2-shock}}
			\begin{pmatrix}
				0.5138\\
				-1.6697
			\end{pmatrix}\xrightarrow{\text{0-wave}}
			\begin{pmatrix}
				0.7\\-1.05
			\end{pmatrix} = W_R.
		\end{equation}
		
		The result is similar to the previous negative supercritical case: from Figure \ref{fig-1r2r-supern} and Table \ref{tab-1r2r-supern}, one can see that the LxF scheme converges to the exact solution in $L^1$ with numerical artifacts at the bottom step, while the cLxF converges to the exact solution in $L^1$ without a spurious spike in $m$ and transition points in $h$ as the mesh is sufficiently refined. 
		\begin{figure}[h!]
			\centering
		\subfloat[cLxF, $h$.]{%
			\includegraphics[width=0.25\textwidth,trim={2cm 0 1cm 0}]{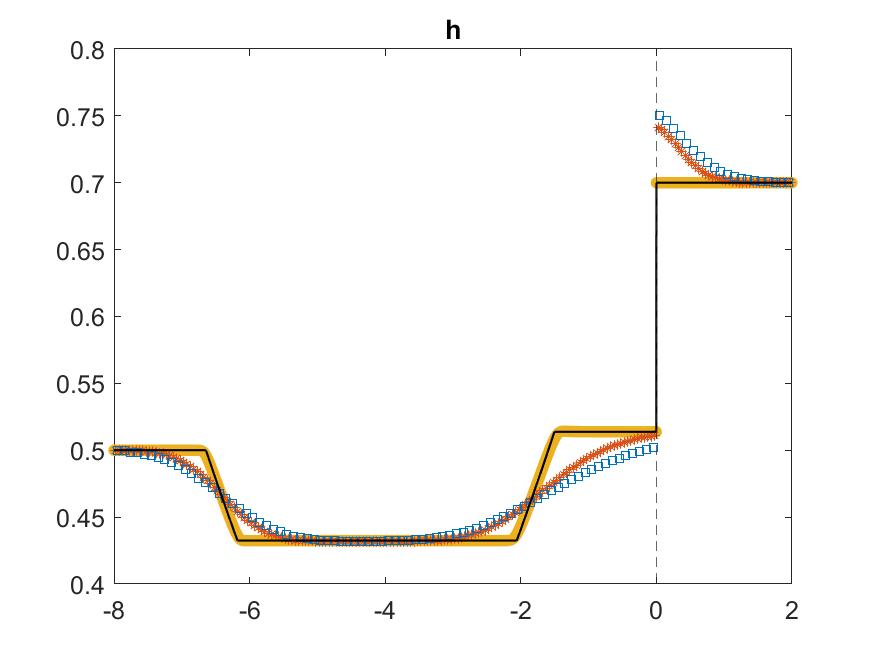}
		}
		\subfloat[cLxF, $m = hu$.]{%
			\includegraphics[width=0.25\textwidth,trim={2cm 0 1cm 0}]{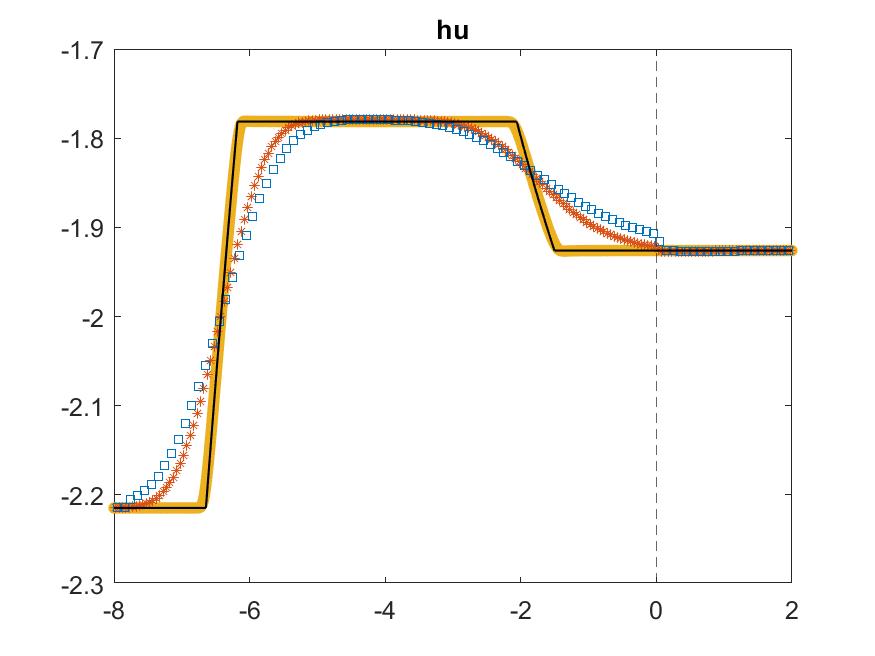}
		}
		\subfloat[LxF, $h$.]{%
			\includegraphics[width=0.25\textwidth,trim={2cm 0 1cm 0}]{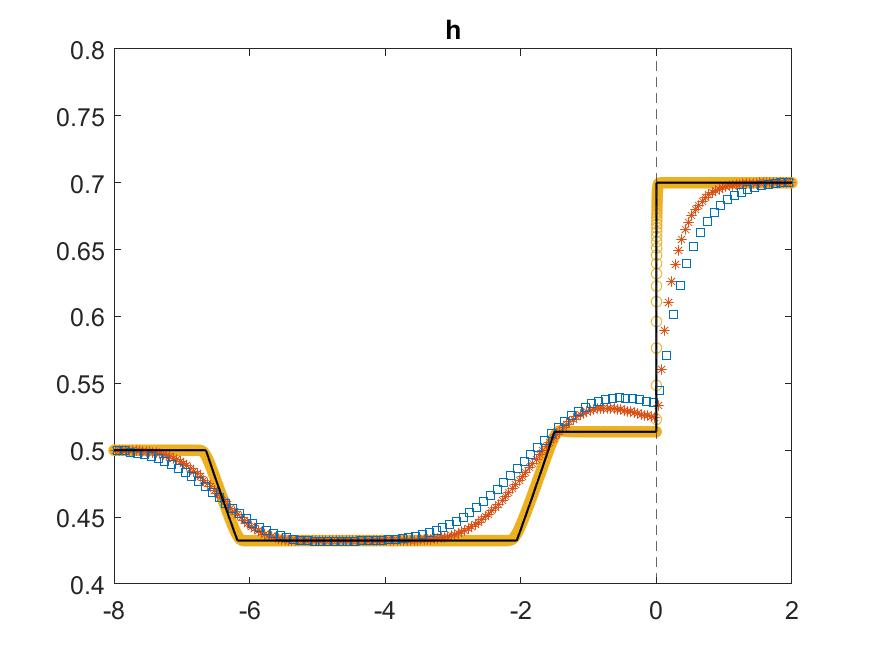}
		}
		\subfloat[LxF, $m = hu$.]{%
			\includegraphics[width=0.25\textwidth,trim={2cm 0 1cm 0}]{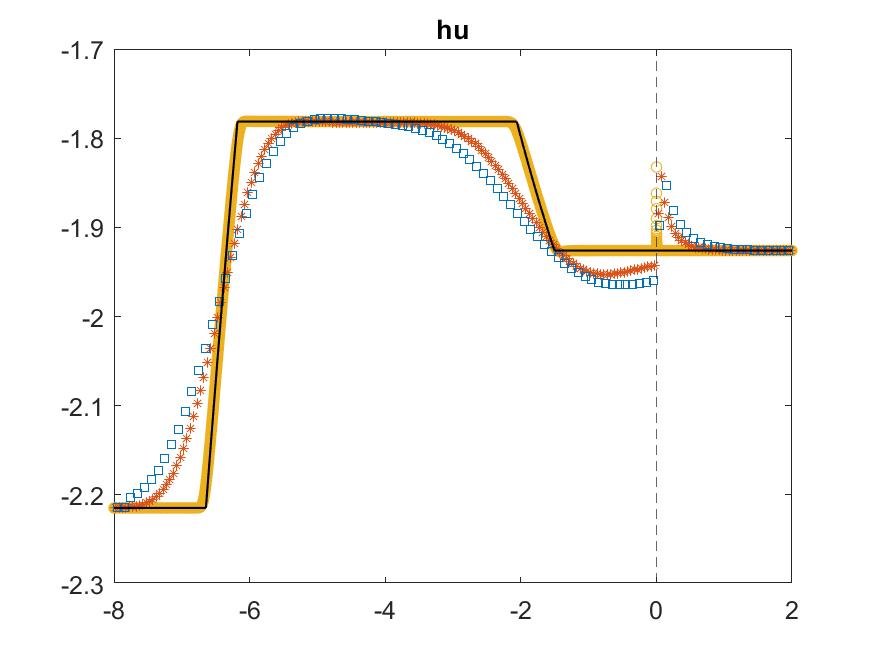}
		}
		\caption{Solutions to Example \ref{ex-1r2r-supern}. Black solid line: exact solution; blue squares: $N = 100$; red stars: $N = 200$, yellow circles: $N = 25600$.}\label{fig-1r2r-supern}
	\end{figure}
	\begin{table}[h!]
		\centering
		\begin{tabular}{c|c|c|c|c|c|c|c|c}
			\hline
			$N$&cLxF $e_{h}$&order&cLxF $e_{\mm}$&order&LxF $e_{h}$&order&LxF $e_{\mm}$&order\\
			\hline
			100&  1.18e-01&  -&  2.85e-01&  -&  1.77e-01&  -&  3.37e-01&  -\\ 
			200&  7.93e-02&  0.58&  1.96e-01&  0.54&  1.15e-01&  0.63&  2.25e-01&  0.58\\ 
			400&  4.86e-02&  0.70&  1.30e-01&  0.59&  7.19e-02&  0.68&  1.47e-01&  0.62\\ 
			800&  2.96e-02&  0.72&  8.62e-02&  0.60&  4.38e-02&  0.71&  9.44e-02&  0.63\\ 
			1600&  1.81e-02&  0.71&  5.61e-02&  0.62&  2.60e-02&  0.75&  5.99e-02&  0.66\\ 
			3200&  1.10e-02&  0.72&  3.57e-02&  0.65&  1.51e-02&  0.78&  3.72e-02&  0.69\\ 
			6400&  6.57e-03&  0.74&  2.21e-02&  0.69&  8.65e-03&  0.81&  2.26e-02&  0.72\\ 
			12800&  3.88e-03&  0.76&  1.33e-02&  0.73&  4.87e-03&  0.83&  1.35e-02&  0.75\\ 
			25600&  2.25e-03&  0.78&  7.82e-03&  0.77&  2.72e-03&  0.84&  7.88e-03&  0.77\\ 
			\hline
		\end{tabular}
	\caption{$L^1$ error table for Example \ref{ex-1r2r-supern}.}\label{tab-1r2r-supern}
	\end{table}
\end{ex}
\begin{ex}\label{ex-db-gamma0}
	In this example, we revisit the test problem in Section \ref{sec:artifact} to examine the choice of $\gamma$ on the convergence of the numerical solution. We have seen that the cLxF scheme with $\gamma = \sgnb$ converges to the correct solution of \eqref{eq:db}. However, if $\gamma = 0$ is used in the cLxF scheme, %in other words, we set the pressure to be $p(0,z) = \rho g(\{h+b\}-z)$ at the bottom discontinuity, 
in other words, we set $\check{h}_{\gamma, j+\hf} = \hbn{h+b}_{\gamma,j+\hf} - \{b\}_{j+\hf}= \{h\}_{j+\hf}$ in \eqref{eq:lxfsource}, 
then the cLxF scheme will converge to a different solution from that produced by the exact Riemann solver. See Figures \ref{fig-db-gamma0-1} and \ref{fig-db-gamma0-2}, as well as Table \ref{fig-db-gamma0}. This mismatch can be expected, since the derivation of the exact solution assumes $\gamma = \sgnb$ but the numerical solution in this test uses $\gamma = 0$. Furthermore, despite the wrong convergence, the solution limit of the cLxF scheme does not form a spike in the numerical momentum, while that of the LxF scheme does. This phenomenon can be predicted with Theorems \ref{thm:spike} and \ref{thm:cLxFm}.
	\begin{figure}[!ht]
		\centering
		\subfloat[cLxF ($\gamma = 0$), $h$.\label{fig-db-gamma0-1}]{%
			\includegraphics[width=0.25\textwidth,trim={2cm 0 1cm 0}]{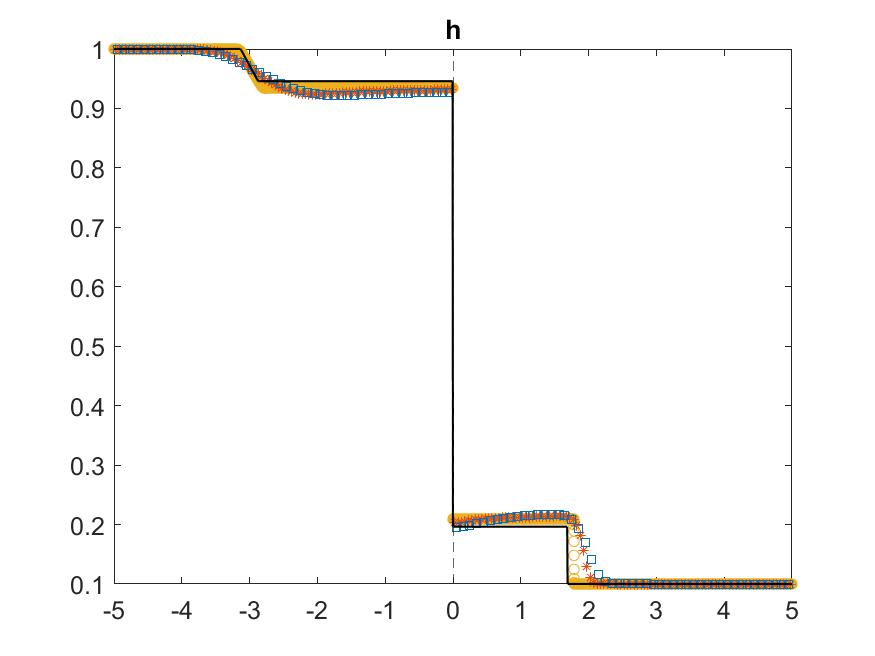}
		}
		\subfloat[cLxF ($\gamma = 0$), $m = hu$.\label{fig-db-gamma0-2}]{%
			\includegraphics[width=0.25\textwidth,trim={2cm 0 1cm 0}]{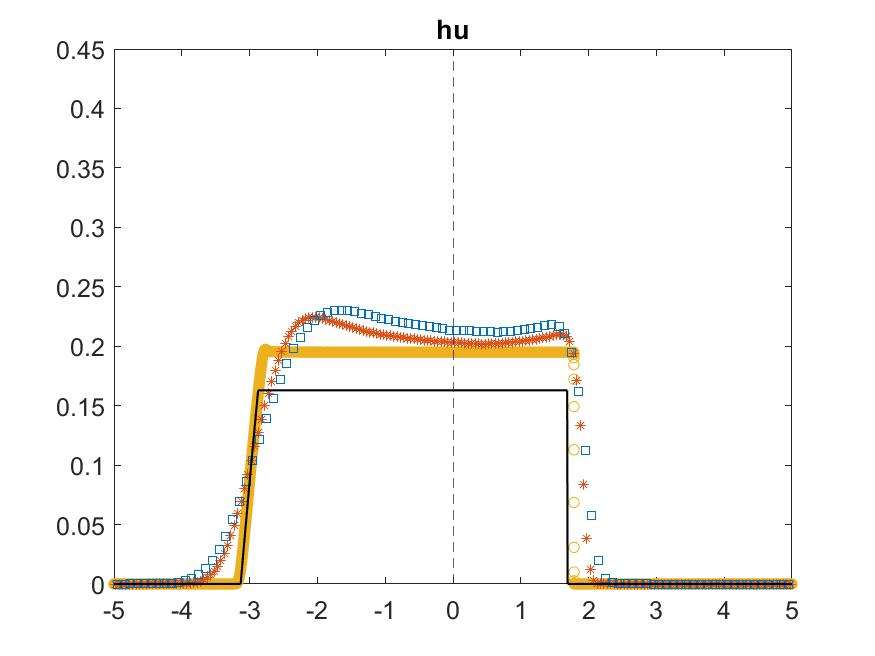}
		}
		\subfloat[LxF ($\gamma = 0$), $h$.]{%
			\includegraphics[width=0.25\textwidth,trim={2cm 0 1cm 0}]{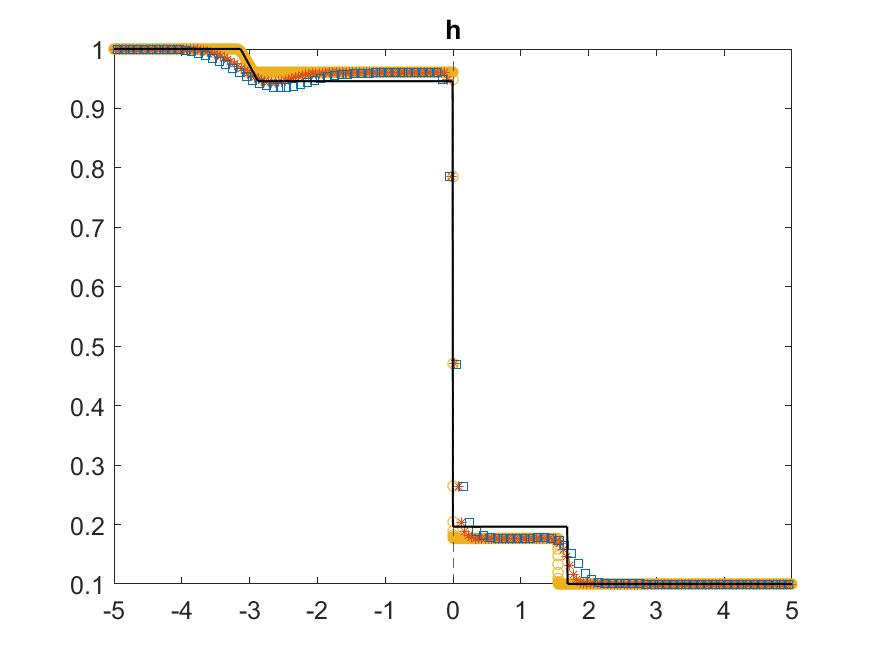}
		}
		\subfloat[LxF ($\gamma = 0$), $m = hu$.]{%
			\includegraphics[width=0.25\textwidth,trim={2cm 0 1cm 0}]{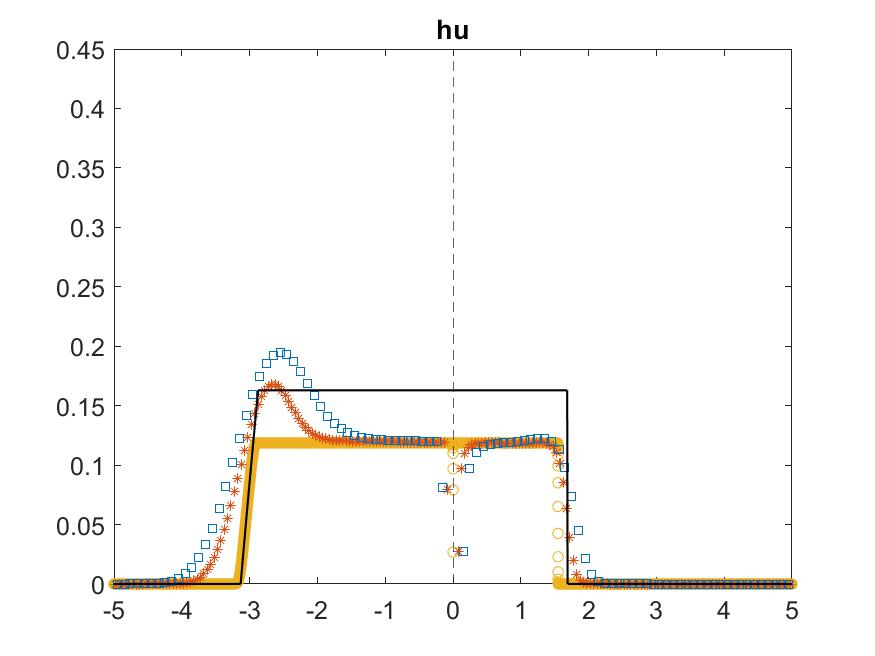}
		}
		\caption{Solutions to Example \ref{ex-db-gamma0}. $\gamma = 0$. Black solid line: exact solution; blue squares: $N = 100$; red stars: $N = 200$, yellow circles: $N = 25600$.}\label{fig-db-gamma0}
	\end{figure}
	\begin{table}[h!]
		\centering
		\begin{tabular}{c|c|c|c|c|c|c|c|c}
			\hline
			$N$&cLxF $e_{h}$&order&cLxF $e_{\mm}$&order&LxF $e_{h}$&order&LxF $e_{\mm}$&order\\
			\hline
			100&  1.18e-01&  -&  3.19e-01&  -&  1.40e-01&  -&  3.67e-01&  -\\ 
			200&  1.03e-01&  0.20&  2.66e-01&  0.27&  1.07e-01&  0.39&  2.81e-01&  0.38\\ 
			400&  8.77e-02&  0.23&  2.23e-01&  0.25&  9.67e-02&  0.14&  2.53e-01&  0.15\\ 
			800&  7.81e-02&  0.17&  1.98e-01&  0.17&  9.29e-02&  0.06&  2.39e-01&  0.08\\ 
			1600&  7.18e-02&  0.12&  1.82e-01&  0.12&  9.14e-02&  0.02&  2.32e-01&  0.04\\ 
			3200&  6.81e-02&  0.08&  1.73e-01&  0.07&  9.06e-02&  0.01&  2.28e-01&  0.02\\ 
			6400&  6.64e-02&  0.04&  1.69e-01&  0.04&  8.97e-02&  0.01&  2.25e-01&  0.02\\ 
			12800&  6.53e-02&  0.02&  1.66e-01&  0.02&  8.93e-02&  0.01&  2.23e-01&  0.01\\ 
			25600&  6.47e-02&  0.01&  1.65e-01&  0.01&  8.90e-02&  0.00&  2.22e-01&  0.01\\ 
			\hline
		\end{tabular}
		\caption{$L^1$ error table for Example \ref{ex-db-gamma0}. $\gamma = 0$.}\label{tab-db-gamma0}
	\end{table}	
\end{ex}

\begin{ex}\label{ex-coincide}
	This Riemann test is taken from \cite[Test SBRS]{rosatti2010riemann}. In this test, a rarefaction wave coincides with the stationary 0-wave and two shock waves are developed on both sides of the step. The critical flow condition $Fr= 1$ is imposed on the right side of the step. This is a resonant (transonic) test in the sense that the flow changes from the subcritical condition ($Fr<1$) to the supercritical condition ($Fr>1$) across the bottom step. The exact solution admits the following states:
	\begin{equation}
		W_L = \begin{pmatrix}
			4\\1.1175
		\end{pmatrix}\xrightarrow{\text{1S}}\begin{pmatrix}
			6.1431\\
			0.5089
		\end{pmatrix}\xrightarrow{0}%\text{0-wave}}
	\begin{pmatrix}
		3.9157\\
		1
	\end{pmatrix}\xrightarrow{R}%\text{0-wave}}
\begin{pmatrix}
	1.9999
	\\
	2.1977
\end{pmatrix}\xrightarrow{2S}%\text{0-wave}}
\begin{pmatrix}
1.0299\\
2.2428
\end{pmatrix} = W_R.
\end{equation}  

\begin{figure}[h!]
%	\centering
\subfloat[cLxF, $h$.]{%
\includegraphics[width=0.25\textwidth,trim={2cm 0 1cm 0}]{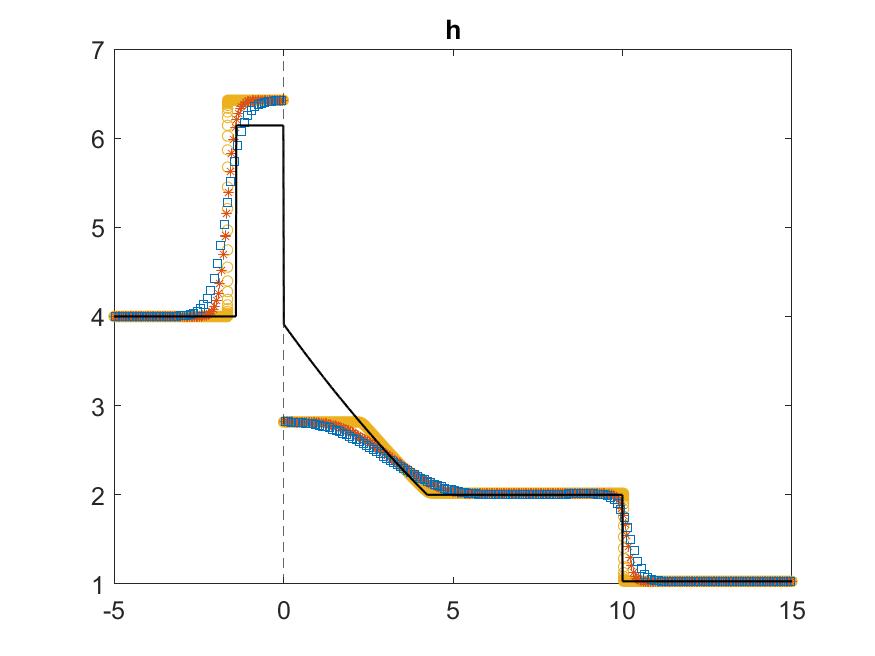}
}\label{fig:h_ex15scheme1a}
\subfloat[cLxF, $m = hu$.]{%
\includegraphics[width=0.25\textwidth,trim={2cm 0 1cm 0}]{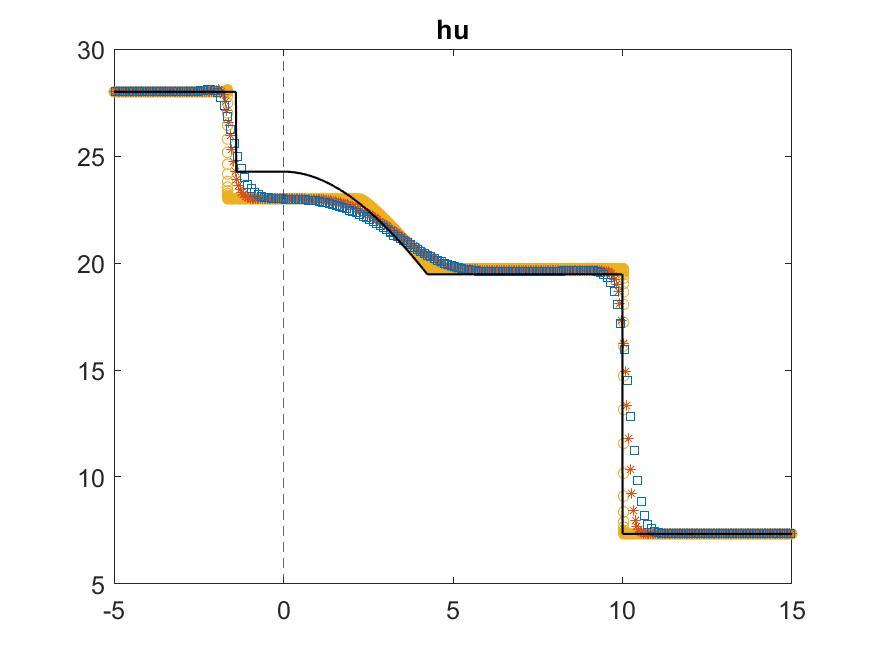}
}\label{fig:hu_ex15scheme1b}
\subfloat[LxF, $h$.]{%
\includegraphics[width=0.25\textwidth,trim={2cm 0 1cm 0}]{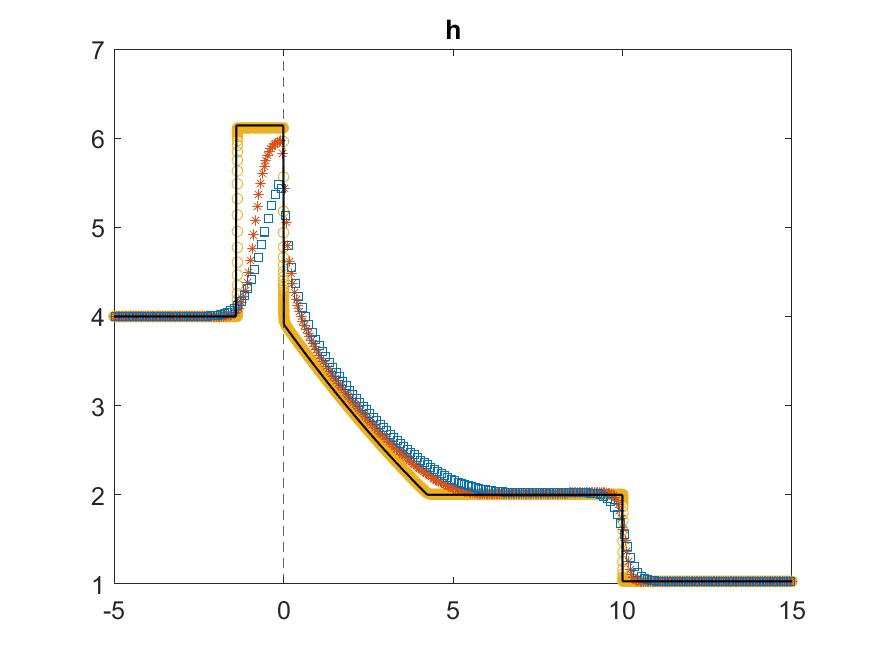}
}
\subfloat[LxF, $m = hu$.]{%
\includegraphics[width=0.25\textwidth,trim={2cm 0 1cm 0}]{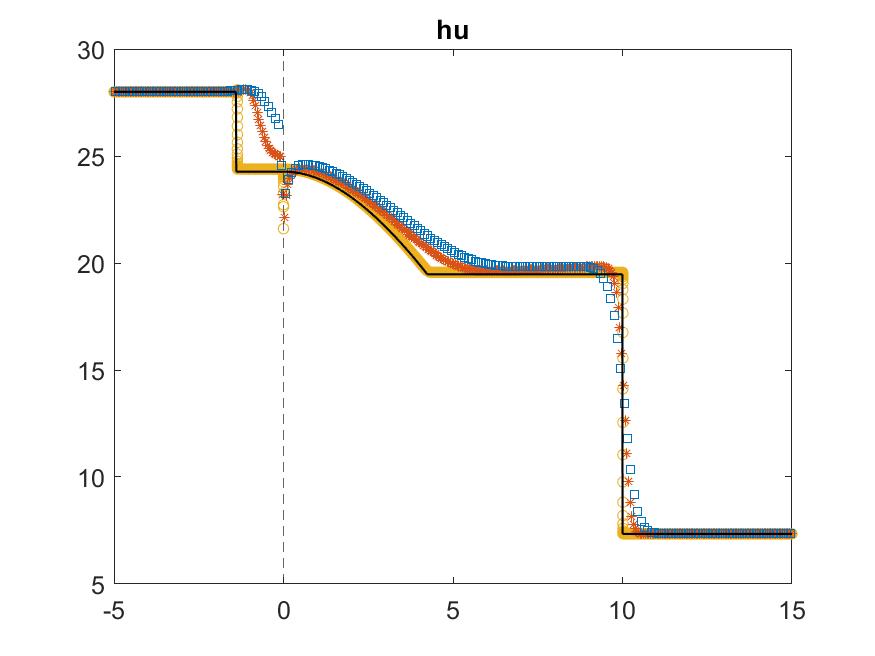}
}
\caption{Solutions to Example \ref{ex-coincide}. Black solid line: exact solution; blue squares: $N = 100$; red stars: $N = 200$, yellow circles: $N = 102400$.}\label{fig-resonant_a}
\end{figure}
\begin{table}[h!]
\centering
\begin{tabular}{c|c|c|c|c|c|c|c|c}
\hline
$N$&cLxF $e_{h}$&order&cLxF $e_{\mm}$&order&LxF $e_{h}$&order&LxF $e_{\mm}$&order\\
\hline
  200&  3.08e+00&  -&  9.64e+00&  -&  3.99e+00&  -&  1.40e+01&  -\\ 
  400&  2.76e+00&  0.16&  8.30e+00&  0.22&  2.71e+00&  0.55&  8.98e+00&  0.64\\ 
  800&  2.60e+00&  0.09&  7.72e+00&  0.10&  1.72e+00&  0.66&  5.48e+00&  0.71\\ 
 1600&  2.53e+00&  0.04&  7.49e+00&  0.04&  1.07e+00&  0.69&  3.54e+00&  0.63\\ 
 3200&  2.46e+00&  0.04&  7.37e+00&  0.02&  7.02e-01&  0.60&  2.56e+00&  0.47\\ 
 6400&  2.43e+00&  0.02&  7.35e+00&  0.00&  4.79e-01&  0.55&  2.02e+00&  0.34\\ 
12800&  2.42e+00&  0.01&  7.34e+00&  0.00&  3.46e-01&  0.47&  1.72e+00&  0.24\\ 
25600&  2.42e+00&  0.00&  7.35e+00& 0.00&  2.77e-01&  0.33&  1.58e+00&  0.12\\ 
51200&  2.41e+00&  0.00&  7.35e+00& 0.00&  2.36e-01&  0.23&  1.50e+00&  0.07\\ 
102400&  2.41e+00&  0.00&  7.35e+00& 0.00&  2.14e-01&  0.14&  1.46e+00&  0.04\\ 
\hline
\end{tabular}
\caption{$L^1$ error table for Example \ref{ex-coincide}.}\label{tab:resonant_a}
\end{table}
From Figure \ref{fig-resonant_a}, we can see that the cLxF scheme converges to a wrong solution which fails to capture the correct wave pattern. The solution limit of the cLxF scheme forms two constant states on both sides of the bottom discontinuity, and a right rarefaction wave is developed away from the bottom step. In contrast, in the exact solution, the right rarefaction wave should be adjacent to the bottom step. This wrong convergence may be attributed to the fact that the cLxF scheme does not have enough numerical viscosity to drive the solution converging towards the physical solution. Similar numerical difficulty also occurs in conservative schemes for hyperbolic conservation laws -- for example, without an entropy fix, the Roe's scheme may converge to a weak solution violating the entropy condition for (tran)sonic rarefaction waves \cite{leveque2002finite}. For the LxF scheme, although at a glance that the scheme captures the correct solution profile, we notice that the convergence rate drops to somewhere close to 0 after the mesh refinement, which indicates that the LxF scheme indeed also converges to a wrong solution. At the same time, we also observe the numerical artifact of developing a spurious spike in the numerical momentum at the bottom discontinuity. The behavior of the LxF scheme for this test problem is similar to those for the subcritical tests. 
\end{ex}

\begin{ex}\label{ex-resonant}
		This Riemann problem is taken from \cite[Test RRBR]{rosatti2010riemann}. In this test, a rarefaction wave coincides with the stationary 0-wave and two rarefaction waves are developed on the left of the bottom step. As that in Example \ref{ex-coincide}, the critical flow condition $Fr=1$ is imposed on the right side of the step. This is also a resonant (transonic) test. The exact solution of the problem has the following states
	\begin{equation}
	W_L = \begin{pmatrix}
	6\\
	-2.0855
	\end{pmatrix}\xrightarrow{\text{1R}}\begin{pmatrix}
	1.9766\\
	-2.1490
	\end{pmatrix}\xrightarrow{\text{2R}}%\text{0-wave}}
	\begin{pmatrix}
	2.5253\\
	-1.6707 
	\end{pmatrix}\xrightarrow{0}%\text{0-wave}}
	\begin{pmatrix}
	3.5556\\
	-1
	\end{pmatrix}\xrightarrow{\text{R}}
	\begin{pmatrix}
	8\\0
	\end{pmatrix} = W_R.
	\end{equation}
	
	From Figure \ref{fig-resonant}, it can be seen that the cLxF scheme converges to a wrong solution. We expect this is again due to the incapability to capture the entropy solution with the central flux at the bottom discontinuity. While the LxF scheme seems to achieve the correct $L^1$ convergence, despite some minor numerical artifacts at the bottom discontinuity. See Figure \ref{fig-resonant} and Table \ref{tab:resonant}. The LxF scheme in this test behaves similarly to those for the negative supercritical tests.
\begin{figure}[h!]
%	\centering
	\subfloat[cLxF, $h$.]{%
		\includegraphics[width=0.25\textwidth,trim={2cm 0 1cm 0}]{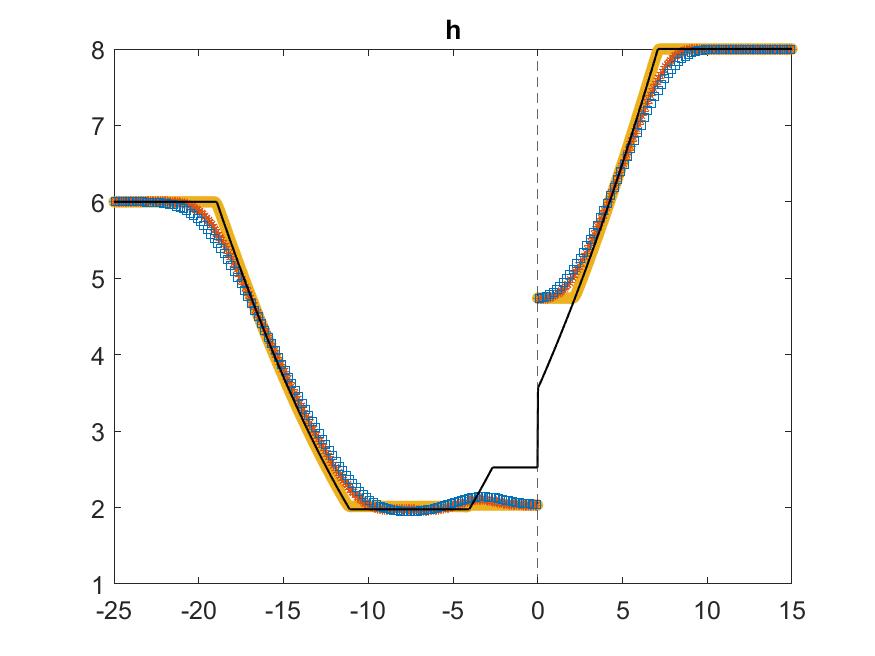}
	}\label{fig:h_ex15scheme1c}
	\subfloat[cLxF, $m = hu$.]{%
		\includegraphics[width=0.25\textwidth,trim={2cm 0 1cm 0}]{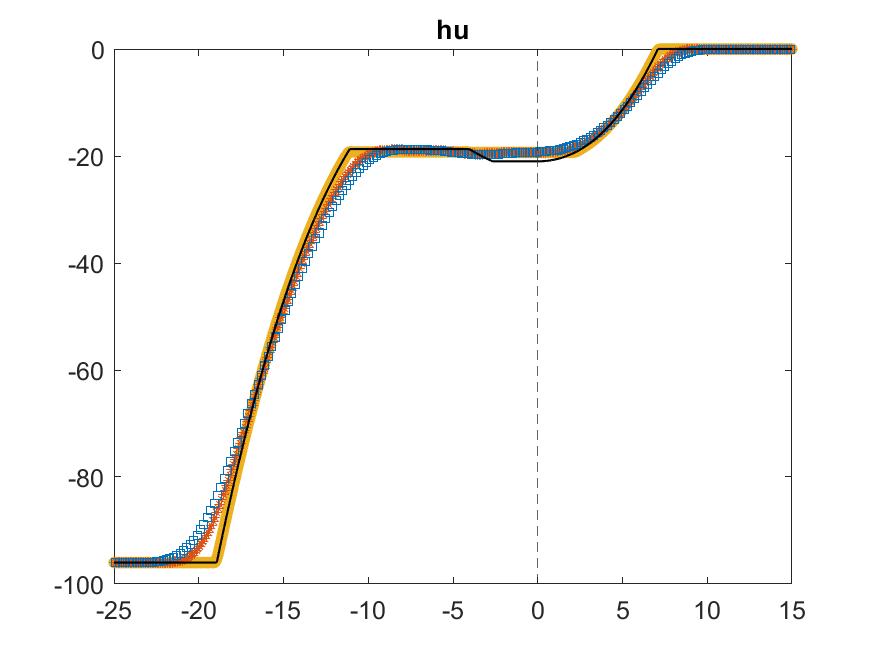}
	}\label{fig:hu_ex15scheme1d}
	\subfloat[LxF, $h$.]{%
		\includegraphics[width=0.25\textwidth,trim={2cm 0 1cm 0}]{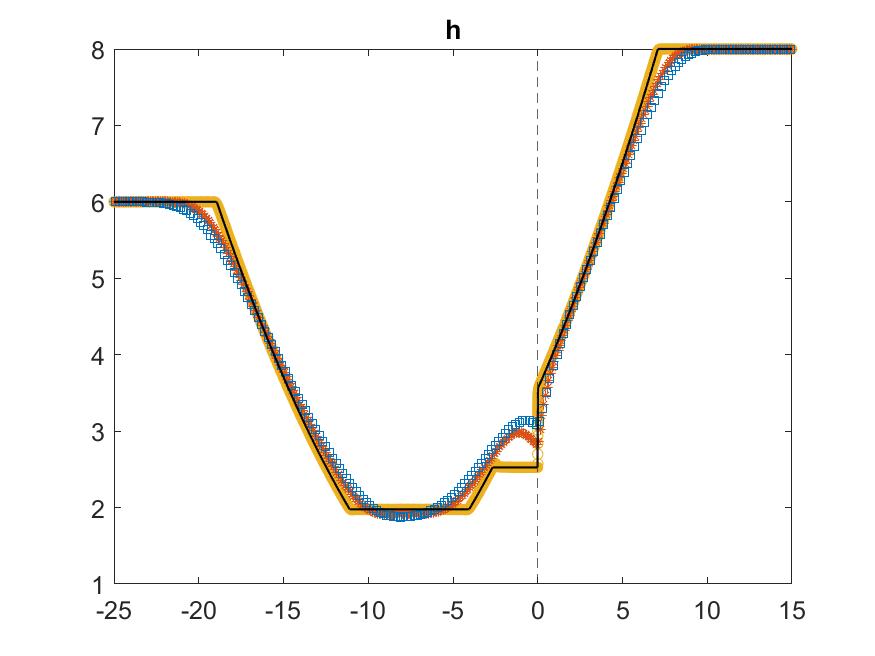}
	}
	\subfloat[LxF, $m = hu$.]{%
		\includegraphics[width=0.25\textwidth,trim={2cm 0 1cm 0}]{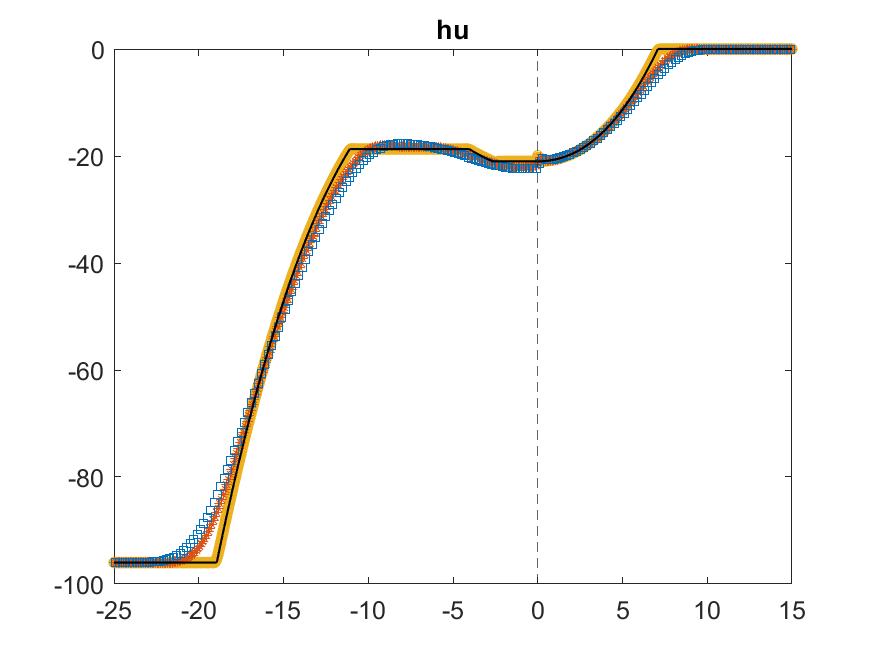}
	}
	\caption{Solutions to Example \ref{ex-resonant}. Black solid line: exact solution; blue squares: $N = 100$; red stars: $N = 200$, yellow circles: $N = 102400$.}\label{fig-resonant}
\end{figure}
% \begin{figure}[h!]
% 		\centering
% 	\subfloat[cLxF, $h$.]{%
% 		\includegraphics[width=0.25\textwidth,trim={2cm 0 0cm 0}]{h_exRRBR_scheme1_zoomin.jpg}
% 	}\label{fig:h_ex15scheme1_zoomedin}
% 	\subfloat[LxF, $h$.]{%
% 		\includegraphics[width=0.25\textwidth,trim={2cm 0 0cm 0}]{h_exRRBR_scheme2_zoomin.jpg}
% 	}
% 	\caption{Zoomed-in views of solutions to Example \ref{ex-resonant}. Black solid line: exact solution; blue squares: $N = 100$; red stars: $N = 200$, yellow circles: $N = 102400$.}\label{fig-resonant_zoomin}
% \end{figure}
\begin{table}[h!]
	\centering
	\begin{tabular}{c|c|c|c|c|c|c|c|c}
		\hline
		$N$&cLxF $e_{h}$&order&cLxF $e_{\mm}$&order&LxF $e_{h}$&order&LxF $e_{\mm}$&order\\
		\hline
		  200&  7.39e+00&  -&  7.20e+01&  -&  6.67e+00&  -&  6.76e+01&  -\\ 
  400&  5.79e+00&  0.35&  4.88e+01&  0.56&  4.46e+00&  0.58&  4.27e+01&  0.66\\ 
  800&  4.77e+00&  0.28&  3.41e+01&  0.52&  2.85e+00&  0.64&  2.63e+01&  0.70\\ 
 1600&  4.14e+00&  0.20&  2.50e+01&  0.45&  1.76e+00&  0.69&  1.59e+01&  0.72\\ 
 3200&  3.76e+00&  0.14&  1.95e+01&  0.36&  1.06e+00&  0.73&  9.46e+00&  0.75\\ 
 6400&  3.53e+00&  0.09&  1.62e+01&  0.27&  6.23e-01&  0.77&  5.51e+00&  0.78\\ 
12800&  3.40e+00&  0.06&  1.42e+01&  0.18&  3.58e-01&  0.80&  3.15e+00&  0.81\\ 
25600&  3.32e+00&  0.03&  1.31e+01&  0.12&  2.03e-01&  0.82&  1.78e+00&  0.83\\ 
51200&  3.27e+00&  0.02&  1.25e+01&  0.07&  1.13e-01&  0.84&  9.90e-01&  0.84\\ 
102400&  3.25e+00&  0.01&  1.22e+01&  0.04&  6.27e-02&  0.86&  5.46e-01&  0.86\\ 
		\hline
	\end{tabular}
		\caption{$L^1$ error table for Example \ref{ex-resonant}.}\label{tab:resonant}
\end{table}

\end{ex}

\section{Conclusions}\label{sec:conclusion}
This paper concerns a numerical artifact that occurred in solving the nonlinear SWEs over a discontinuous riverbed. Several first-order schemes, including the well-balanced and non-well-balanced LxF schemes, the HR scheme, and the XS scheme, are considered. We observe that the solutions of these schemes would form a spurious spike in the numerical momentum, which may prevent the numerical solution from converging to the exact weak solution of the SWEs. To explain the cause of this numerical artifact, we rewrite the above-mentioned first-order schemes into a unified form and establish a Lax--Wendroff type theorem to study the convergence of the methods. Based on the Lax--Wendroff type theorem, we are able to show that the spurious spike in the momentum is caused by the numerical viscosity in the equation of mass conservation at the bottom discontinuity and the height of the spike is proportional to the viscosity constant of the LxF flux. Furthermore, we note that by setting the LxF constant to be zero at the bottom discontinuity in the equation of mass conservation, or in other words, by adopting the central flux for $\wh{m}$, one can avoid the formation of the one-sided spurious spike. The resulting method is referred to as the cLxF scheme. Numerically we validate that the cLxF scheme has the correct convergence for nontransonic problems, although it may fail for the transonic tests. In our future works, we will investigate corrective procedures for retrieving correct convergence for cLxF schemes for transonic problems.

\appendix

\section{Derivation of \eqref{eq-mom-cont}}\label{app:weaksoln} %\YX{\cite{ROSATTI200810058,bernetti2008exact}} \YX{\cite{bernetti2008exact} seems to have some assumptions, which needs to be included. Or we can just refer to \cite{bernetti2008exact} without discussing the details.}  \zs{This is not exactly the same as \cite{bernetti2008exact}. I am not very confident on how assumptions should be made rigorously... I call it a ``heuristic" derivation to avoid discussing the technical assumptions... }

Use the fact $p(x,z) \equiv 0$ at the air-water interface. One can explicitly write down $\Omega(t)$ and $\partial \Omega(t)$ in \eqref{eq-mom} to obtain
\begin{equation}\label{eq-ctrl}
\begin{aligned}
\frac{\dd}{\dd t}\int_{x_L(t)}^{x_R(t)}\int_{b(x)}^{(h+b)(x)}\rho u \dd z \dd x + \int_{b(x_L)}^{(h+b)(x_L)} (p \nux)(x_L,z) \dd z + \int_{b(x_R)}^{(h+b)(x_R)} (p\nux)(x_R,z)\dd z \\
+\int_{x_L}^{x_R} (p\nux)(x,b(x)) \dd x + \sgnb \int_{b^-}^{b^+} (p\nux)(0,z) \dd z = 0.
\end{aligned}
\end{equation}
Note that $\nux$ admits to the following values
\begin{equation}
\nux = \left\{
\begin{array}{cl}
1,& (x,z) = (x_L,z),\\
-1,& (x,z) = (x_R,z),\\
b_x, & (x,z) = (x, b(x)),\\
\sgnb,& (x,z) = (x_\star,z). 
\end{array}
\right.
\end{equation}
Under the assumption of the hydrostatic pressure distribution, we have  
\begin{equation}\label{eq-pre}
p(x,z) = \rho g (h(x)+b(x)-z), \quad \text{if } x\neq 0.
\end{equation}
We substitute the pressure into \eqref{eq-ctrl}. After integration with respect to $z$ and dividing by $\rho$ on both sides of the equation, it gives
\begin{equation}\label{eq-mom1}
\frac{\dd}{\dd t}\int_{x_L(t)}^{x_R(t)} h u  \dd x - \hf gh^2\bigg|_{x = x_L} +\hf g h^2\bigg|_{x = x_R} + \int_{x_L}^{x_R}g h b_x \dd x + \int_{b^-}^{b^+} \frac{p(0,z)}{\rho} \dd z = 0.
\end{equation}
Note that 
\begin{equation}\label{eq-mom2}
\begin{aligned}
\frac{\dd}{\dd t}\int_{x_L(t)}^{x_R(t)} h u  \dd x =\;& \int_{x_L}^{x_R} (h u)_t  \dd x + \left(\frac{\dd x_R}{\dd t} hu\right)\bigg|_{x=x_R} -  \left(\frac{\dd x_L}{\dd t} hu\right)\bigg|_{x=x_L} \\
=\;& \int_{x_L}^{x_R} (h u)_t  \dd x + hu^2\bigg|_{x=x_R} -   hu^2\bigg|_{x=x_L}.
\end{aligned}
\end{equation}
We can substitute \eqref{eq-mom2} into \eqref{eq-mom1} and rearrange terms to obtain \eqref{eq-mom-cont}.

%\YX{also in \cite[Eqs. (17)-(18)]{ROSATTI200810058}, \cite[Appendix C]{MURILLO20104327}. Double check for the validity of the definition of weak solution, the literature seems to be for step function $b$ only???}\zs{Yes...}

\section{Proof of Theorem \ref{thm-HR}}\label{app-HR}

\begin{proof}
	Note that 
	\begin{equation}
		\wh{{F}}^* = 
		\wh{{F}} + \left(\widehat{{F}}^*- \wh{{F}}\right) = \wh{{F}} + \left\{{F}(U^*)-{F}(U)\right\} -\frac{1}{2}A\left[U^*-U\right].
	\end{equation}
	We denote by $\delta = h^* - h = b-\max(b^-,b^+)$. Then 
	\begin{equation}
		U^* - U = \left(
		\begin{array}{c}
			\delta\\
			0
		\end{array}\right),
	\end{equation}
	and 
	\begin{equation}
		{F}(U^*)-{F}(U) = \left(
		\begin{array}{c}
			\mm\\
			\frac{\mm^2}{h^*} + \frac{g}{2} (h^*)^2
		\end{array}\right) - 
		\left(
		\begin{array}{c}
			\mm\\
			\frac{\mm^2}{h} + \frac{g}{2} h^2
		\end{array}\right) = 
		\left(
		\begin{array}{c}
			0\\
			\frac{g}{2}\left(\left(h^*\right)^2-h^2\right)-\frac{\mm^2\delta}{h h^*} 
		\end{array}\right).	
	\end{equation}
	Hence
	\begin{equation}\label{eq:FstarF}
		\wh{F}^* = \wh{F} + 
		\left(
		\begin{array}{c}
			-\frac{\alpha_1}{2}[\delta]\\
			\frac{g}{2}\left\{(h^*)^2-h^2\right\}-\left\{\frac{\mm^2 \delta}{hh^*}\right\}
		\end{array}\right).
	\end{equation}
	Therefore, after substituting \eqref{eq:FstarF} into \eqref{eq:hrFhat}, one can get
	\begin{equation}\label{eq-hrm}
		\begin{aligned}
			\wh{F}^{*,\pm} =&\, \wh{F} +  \left(
			\begin{array}{c}
				-\frac{\alpha_1}{2}[\delta]\\
				\frac{g}{2}\left\{(h^*)^2-h^2\right\}-\left\{\frac{\mm^2\delta}{hh^*}\right\}
			\end{array}\right)
			+ \left(
			\begin{array}{c}
				0\\
				\frac{g}{2}\left(\left(h^\pm\right)^2 - (h^{*,\pm})^2\right)
			\end{array}\right)\\
			=&\, \wh{{F}} - \left(
			\begin{array}{c}
				\frac{\alpha_1}{2}[\delta]\\
				\pm \frac{g}{4}\left[\left(h^*\right)^2-h^2\right]+\left\{\frac{\mm^2\delta}{h h^*}\right\}
			\end{array}\right).
		\end{aligned}
	\end{equation}
	Let us introduce the notation 
	\begin{equation}\label{eq:sgnavg}
		\hbn{v} = \hbn{v}_\sgnb = 
		\left\{\begin{array}{cc}
			v^- & [b]>0\\
			v^+ & [b]<0
		\end{array}\right.. 
	\end{equation}
	Note we have
	\begin{equation}\label{eq:deltabjump}
		\delta^+ = 
		\left\{
		\begin{matrix}
			0,&[b]>0\\
			[b],&[b]<0
		\end{matrix}\right.\quand 		
		\delta^- = \left\{
		\begin{matrix}
			-[b],&[b]>0\\
			0,&[b]<0
		\end{matrix}\right..
\end{equation}
	Hence it can be verified that $[\delta v] = [b]\hbn{v}$ and $\hbn{\delta} = 2\left(\hbn{b}-\{b\}\right)$. As a result, it gives
	\begin{equation}\label{eq-hstar2}
		[(h^*)^2-h^2] =	[(h+\delta)^2-h^2] =  [\delta ( 2h+\delta)] = [b]\hbn{2h+\delta} = 2[b]\left(\hbn{h + b} -\{b\}\right). 
	\end{equation}
	In addition,
	with \eqref{eq:sgnavg} and \eqref{eq:deltabjump}, we can show that 
	\begin{equation}\label{eq:hr-redum}
		\left\{\frac{\mm^2\delta}{hh^*}\right\} = \left\{
		\begin{matrix}
			\left(\frac{\mm^2}{2h}\right)^-\frac{-[b]}{h^--[b]},&[b]>0\\
			\left(\frac{\mm^2}{2h}\right)^+\frac{[b]}{h^++[b]},&[b]<0\\		
		\end{matrix}\right. = \hbn{\frac{\mm^2}{2h}}\frac{|[b]|}{|[b]|-\hbn{h}}.
	\end{equation}
	Substituting \eqref{eq-hstar2} and \eqref{eq:hr-redum} into \eqref{eq-hrm}, together with the fact $[\delta] = [b]$, we can obtain 
	\begin{equation}\label{eq:hrflux}
		\begin{aligned}
			\wh{F}^{*,\pm} =& \wh{F} -\begin{pmatrix}
				\frac{\alpha_1}{2}[b]\\
				\pm \frac{g}{2}\left(\hbn{h + b} -\{b\}\right)[b] + \hbn{\frac{\mm^2}{2h}}\frac{|[b]|}{|[b]|-\hbn{h}}
			\end{pmatrix}\\
			=& \widehat{F} - \left(\pm\frac{g}{2}\begin{pmatrix}
				0\\\hbn{h + b} -\{b\}
			\end{pmatrix}[b] + \hf\begin{pmatrix}
				\alpha_1\\
				\hbn{\frac{\mm^2}{h}}\frac{\sgnb}{|[b]|-\hbn{h}}
			\end{pmatrix}[b]\right).
		\end{aligned}
	\end{equation}
	Substituting \eqref{eq:hrflux} into \eqref{eq:hr}, we get (omitting the superscript $n$ on the right)
	\begin{equation}
		\begin{aligned}
			&\frac{U_j^{n+1}-U_j^n}{\Delta t^n} + \frac{\wh{F}_{j+\hf}^n-\wh{F}_{j-\hf}^n}{\dx_j}\\
			= &\frac{1}{\Delta x_j}\left({-\frac{g}{2}\begin{pmatrix}
					0\\\hbn{h + b} -\{b\}
				\end{pmatrix}_{j+\hf}[b]_{j+\hf}-\frac{g}{2}\begin{pmatrix}
					0\\\hbn{h + b} -\{b\}
				\end{pmatrix}_{j-\hf}[b]_{j-\hf}}\right)\\
			&+ \frac{1}{\Delta x_j}\left({\hf\begin{pmatrix}
					\alpha_1\\
					\hbn{\frac{\mm^2}{h}}\frac{\sgnb}{|[b]|-\hbn{h}}
				\end{pmatrix}_{j+\hf}[b]_{j+\hf}- \hf\begin{pmatrix}
					\alpha_1\\
					\hbn{\frac{\mm^2}{h}}\frac{\sgnb}{|[b]|-\hbn{h}}
				\end{pmatrix}_{j-\hf}[b]_{j-\hf}}\right).
		\end{aligned}
	\end{equation}
	The values of $\gamma$ and $\widehat{N}$ can be read from the above reformulated numerical scheme. 
\end{proof}
\section{Proof of Theorem \ref{thm-XS}}\label{app-XS}

\begin{proof}
	Noting that 
	\begin{equation}
		\left\{v\right\}_{j+\hf} - \left\{v\right\}_{j-\hf} = \hf\left(v_{j+1}-v_{j-1}\right) 
		= \hf\left(v_{j+1}-v_j\right) + \hf \left(v_j - v_{j-1}\right) = \hf\left[v\right]_{j+\hf} + \hf\left[v\right]_{j-\hf},
	\end{equation}
	and the identity $\left[b^2\right]/2 = \left\{b\right\}\left[b\right]$, the source term approximation \eqref{eq:fsrS} can be rewritten as
%	\begin{equation}
%		s_j = \frac{g}{2}\left(\hf\left[b^2\right]_{j+\hf} + \hf\left[b^2\right]_{j-\hf}\right) - \frac{g}{2}\left(h_j+b_j\right)\left(\left[b\right]_{j+\hf} + \left[b\right]_{j-\hf}\right).
%	\end{equation}
%	Using the identity $\left[b^2\right]/2 = \left\{b\right\}\left[b\right]$, it yields
	\begin{equation}
		\begin{aligned}
			s_j =& \frac{g}{2}\left(\hf\left[b^2\right]_{j+\hf} + \hf\left[b^2\right]_{j-\hf}\right) - \frac{g}{2}\left(h_j+b_j\right)\left(\left[b\right]_{j+\hf} + \left[b\right]_{j-\hf}\right) \\
			       =& \frac{g}{2}\left(\left\{b\right\}_{j+\hf}\left[b\right]_{j+\hf} + \left\{b\right\}_{j-\hf}\left[b\right]_{j-\hf}\right) - \frac{g}{2}\left(h_j+b_j\right)\left(\left[b\right]_{j+\hf} + \left[b\right]_{j-\hf}\right)\\
			=&-\frac{g}{2}\left(h_j+b_j-\{b\}_{j+\hf}\right)\left[b\right]_{j+\hf} -\frac{g}{2}\left(h_j+b_j-\left\{b\right\}_{j-\hf}\right)\left[b\right]_{j-\hf}.
		\end{aligned}
	\end{equation}
	With $v_j = \left\{v\right\}_{j+\hf} - \hf[v]_{j+\hf} = \{v\}_{j-\hf} + \hf[v]_{j-\hf}$, one can get
	\begin{equation}
		\begin{aligned}
			s_j
			=& -\frac{g}{2}\left(\left\{h+b\right\}_{j+\hf} -\hf[h+b]_{j+\hf}-\{b\}_{j+\hf}\right)[b]_{j+\hf}\\& -\frac{g}{2}\left(\{h+b\}_{j-\hf} + \hf[h+b]_{j-\hf}-\{b\}_{j-\hf}\right)[b]_{j-\hf}\\
			=& -\frac{g}{2}\left(\{h+b\}-\{b\}\right)_{j+\hf}[b]_{j+\hf} -\frac{g}{2}\left(\{h+b\} -\{b\}\right)_{j-\hf}[b]_{j-\hf} \\
			&+ \frac{g}{4}[h+b]_{j+\hf}[b]_{j+\hf} - \frac{g}{4}[h+b]_{j-\hf}[b]_{j-\hf}.
		\end{aligned}
	\end{equation}
	As a result, one can rewrite \eqref{eq:fsr} as  %\YX{Does the first term on RHS need to be $\hbn{h + b} -\{b\}$?}
	\begin{equation}
		\begin{aligned}
			&\frac{U_j^{n+1}-U_j^n}{\Delta t^n} + \frac{\wh{F}_{j+\hf}^n - \wh{F}_{j-\hf}^n}{\dx_j}\\
			=& \frac{1}{\dx_j}\left(-\frac{g}{2}\begin{pmatrix}
				0\\
				\{h+b\}-\{b\}
			\end{pmatrix}_{j+\hf}[b]_{j+\hf} 
			-\frac{g}{2}\begin{pmatrix}
				0\\
				\{h+b\} -\{b\}
			\end{pmatrix}_{j-\hf}[b]_{j-\hf}\right)\\
			&+ \frac{1}{\dx_j}\left(
			\hf\begin{pmatrix}
				{\alpha_1}\\
				\frac{g}{2}[h+b]
			\end{pmatrix}_{j+\hf}[b]_{j+\hf} - 
			\hf \begin{pmatrix}
				{\alpha_1}\\
				\frac{g}{2}[h+b]
			\end{pmatrix}_{j-\hf}[b]_{j-\hf}
			\right).
		\end{aligned}
	\end{equation}
	The values of $\gamma=0$ and $\widehat{N}$ can be read from the above reformulated numerical scheme. 
\end{proof}

\bibliography{refs}
\bibliographystyle{abbrv}

\end{document}